\documentclass[11pt]{amsart}

\usepackage{pgf,tikz,hyperref,comment,amsthm}
\usepackage{amsmath, graphicx}
\usepackage{amssymb,color}
\usepackage{enumitem}
\newcommand*{\DashedArrow}[1][]{\mathbin{\tikz [baseline=-0.25ex,-latex, dashed,#1] \draw [#1] (0pt,0.5ex) -- (1.3em,0.5ex);}}%

\newtheoremstyle{thmstyle3}
{18pt plus2pt minus1pt}
{18pt plus2pt minus1pt}
{\small\normalfont}
{0pt}
{\small\bfseries}
{}
{.5em}
{\thmname{#1}\thmnumber{ {#2}}%
	\thmnote{ {\the\thm@notefont(#3)}}}
\newtheoremstyle{thmstyle2}
{18pt plus2pt minus1pt}
{18pt plus2pt minus1pt}
{\normalfont}
{0pt}
{\bfseries}
{}
{.5em}
{\thmname{#1}\thmnumber{ {#2}.}%
	\thmnote{ {\the\thm@notefont(#3)}}}

\theoremstyle{thmstyleone}%
\newtheorem{theorem}{Theorem}[subsection]%
\newtheorem{proposition}[theorem]{Proposition}%
\newtheorem{lemma}[theorem]{Lemma}%
\newtheorem{corollary}[theorem]{Corollary}%
\theoremstyle{thmstyle2}%
\newtheorem{question}[theorem]{Question}%
\newtheorem{definition}[theorem]{Definition}%
\newtheorem{example}[theorem]{Example}%
\newtheorem{remark}[theorem]{Remark}%

\definecolor{mygreen}{rgb}{0.0, 0.5, 0.0}

\newcommand{\PP}{ \ensuremath{\mathbb{P}}}
\newcommand{\field}{K}

\newcommand{\call}{ \ensuremath{\mathcal{L}}}

\begin{document}

\title{Combinatorics of skew lines in $\mathbb P^3$ with an application to algebraic geometry}

\author[\tiny Chiantini]{Luca Chiantini}
\address[L.~Chiantini]{Dipartmento di Ingegneria dell'Informazione e Scienze Matematiche, Universit\`a di Siena, Italy}
\email{luca.chiantini@unisi.it}

\author[Farnik]{{\L}ucja Farnik}
\address[{\L}.~Farnik]{Department of Mathematics, University of the National Education Commission, Krakow,
   Podcho\-r\c a\.zych~2,
   PL-30-084 Krak\'ow, Poland}
\email{lucja.farnik@gmail.com}

\author[Favacchio]{Giuseppe Favacchio}
\address[G.~Favacchio]{Dipartimento di Ingegneria, Universit\`a degli studi di Palermo,
Viale delle Scienze,  90128 Palermo, Italy}
\email{giuseppe.favacchio@unipa.it}

\author[Harbourne]{Brian Harbourne}
\address[B.~Harbourne]{Department of Mathematics,
University of Nebraska,
Lincoln, NE 68588-0130 USA}
\email{brianharbourne@unl.edu}

\author[Migliore]{Juan Migliore} 
\address[J.~Migliore]{Department of Mathematics,
University of Notre Dame,
Notre Dame, IN 46556 USA}
\email{migliore.1@nd.edu}

\author[Szemberg]{Tomasz Szemberg}
\address[T.~Szemberg]{Department of Mathematics, University of the National Education Commission, Krakow,
   Podcho\-r\c a\.zych~2,
   PL-30-084 Krak\'ow, Poland}
\email{tomasz.szemberg@gmail.com}

\author[Szpond]{Justyna Szpond}
\address[J.~Szpond]{Department of Mathematics, University of the National Education Commission, Krakow,
   Podcho\-r\c a\.zych~2,
   PL-30-084 Krak\'ow, Poland}
\email{szpond@gmail.com}


\begin{abstract}
This article introduces a previously unrecognized combinatorial structure underlying configurations of skew lines in $\PP^3$, and reveals its deep and surprising connection to the algebro-geometric concept of geproci sets.

Given any field $\field$ and a finite set $\mathcal L$ of 3 or more skew lines in $\PP^3_\field$, we associate to it a group $G_{\mathcal L}$ and a groupoid $C_{\mathcal L}$ whose action on the union $\cup_{L\in\mathcal L}L$ provides orbits which have a rich combinatorial structure. We characterize when $G_{\mathcal L}$ is abelian and give partial results on its finiteness. The notion of \emph{collinearly complete} subsets is introduced and shown to correspond exactly to unions of groupoid orbits.

In the case where $\field$ is a finite field and $\mathcal L$ is a full spread in $\PP^3_\field$ (i.e., every point of $\PP^3_\field$ lies on a line in $\mathcal L$), we prove that $G_{\mathcal L}$ being abelian characterizes the classical spread given by the fibers of the Hopf fibration.

Over any algebraically closed field, we establish that finite unions of $C_{\mathcal L}$-orbits are geproci sets—that is, finite sets whose general projections to a plane are complete intersections. Furthermore, we prove a converse: if $\field$ is algebraically closed and $Z \subset \PP^3_\field$ is a geproci set consisting of $m$ points on each of $s \geq 3$ skew lines $\mathcal L$ where the general projection of $Z$ is a complete intersection of type $(m, s)$, then $Z$ is a finite union of orbits of $C_{\mathcal L}$.

This work thus uncovers a profound combinatorial framework governing geproci sets, providing a new bridge between incidence combinatorics and algebraic geometry.
\end{abstract}



\date{June 4, 2025} 
\maketitle


\section{Introduction}\label{sec:Intro}
We work in the projective space $\PP^3_\field$ over an arbitrary field $\field$. We omit the subscript if the field is understood or irrelevant.

Finite sets of lines in $\PP^3$ have been studied from various points of view. For example, 
\cite{hartley1993, hartley1994}
study invariants of 4 skew lines
from the perspective of computer vision
and image reconstruction, 
\cite{huang1991}
studies rings of invariants
for general sets of skew lines, and 
\cite{srinivasan2021}
obtains enumerative geometric results related to transversals of
4 skew lines in $\PP^3$ over finite extensions of the prime field.

In this work, we uncover and investigate previously unrecognized intrinsic combinatorial and algebraic structures arising from any finite set $\call = \{L_1, \ldots, L_s\}$ of $s \geq 3$ skew lines in $\PP^3$. Specifically, we introduce a groupoid $C_\call$ (i.e., a category with all morphisms invertible; see Definition~\ref{groupoidCL}) and a group $G_\call$, both naturally associated to $\call$. The group $G_\call$ acts on the points of each individual line $L_i$, while the groupoid $C_\call$ acts compatibly on the union $\bigcup_i L_i$. This groupoid and its action is interesting not only combinatorially but also because the orbits, when finite, have a rich and geometrically meaningful structure: they are \emph{geproci} sets.
In this paper we obtain results about the combinatorics of skew lines (discussed in more detail in our Overview of Results), which we then apply to obtain results about geproci sets.

Recall that a finite set of points $Z$ in $\PP^3$ is an $\{a,b\}$-geproci set (with $a,b$ positive integers), if its projection from a general point to a plane in $\PP^3$ is a complete intersection of two curves of degrees $a$ and $b$. 

If the set $Z$ is itself a complete intersection of two curves of degree $a$ and $b$ contained in some plane $H_Z\subset\PP^3$ (i.e., $Z$ is \emph{degenerate}), then it is an $\{a,b\}$-geproci set in a trivial way. Nondegenerate geproci sets are more interesting. One of our main results is Proposition \ref{curve prop}, which asserts that 
\begin{center}
\begin{minipage}{0.9\textwidth}
  {\it   If the least degree $c$ of a curve $C$ containing a nondegenerate $\{a,b\}$-geproci set $Z$ satisfies $c\leq \max(a,b)$, then $c=a$ or $c=b$, and $C$ is a union of $c$ disjoint lines, each of which contains the same number, $ab/c$, of points of $Z$.}
\end{minipage}
\end{center}
As a consequence of this result, the combinatorics of skew lines plays an important role when considering geproci sets. To assist in this consideration, we will say that $Z$ is an $[a,b]$-geproci set if $Z$ is an $\{a,b\}$-geproci set and it is the union of $a$ points on each of $b$ disjoint lines. 

If $Z$ is $[a,b]$- and $[b,a]$-geproci, then $Z$ is a \emph{grid}, i.e., it is the intersection of $a$ skew lines each of which meet each of $b$ skew lines in a single point. When $a,b\geq3$, these two sets of skew lines belong to
separate rulings on a smooth quadric in $\PP^3$. Grids are easy to classify \cite{CM,CDFPGR,Polizzi}.

If $Z$ is $[a,b]$- but not $[b,a]$-geproci, then is called a \emph{half grid}. Related to this, one of our main results is Theorem \ref{CHGthm} which asserts that
\begin{center}
\begin{minipage}{0.9\textwidth}
   \it{A finite set of points $Z$ consisting of $a$ points on each of $b\geq 3$ skew lines $\mathcal L$ is $[a,b]$-geproci if and only if it is a union of  $C_{\mathcal L}$ orbits.}
\end{minipage}
\end{center}
By this result understanding half grids 
boils down to understanding the combinatorics of skew lines. Although nondegenerate sets which are $\{a,b\}$-geproci but not $[a,b]$- or $[b,a]$-geproci (viz., nondegenerate
nongrid non-half grids) are still mysterious, especially in characteristic $0$,
our results here may be an important first step in understanding them. 
So far in characteristic $0$ only three such sets of points have been discovered. However these three decompose as disjoint unions of half grids; when such a decomposition exists, the theory developed here will be relevant. In finite characteristics more examples are known \cite{kettinger2023} based on combinatorial results about spreads, but if such decompositions are intrinsic to such sets, our results may be relevant for nondegenerate
nongrid non-half grids quite generally. 
Moreover, the combinatorial structures on $\call$ also suggest
a potential for contributing towards the classification of spreads in finite characteristics, see Theorem \ref{Thm:HopfFibr}.

Thus there is an array of ambitious questions for which it is essential to develop a better understanding of the algebraic and combinatorial properties of $G_\call$ and $C_\call$. In the next section, we summarize our findings in this direction.

\subsection{Overview of Results}
We now outline the main questions that have guided our research and we provide a roadmap indicating where in this manuscript to find detailed arguments and proofs that address these questions.
\vskip\baselineskip

$\bullet$ When is $G_{\mathcal L}$ abelian? 
\vskip\baselineskip

We provide two answers to this question.
Since $G_\call$ has a finite set of generators (as given in Proposition~\ref{Prop1}(a)), it is abelian if and only if its generators commute.

An alternative characterization is given in Theorem~\ref{TransThm}, stated in terms of transversals—i.e., lines that intersect all lines in $\call$. This theorem asserts that $G_\call$ is abelian if and only if $\call$ admits two or more transversals, counted with multiplicity.
\vskip\baselineskip

$\bullet$ What is the order of $G_{\mathcal L}$? When is it infinite?
What groups arise?
\vskip\baselineskip

We show that $G_\call$ is trivial (i.e., has order 1) if and only if all lines $L_i \in \call$ lie on a single smooth quadric surface (Corollary~\ref{transCor}). By Proposition~\ref{Prop1}(b), $G_\call$ is infinite if and only if it contains an element of infinite order. In particular, if $G_\call$ is abelian, then it is finite exactly when all its generators have finite order.

We partially address infiniteness over algebraic extensions of prime fields. If $\call$ is defined over a finite field $\field$ of characteristic $p$, then $G_\call$ embeds into the finite group $PGL_2(\field)$ and is finite. For $\overline{\field} = \overline{\mathbb{Q}}$, we give an algorithmic test for finiteness (see Remark~\ref{InfRem}), but a full characterization remains open.

Regarding structure, we classify all finite abelian groups that arise as $G_\call$ over an algebraically closed field. If $G_\call$ is abelian, it is a subgroup of either the additive or multiplicative group of $\overline{\field}$ (Theorem~\ref{TransThm}, Propositions~\ref{twoTrans} and~\ref{oneTransProp}). Moreover, every such finite abelian subgroup occurs for some configuration $\call$ (Theorem~\ref{AbelianGrpClassifThm}). In particular, these groups are either cyclic or finite-dimensional $\mathbb{Z}/p\mathbb{Z}$-vector spaces.

In the more general (nonabelian) case, if $\field \subset \mathbb{C}$ and $G_\call$ is finite, then it must be cyclic, dihedral, or one of the exceptional groups of order $24$, $48$, or $120$ (cf. \cite[Theorem 6.11]{LW}). Whether all such groups occur as $G_\call$ is unknown. The same question remains open in positive characteristic
(see \cite{F} for a classification of 
the finite subgroups of $PGL_2(\field)$).
\vskip\baselineskip
$\bullet$ How is $|G_{\call}|$ related to $|\call|$? For each $b$, among all sets $\call$ of $b$ skew lines with $|G_{\call}| > 1$, is $G_{\call}$ abelian when its order is minimal?
\vskip\baselineskip

Our intuition is that it is difficult for $|G_{\call}|$ to be small when $|\call|$ is large. Indeed, if $\call$ consists of $b \geq 3$ lines from the same ruling on a smooth quadric, then $|G_{\call}| = 1$, so the group order does not, in general, bound the number of lines.

However, we provide explicit lower bounds in Corollary~\ref{LineCountCor} under the assumption that $|G_{\call}| > 1$ and that $\call$ has two transversals (in which case $G_{\call}$ is abelian). These bounds are likely not optimal, but current evidence suggests the following:

\begin{itemize}
\item In \textbf{positive characteristic}, it may always hold that $|G_{\call}| \geq \sqrt{|\call|-1} + 1$, with equality attained for spreads given by the Hopf fibration (see \cite{Ganger}). In these cases, $G_{\call}$ is abelian.
\item In \textbf{characteristic zero}, we expect that $|G_\call| \geq |\call| - 2$. Equality occurs in the so-called standard construction from \cite{POLITUS}, as shown in \cite{Ganger}, where again $G_{\call}$ is abelian and $\call$ consists of Hopf fibration fibers.
\end{itemize}
\vskip\baselineskip
$\bullet$ When are groupoid orbits finite? Can there be two orbits of different orders (one possibly infinite), both larger than the number of lines?
\vskip\baselineskip
We provide a complete answer in the case when $G_{\mathcal L}$ is infinite.
By Theorem~\ref{Thm:InfOrbits}, in this case, every orbit is either infinite or has order equal to the number of lines. Moreover, there can be at most two finite orbits of this kind, and these arise from transversals, as we now explain.
It is immediate from the construction of the groupoid that every orbit meets every line. Thus, if $\mathcal L = \{L_1, \ldots, L_s\}$ with $s \geq 3$, then every orbit has order at least $s$. Furthermore, by Corollary~\ref{transCor}, each orbit of order exactly $s$ consists of the points $T \cap \bigcup_i L_i$ for some transversal $T$.

In the case where there are no infinite orbits, we have partial results.
By Remark~\ref{rem:order of orbit}, the order of every finite orbit is divisible by $s$. Moreover, if $G_{\mathcal L}$ is a nontrivial abelian group (by Theorem~\ref{TransThm}, Corollary~\ref{oneTransCor}, Proposition~\ref{twoTrans}(a), and Corollary~\ref{transCor}), then there is either one or two orbits of order $s$, and all other orbits are projectively equivalent, each having order $|G_{\mathcal L}| \cdot s$.

\vskip\baselineskip
$\bullet$ How can one computationally check if a set $Z$ of $m$ points on each of $s$ skew lines $\mathcal L$ is $[m,s]$-geproci?
\vskip\baselineskip

It can be difficult to perform the computations required to directly verify whether such a set $Z$ is geproci. However, it is much easier (using Proposition~\ref{p. chg and finite union orbits}(a,b)) to check whether $Z$ is a union of orbits. One can then apply one of our main results, Theorem~\ref{CHGthm}, which states that $Z$ is a finite union of $C_{\mathcal L}$-orbits if and only if $Z$ is $[m,s]$-geproci.

\vskip\baselineskip
$\bullet$ If $Z$ is $[a,b]$-geproci for $b \geq 3$ skew lines $\mathcal L$, what can the group $G_{\mathcal L}$ tell us about $Z$?
\vskip\baselineskip

As a first step, by Corollary~\ref{grids and groups Cor}, $Z$ is a grid if and only if $|G_{\mathcal L}| = 1$.

\vskip\baselineskip
$\bullet$ How many projective equivalence classes of single-orbit $[m,s]$-half grids are there in terms of $m$ and $s$?
\vskip\baselineskip

When $s = 4$, our results yield an algorithm for computing the number of complex projective equivalence classes for any given $m$. We apply this in Theorem~\ref{PrimeThm} to show that there are exactly $(m^2 - 1)/12$ such classes of single-orbit complex $[m,4]$-half grids when $m \geq 5$ is prime.

\section{Combinatorics of skew lines}\label{sec: combinatorics of skew lines}
\subsection{The groupoid associated to a finite set of skew lines}\label{sec:Combinatorics of skew lines}
Consider three distinct skew lines, $L_1,L_2,L_3\subset \PP^3_\field$. We can define a 
linear isomorphism $f_{123}\colon L_1\to L_2$ as follows. For each $p\in L_1$, let
$q=f_{123}(p)$ be the point where the plane spanned by $p$ and $L_3$
intersects $L_2$. 
Alternatively, there is a unique smooth quadric $Q_{123}$ containing all three lines in one ruling. 
Taking the line $M_p$ passing through $p$ from the other ruling we have $q=L_2\cap M_p$.

More generally, let ${\mathcal L}=\{L_1,\ldots,L_s\}$ be a set of $s\geq 3$ distinct skew lines $L_i\subset \PP^3_\field$.
Each choice $L_{i_1},L_{i_2},L_{i_3}$ of three of the lines 
then determines the map $f_{i_1i_2i_3}\colon L_{i_1}\to L_{i_2}$ as above.
The composition $f_{lmn}f_{ijk}\colon L_i\to L_m$ makes sense as long as $l=j$
(so that the image of the first map is the domain of the second), but usually
$f_{jmn}f_{ijk}$ is not of the form
$f_{imr}$ for any $L_r\in {\mathcal L}$.
It follows easily from the definition
that $f_{ijk}$ is inverse to $f_{jik}$ and $f_{kij}f_{jki}f_{ijk}$ is the identity
(which one can see keeping in mind
that the quadric $Q_{i_1i_2i_3}$ containing the lines $L_{i_1},L_{i_2},L_{i_3}$
remains the same under permutation of the three indices).

\begin{remark}
Let $L_1,\ldots,L_s\subset \PP^3_\field$ be skew lines.
For $j\neq k$, the projection from $L_k$ to $L_j$ is a 
rational map $F_{jk}\colon \PP^3\DashedArrow L_j$ defined off $L_k$. 
In fact, the map $f_{ijk}$ is just the restriction of $F_{jk}$ to $L_i$.
For computational purposes, it is useful to note the following.
Let $p=(p_1:p_2:p_3:p_4)\in \PP^3$ not be on $L_k$, let $q=(q_1:q_2:q_3:q_4)\in L_j$
and, for $d=1,2$, let $s_d=(s_{d1}:s_{d2}:s_{d3}:s_{d4})$ 
be distinct points of $L_k$. Then $F_{jk}(p)=q$ if and only if
$$\det \begin{pmatrix}
s_{11} & s_{12} & s_{13} & s_{14}\\
s_{21} & s_{22} & s_{23} & s_{24}\\
p_1 & p_2 & p_3 & p_4\\
q_1 & q_2 & q_3 & q_4\\
\end{pmatrix}=0.
\eqno\qed$$
\end{remark}
\vskip\baselineskip

Recall that a {\em groupoid} is a category in which every morphism is invertible.

\begin{definition}\label{groupoidCL}
To a finite set of skew lines $L_1,\ldots,L_s\subset \PP^3_\field$, $s\geq3$, we  associate the groupoid $C_{\mathcal L}$ 
whose objects are the lines $L_i$ and whose arrows
are the maps $f_{ijk}$ and all possible compositions of these maps.
\end{definition}

Note that since $f_{klm}$ is the inverse of $f_{lkm}$, and since the elements of
${\rm Hom}_{C_{\mathcal L}}(L_i,L_j)$ are invertible maps and compositions of invertible maps, 
the elements all are invertible.
In particular, $G_i({\mathcal L})={\rm Hom}_{C_{\mathcal L}}(L_i,L_i)$ 
is a subgroup of the general linear group ${\rm Aut}(L_i)\cong PGL(2,\field)$. If $|\field|<\infty$,
then $G_i({\mathcal L})$ is a subgroup of the group of permutations
of the finitely many points of $L_i$. Since the groups $G_i({\mathcal L})$ all are (noncanonically) isomorphic
(see (c) of the next result), we write $G_{\mathcal L}$ for these groups,
which we refer to as the group of the groupoid.
When ${\mathcal L}$ is clear, we will sometimes abbreviate $G_i({\mathcal L})$ as $G_i$.

\begin{proposition}\label{Prop1}
Let ${\mathcal L}$ be a set of $s\geq 3$ distinct skew lines in $\PP^3_\field$.
\begin{enumerate}[label=\textnormal{(\alph*)}]
\item $G_i$ is generated by the maps of the form $f_{jil}f_{ijk}$ and
$f_{kij}f_{jkl}f_{ijk}$ (and hence $G_i$ is finitely generated). 
\item $G_i$ is finite if and only if every element has finite order.
\item The groups $G_i$ are all isomorphic and the sets 
${\rm Hom}_{C_{\mathcal L}}(L_i,L_j)$ all have the same cardinality, so the set of arrows in $C_{\mathcal L}$ is finite
if and only if $G_i$ is finite for some $i$. 
\item $G_i$ is finite abelian if and only if the elements of the form $f_{jis}f_{ijk}$ and
$f_{kij}f_{jkl}f_{ijk}$ commute and have finite order.
\item Let $\field\subset \field'$ be an extension of fields. Let ${\mathcal L'}$ be the lines 
in $\PP^3_{\field'}$ coming from the lines ${\mathcal L}$. The bijection
${\mathcal L}\to{\mathcal L}'$ given by $L\in {\mathcal L}\mapsto L'\in {\mathcal L'}$
(where $L'$ is the line coming from $L$) induces a canonical isomorphism 
$C_{\mathcal L}'\cong C_{\mathcal L}$ of groupoids. In particular,
the groups of the groupoids are canonically isomorphic.
\item Let $\field\subset \field'$ be an extension of fields, let
${\mathcal L}'$ be the lines in $\PP^3_{\field'}$ coming from 
${\mathcal L}$ and let ${\mathcal L}''$ be a finite set of skew lines in $\PP^3_{\field'}$
with ${\mathcal L}'\subseteq{\mathcal L}''$.
Then $C_{{\mathcal L}'}$ is a subgroupoid of $C_{{\mathcal L}''}$, so
$G_{\mathcal L}$ is isomorphic to a subgroup of $G_{{\mathcal L}''}$.
\item If $\field$ is finite, then $G_i$ is a subgroup
of the group of permutations on the points of $L_i$,
hence $G_i$ is finite.
\end{enumerate}
\end{proposition} 

\begin{proof}
(a) Let $G$ be the group generated by the maps of the form
$f_{jil}f_{ijk}$ and $f_{kij}f_{jkl}f_{ijk}$.
Clearly, $G\subseteq G_i$, so to show $G=G_i$ it is enough to show
$G_i\subseteq G$. 

Let $h\in G_i$ be a product of $r$ maps, each of which is of the form $f_{jkl}$. 
There are no maps in $G_i$ of the form
$f_{jkl}$, so we must have $r\geq2$.
When $r=2$, we have by definition $h\in G$.
If $r>2$, and if we assume that all products in $G_i$ of $r'<r$ maps of the form $f_{jkl}$
are in $G$, we will show that $h\in G$. Then $G_i\subseteq G$ follows by induction on $r$.
 The last two maps in $h$ must be of the form 
$f_{jil}f_{tjk}$ since we end up back at $L_i$. Thus $h=f_{jil}f_{tjk}g$
for some product $g$ of $r-2$ of the maps. 

If $t=i$, then $f_{jil}f_{tjk}=f_{jil}f_{ijk}\in G$, so
$g=(f_{jil}f_{ijk})^{-1}h\in G_i$, so $g\in G$ by induction,
and hence also $h=(f_{jil}f_{tjk})g\in G$.

Say $t\neq i$. Then $f_{jil}f_{tjk}=f_{jil}f_{tjk}f_{itj}f_{tij}
=f_{jil}f_{ijt}f_{jit}f_{tjk}f_{itj}f_{tij}$ so
$h=(f_{jil}f_{ijt})(f_{jit}f_{tjk}f_{itj})(f_{tij}g)$
but $f_{jil}f_{ijt}, f_{jit}f_{tjk}f_{itj}\in G$ and
$f_{tij}g$ is a composition of $r-1$ maps and hence in $G$ by induction,
so $h\in G$.

(b) Let $\overline{G}\subset \PP {\rm GL}(n,\overline{\field})$
be a finitely generated torsion subgroup. It is well known to experts that 
$\overline{G}$ must be finite
(see \url{https://math.stackexchange.com/questions/4071442}).
For lack of a  suitable reference we sketch a proof. 
We have $\PP {\rm SL}(n,\overline{\field})= \PP {\rm GL}(n,\overline{\field})$.
Let $\pi\colon {\rm SL}(n,\overline{\field})\to \PP {\rm SL}(n,\overline{\field})$ 
be the canonical quotient. 
For $g\in {\rm SL}(n,\overline{\field})$, if $\pi(g)$ has finite order,
then so does $g$. Thus $G=\pi^{-1}(\overline{G})$
is a torsion group. Let $G'\subset G$ be the group generated by a finite set of elements
of $G$ that map to generators of $\overline{G}$.
Thus $G'$ is a finitely generated torsion group whose image under $\pi$ is $\overline{G}$. 
By Burnside's Theorem 
\cite[Corollary 3]{Alperin}
(or 
\cite[Theorem G, p. 105]{Kaplansky}) 
$G'$ is finite and hence $\overline{G}$ is finite.

(c) Since $G_i\subseteq C_{\mathcal L}$ for all $i$, if $C_{\mathcal L}$ is finite, then so is $G_i$.

The idea of the converse is that
all of the hom sets ${\rm Hom}_{C_{\mathcal L}}(L_i,L_j)$
have the same cardinality, and there are only finitely many
of them, so the union is finite
if any hom set is finite. 
We now give the details.
Since $g\mapsto f_{ijk}gf_{jik}$
is an isomorphism $G_i\to G_j$,
all the hom groups $G_i$ have the same
cardinality.
Likewise, $g\mapsto f_{ijk}g$ gives a bijection $G_i\to {\rm Hom}_{C_{\mathcal L}}(L_i,L_j)$
(with inverse $g\mapsto f_{jik}g$),
so all the hom sets have the same cardinality too. So if $G_i$ is finite for any $i$, so is the set of all arrows in
$C_{\mathcal L}$.

(d) Assume $G_i$ is finite abelian. 
Then the elements commute and have finite order.
Conversely, if the generators of a group commute, then the group is abelian,
and if they also have finite order then (since the group is abelian)
the group has finite order (which divides the least common multiple
of the orders of the generators).

(e) The lines in ${\mathcal L'}$ are defined over $\field$, and so are the maps
$f'_{ijk}$, hence the restriction of $f'_{ijk}$ to $L_i$ is $f_{ijk}$ and 
$f_{ijk}$ extends canonically to $f'_{ijk}$.
Thus the map from
${\rm Hom}_{C_{\mathcal L'}}(L_i',L_j')$
to ${\rm Hom}_{C_{\mathcal L}}(L_i,L_j)$ is given by restriction, and the inverse
is given by the canonical extensions.

(f) The maps $f'_{ijk}$
coming from ${\mathcal L}'$ are a subset of those coming from 
${\mathcal L}'$ so ${\rm Hom}_{C_{{\mathcal L}'}}(L_i,L_j)$
is a subset of ${\rm Hom}_{C_{{\mathcal L}''}}(L_i,L_j)$.

(g) The maps in ${\rm Hom}_{C_{\mathcal L}}(L_i,L_i)$
are bijective and defined over $\field$ so permute the finitely many
points of $L_i$.
\end{proof}

\begin{definition}\label{rem:order of orbit}
Given a point $p\in \cup_kL_k$, so $p\in L_i$ for some $i$,
we define its $C_{\mathcal L}$-orbit $[p]$ to be
$$[p]=\bigcup_j \{gp : g\in {\rm Hom}_{C_{\mathcal L}}(L_i,L_j)\}.$$
\end{definition}

Note that if $q\in [p]$, then $[q]=[p]$.
We will denote $[p]\cap L_i$ by $[p]_i$. If $p\in L_i$,
note that $[p]_i=G_ip$; i.e., $[p]_i$ is the orbit of $p$ under the action of the group
$G_i$ on $L_i$.

When $T$ is a transversal for $\mathcal L$, the next result shows that
the orbit $[p]$ for any point $p\in T\cap (\cup_i L_i)$ is very simple.

\begin{corollary}\label{transCor}
Consider $s\geq 3$ distinct skew lines $\mathcal L=\{L_1,\ldots,L_s\}$ in $\PP^3$.
Let $p\in L_i$. Then $[p]_i=\{p\}$ if and only if
$p\in T$ for some transversal $T$ for $\mathcal L$,
in which case $[p]=T\cap(\cup L_i)$.
In particular, $|(G_\mathcal L)_i|=1$ if and only if the lines $L_j$ are all contained in a smooth quadric.
\end{corollary}

\begin{proof} Let $G_i=(G_\mathcal L)_i$.
Without loss of generality we may assume $L_i$ is $L_1$.

If $p\in L_1\cap T$ for some transversal $T$, then the plane $\Pi_k$ spanned by $p$ and $L_k$
for any $k\neq 1$ contains $T$, so $f_{1jk}(p)=L_j\cap \Pi_k=L_j\cap T$.
In particular, $[p]=T\cap(L_1\cup\cdots\cup L_s)$, so 
$\{p\}\subseteq[p]_1\subseteq [p]\cap L_1=T\cap L_1=\{p\}$.
Thus we see that $[p]_1=\{p\}$.

Now assume there is no transversal through $p$. Let $p_2=f_{123}(p)$ and $p_3=f_{132}(p)$;
then $p,p_2,p_3$ are collinear (contained in the line $L$ through $p$ meeting
$L_1,L_2,L_3$). Since $L$ is not a transversal for all
of the lines, there is a line $L_j$ which does not meet $L$.
Let $p_j=f_{2j1}(p_2)$.
Then $p_j$ and $p_2$ are collinear with a point $q\in L_1$, namely $q=f_{j12}(p_j)$. Since $p_j\in L_j$
we know $p_j$ is not on $L$ and hence $q$ is not on $L$, so $q\neq p$.
But $q\in[p]_1$, so $[p]_1\neq \{p\}$.

Finally, if all of the lines $L_i$ are contained in a smooth quadric, then every point $p\in L_i$
is on a transversal $T$, so $|[p]_i|=1$ for all $p\in L_i$, hence
$|G_i|=1$. Conversely, assume $|G_i|=1$ for some $i$. Then $|G_j|=1$ for all $j$,
so $|[p]_j|=1$ for all $p\in L_j$,
hence every point $p\in L_1$ is on a transversal to $\mathcal L$.
The union of the transversals for $L_2,L_3,L_4$ is the quadric $Q_{234}$. Thus every 
transversal to $\mathcal L$ lies in $Q_{234}$. This implies that $L_1\subset Q_{234}$.
By symmetry, we conclude that $Q_{234}$ contains all of the lines $L_i$.
\end{proof}

The preceding result characterizes orbits which are as small as possible.
The next result characterizes the situation when an orbit is infinite.

\begin{theorem}\label{Thm:InfOrbits}
Consider a set of skew lines ${\mathcal L}=\{L_1,\ldots,L_s\}$ in $\PP^3_{\overline{\field}}$.
The groupoid $C_{\mathcal L}$ has an infinite orbit if and only if
$|G_{\mathcal L}|=\infty$, in which case all but at most two orbits are infinite,
and the finite orbits, if any, have order $s$ and come from transversals for $\mathcal L$.
\end{theorem}

\begin{proof}
If $C_{\mathcal L}$ has an infinite orbit, then $G_1$ is infinite. Conversely, if $G_1$ is infinite, then
by Proposition \ref{Prop1}, $G_1$ has an element $g$ of infinite order.
Regarded as an automorphism of $L_1$, $g$ is represented by a 
$2\times2$ matrix $M_g$ and has either one eigenspace of dimension 1 or 2, or two eigenspaces of dimension 1.
If $M_g$ has one eigenspace of dimension 2, $M_g=cI_2$ for a nonzero scalar $c$ and so 
$g$ is the identity of $G_1$.
If $M_g$ has one eigenspace of dimension 1, then up to similarity we can choose $M_g$ to be
$\begin{pmatrix}
1 & 0 \\
b & 1 \\
\end{pmatrix}$ with $b\neq0$.
This has infinite order if ${\rm char}(\overline{K})=0$ but only finite order equal to ${\rm char}(\overline{K})$ in positive characteristics,
so we may assume ${\rm char}(\overline{K})=0$.
In essence $g$ has a fixed point $p$ which we take to be infinity
and then $g$ is a translation of the affine part of the line by $b$.
Thus ${\rm char}(\overline{K})=0$ and the orbit of every point other than $p$
is infinite, while $[p]_1=\{p\}$ (in which case the orbit $[p]$ comes from a transversal
by Corollary \ref{transCor}).
If $M_g$ and has two eigenspaces of dimension 1,
then $g$ has two fixed points, say $p$ and $q$
which we can regard as 0 and $\infty$.
Thus $g$ corresponds to a scaling by a nonzero element
of $\overline{\field}$ of infinite multiplicative order, hence
again the orbit of every point other than $p$ and $q$ 
is infinite. The finite orbits are either $[p]_1=\{p\}$ and $[q]_1=\{q\}$ (in which case the orbits
come from transversals) or $[p]_1=[q]_1=\{p,q\}$, but this is impossible by the next result.
\end{proof}
We now show that groupoid orbits for $s\geq3$ skew lines $\mathcal L$ cannot have order
$2s$, and that $|G_{\mathcal L}|$ cannot be 2.

\begin{proposition}\label{Prop:not 2}
Let ${\mathcal L}$ be $s\geq3$ distinct skew lines in $\PP^3$.
If $p\in \cup_{L\in\mathcal L} L$ is not on a transversal for ${\mathcal L}$, then $|[p]|>2s$
and hence $|G_{\mathcal L}|>2$.
\end{proposition}

\begin{proof}
Since $p$ is not on a transversal for ${\mathcal L}$ we have $s\geq4$ so we may pick four lines,
$L_1,L_2,L_3,L_4\in{\mathcal L}$ such that $p\in L_1$ but $L_4$ is not in the quadric $Q_{123}$ determined by 
$L_1,L_2,L_3$. There is a unique transversal $T$ for $L_1,L_2,L_3$ through $p$.
(Recall there are two lines through $p$ contained in $Q_{123}$, one from each ruling; 
$L_1$ is one and $T$ is the other.) Let $p_i$ be the points $T\cap L_i$ for $i=1,2,3$, so $p=p_1$
and $p_i\in [p]$. Let $p_{12}=f_{214}f_{123}(p_1)$, so $p_{12}$ is the point where the plane
$\langle p_2,L_4\rangle$ spanned by $p_2$ and $L_4$ meets $L_1$, and let
$p_{13}=f_{314}f_{132}(p_1)$, so $p_{13}$ is the point where the plane
$\langle p_3,L_4\rangle$ spanned by $p_3$ and $L_4$ meets $L_1$.
Note that the line $\langle p_{12},p_2\rangle$ meets $L_4$, so
if $p_1 = p_{12}$, then $T$ meets $L_4$ hence $T$ is a transversal for $L_1,L_2,L_3,L_4$,
contrary to assumption. Thus $p_1\neq p_{12}$ and similarly $p_1\neq p_{13}$.
If $p_{12}= p_{13}$, then $\langle p_2,L_4\rangle=\langle p_{12},L_4\rangle=
\langle p_{13},L_4\rangle=\langle p_3,L_4\rangle$, which again means $T=\langle p_2,p_3\rangle$
meets $L_4$. Thus $p_{12}\neq p_{13}$, hence $[p]$ meets $L_1$ in at least three
distinct points, so $|G_{\mathcal L}|>2$ and $|[p]|>2s$.
\end{proof}

We note that orbits either span a line or they span all of $\PP^3$:

\begin{proposition}\label{prop:projEqOverKbarVsK}
Let ${\mathcal L}$ be $s\geq3$ distinct skew lines in $\PP^3_\field$.
Let $O_1,O_2\subseteq\PP^3_\field$ be $C_{\mathcal L}$ orbits
which are projectively equivalent in $\PP^3_{\overline\field}$.
Then they are projectively equivalent in $\PP^3_\field$.
Moreover, either $O_1$ is contained in a transversal for ${\mathcal L}$
or $O_1$ contains 5 linearly general points.
\end{proposition}

\begin{proof}
Let $\Phi$ be the element of $PGL_4(\overline\field)$ taking
$O_1$ to $O_2$. Let $p_i$ be the points of $O_1$
with respect to some index set $I$,
and let $q_i=\Phi(p_i)$ be the points of $O_2$.

If $|O_1|=s$, then $|O_2|=s$ and by Corollary \ref{transCor} 
there are transversals $T_i$ with $O_i = T_i\cap (\cup_{L\in\mathcal L})$.
Thus we can assume $I=\{1,2,\dots,s\}$. There is an element
of $\phi\in PGL_4(\field)$ taking $p_i$ to $q_i$ for $i=1,2,3$,
since the points are defined over $\field$. This defines
a map $\phi_{T_1}\colon T_1\to T_2$ and must be the same map as 
$\Phi_{T_1}\colon T_1\to T_2$, hence $\phi$ takes $O_1$ to $O_2$.

Now say $|O_1|>s$. Then $|O_1\cap L|\geq3$ for each $L\in {\mathcal L}$
by Proposition \ref{Prop:not 2}. Pick three lines $L_1,L_2,L_3\in{\mathcal L}$
and three points $p_1,p_2,p_3\in O_1\cap L_1$.
Then we have three lines $T_i$ transversal for $L_1,L_2,L_3$
with $p_i\in O_1\cap T_i$. Let $p_{ij}$ be the point $T_i\cap L_j$. 
Then the points $p_{11},p_{13},p_{31},p_{33},p_{22}$ are linearly general
(see Figure \ref{GridOnTetrahedron}) hence so are $q_{11},q_{13},q_{31},q_{33},q_{22}$
where $q_{ij}=\Phi(p_{ij})$.
Since the points are defined over $\field$ and there is a unique linear map
taking $p_{11},p_{13},p_{31},p_{33},p_{22}$ to 
$q_{11},q_{13},q_{31},q_{33},q_{22}$, this map, namely $\Phi$,
is defined over $\field$.
\end{proof}

\begin{figure}[ht!]
\definecolor{ffffff}{rgb}{1.,1.,1.}
\definecolor{uuuuuu}{rgb}{0.26666666666666666,0.26666666666666666,0.26666666666666666}
\definecolor{ududff}{rgb}{0.30196078431372547,0.30196078431372547,1.}
\begin{tikzpicture}[x=1cm,y=1cm]
\clip(-2.9251177150754075,-0.3643893007173551) rectangle (4.096570195978351,5.0108338587789705);
\draw [line width=1pt] (-2.25376510915961,1.4988893252096052)-- (1.2680846267921186,0.4863575261234834);
\draw [line width=1pt] (1.2680846267921186,0.4863575261234834)-- (2.92,1.62);
\draw [line width=1pt] (2.92,1.62)-- (0.64,4.32);
\draw [line width=1pt] (0.64,4.32)-- (-2.25376510915961,1.4988893252096052);
\draw [line width=1.pt,dash pattern=on 3pt off 3pt] (-2.25376510915961,1.4988893252096052)-- (2.92,1.62);
\draw [line width=1.pt,color=white] (0.95,1.574)-- (1.5,1.587);
\draw [line width=1.pt,color=white] (0.35,1.56)-- (-0.1,1.549);
\draw [line width=1.pt] (-0.4928402411837456,0.9926234256665443)-- (1.78,2.97);
\draw [line width=1.pt] (2.0940423133960593,1.0531787630617417)-- (-0.8068825545798048,2.9094446626048027);
\draw [line width=1.pt] (0.877978830884413,2.1852393744871916)-- (1.1147070976400486,2.391193546136915);
\draw [line width=3pt,color=white] (0.9571298794272152,1.780674950653604)-- (1.1510884940726007,1.6565632350199406);
\draw [line width=3pt,color=white] (0.877978830884413,2.1852393744871916)-- (1.1147070976400486,2.391193546136915);
\draw [line width=1pt] (0.64,4.32)--(1.2680846267921186,0.4863575261234834);
\draw [fill=black] (-2.25376510915961,1.4988893252096052) circle (4pt);
\draw[color=black] (-2.392437942512709,1.2154904792697407) node {$p_{11}$};
\draw [fill=black] (1.2680846267921186,0.4863575261234834) circle (4pt);
\draw[color=black] (1.25,0.2) node {$p_{31}$};
\draw [fill=black] (2.92,1.62) circle (4pt);
\draw[color=black] (3.3,1.9) node {$p_{33}$};
\draw [fill=black] (0.64,4.32) circle (4pt);

\draw[color=black] (0.8278534097981529,4.7) node {$p_{13}$};
\draw [fill=white] (-0.8068825545798048,2.9094446626048027) circle (4pt);

\draw[color=black] (-1.2,2.1) node {$T_1$};
\draw[color=black] (.25,1.25) node {$T_2$};
\draw[color=black] (1.9,0.6) node {$T_3$};
\draw[color=black] (-1.4,1) node {$L_1$};
\draw[color=black] (0,2.7) node {$L_2$};
\draw[color=black] (1.5,3.8) node {$L_3$};

\draw[color=black] (-1.2544402465833064,3.0677633247718528) node {$p_{12}$};
\draw [fill=white] (-0.4928402411837456,0.9926234256665443) circle (4pt);
\draw[color=black] (-0.7096541155532734,0.7) node {$p_{21}$};
\draw [fill=white] (2.0940423133960593,1.0531787630617417) circle (4pt);
\draw[color=black] (2.6,0.9733633099230593) node {$p_{32}$};
\draw [fill=white] (1.78,2.97) circle (4pt);
\draw[color=black] (2.2,3.25) node {$p_{23}$};
\draw [fill=black] (0.6435798794081273,1.9813117128332722) circle (4pt);
\draw[color=black] (0.08936554329077495,2.014510138113789) node {$p_{22}$};
\end{tikzpicture}
\caption{Nine points $p_{ij}$ of intersection of three skew lines $L_j$ and three transversals $T_i$ to the $L_j$,
have 5 (shown as black dots) of the 9 points which are linearly general.}
\label{GridOnTetrahedron}
\end{figure}

The next result shows that projectively equivalent sets of skew lines have isomorphic groupoids
and hence corresponding points have projectively equivalent orbits.

\begin{proposition}\label{FunctorCor}
Let $L_1,\ldots,L_s\subset \PP^3_\field$ be $s\geq3$ distinct skew lines, let $f\in {\rm PGL}(\field)$ be 
an automorphism of $\PP^3_\field$ and let $L_i'=f(L_i)$. Denote $\{L_1,\ldots,L_s\}$ by ${\mathcal L}$ and
$\{L_1',\ldots,L_s'\}$ by ${\mathcal L'}$.
Then $f$ induces a canonical isomorphism $\phi_f\colon C_{\mathcal L}\to C_{\mathcal L'}$
and for any point $p\in L_1\cup\cdots\cup L_s$, the orbits $[p]$ and $[f(p)]$
are projectively equivalent.
\end{proposition}

\begin{proof}
Note that $\phi_f$ is a functor, where $\phi_f(L_i)=L_i'$ and for $g\in{\rm Hom}_{C_{\mathcal L}}(L_i,L_j)$
we have $\phi_f(g)$ being the map $\phi_f(g)=fgf^{-1}\colon L_i'\to L_j'$. It is an isomorphism since
its inverse is $\phi_{f^{-1}}$. If $p\in L_i$, then $[p]=\{g(p): g\in \bigcup_j{\rm Hom}_{C_{\mathcal L}}(L_i,L_j)\}$,
and $[f(p)]=\{g(f(p)): g\in \bigcup_j{\rm Hom}_{C_{\mathcal L}}(L_i',L_j')\}$ but $f(g(p))=\phi_f(g)(f(p))$ so 
$$[f(p)]=\{g(f(p)): g\in \bigcup_j{\rm Hom}_{C_{\mathcal L}}(L_i',L_j')\}=\{\phi_f(g)(f(p)): g\in \bigcup_j{\rm Hom}_{C_{\mathcal L}}(L_i,L_j)\}$$
$$=\{f(g(p)): g\in \bigcup_j{\rm Hom}_{C_{\mathcal L}}(L_i,L_j)\}=f(\{g(p): g\in \bigcup_j{\rm Hom}_{C_{\mathcal L}}(L_i,L_j)\}=f([p]).\eqno\qedhere$$
\end{proof}

\begin{definition}
Let ${\mathcal L}=\{L_1,\ldots,L_s\}\subset\PP^3$ be $s\geq3$ distinct skew lines.
Let $Z$ be a nonempty but possibly infinite subset of $\cup_iL_i$.
We will say $Z$ is \emph{collinearly complete with respect to ${\mathcal L}$} if whenever 
$T$ is a transversal for 3 or more of the lines $L_i$ such that $T\cap Z$ is nonempty, then $T\cap(\cup L_i)\subset Z$.
\end{definition}

\begin{proposition}\label{p. chg and finite union orbits}
Let ${\mathcal L}=\{L_1,\ldots,L_s\}\subset\PP^3$ be $s\geq3$ distinct skew lines.
Let $Z$ be a subset of $\cup_iL_i$.
The following are equivalent:
\begin{enumerate}[label=\textnormal{(\alph*)}]
\item $Z$ is collinearly complete with respect to ${\mathcal L}$;
\item $Z$ is a union of orbits for ${\mathcal L}$.
\end{enumerate}
Let $p\in\cup_iL_i$ and let $1\leq j<k\leq s$. We have the following:
\begin{enumerate}
\item If $1\leq j,k,l \leq s$ are distinct, then $f_{jkl}([p]_j)=[p]_k$, hence $|[p]_j|=|[p]_k|$.
\item All orbits $[p]$ are finite if and only if $G_i$ is finite for some $i$ (or equivalently, every $i$).
\item If $Z$ is collinearly complete with respect to ${\mathcal L}$, then $|Z\cap L_j|=|Z\cap L_k|$.
\end{enumerate}
\end{proposition}

\begin{proof}
We first prove (1).
For each $q\in [p]_j$, we have $f_{jkl}(q)\in [p]\cap L_k=[p]_k$ so $f_{jkl}([p]_j)\subseteq [p]_k$.
Similarly, $f_{kjl}([p]_k)\subseteq [p]_j$ but $f_{jkl}$ and $f_{kjl}$ are inverse bijections,
hence $f_{jkl}([p]_j)=[p]_k$. But $f_{jkl}$ is bijective, so $|[p]_j|=|[p]_k|$.

We now prove (2). If $G_i$ is finite, then ${\rm Hom}_{C_{\mathcal L}}(L_i,L_j)$ is finite for each $j$, 
so $[p]_j$ is finite for each $j$, so $[p]$ is finite.
To show that finite orbits imply $G_i$ is finite,
let $U$ be the union of the orbits of 3 points of $L_i$, and let $U_i=U\cap L_i$. 
Now $U_i$ is finite, so given any $g\in G_i$, 
some power of $g$ acts trivially on $U_i$, but $U_i$ contains at least 
3 points of $L_i$ so the power of $g$ is trivial. Thus $g$ has finite order,
so $G_i$ is finite by Proposition \ref{Prop1}(b).

We next show (a) implies (b).
Assume $Z$ is collinearly complete with respect to ${\mathcal L}$.
It is enough to show $[p]\subseteq Z$ whenever $p\in Z$,
and for this it is enough to show $f_{ijk}(Z_i)\subseteq Z$
for each triple  $i,j,k$ of distinct indices $1\leq i,j,k \leq s$
where $Z_i=Z\cap L_i$. To show $f_{ijk}(Z_i)\subseteq Z$,
let $p\in Z_i$. Let $T$ be the transversal for $L_i,L_j,L_k$
through $p$. Then $f_{ijk}(p)=T\cap L_j\in Z_j\subseteq Z$ by completeness,
so $f_{ijk}(Z_i)\subseteq Z$ follows.

Now we show (b) implies (a), so assume $Z$ is a union of orbits
(and hence if $p\in Z$ then $[p]\subseteq Z$).
Let $i,j,k$ be distinct indices $1\leq i,j,k \leq s$.
Let $T$ be a transversal for $L_i,L_j,L_k$ through a point $p$ of $Z_i$. 
Let $l$ be an index such that $T$ meets $Z_l$ and let $q$ be the point $T\cap Z_l$.
If $l=i$, then $q=p$ so $q\in [p]\subseteq Z$.
If $l=j$, then $q=f_{ijk}(p)$ so $q\in [p]\subseteq Z$.
If $l\neq i,j$, then $q=f_{ilj}(p)$ so $q\in [p]\subseteq Z$.
Thus $Z$ is collinearly complete with respect to ${\mathcal L}$.

Finally, we prove (3). 
If $Z$ is collinearly complete with respect to ${\mathcal L}$, then 
$Z$ is a disjoint union of orbits. Thus for some index set $I$ 
and points $p_i\in Z$, $i\in I$, we have a disjoint union $Z=\cup_{i\in I} [p_i]$.
Thus $Z\cap L_j=\cup_{i\in I} ([p_i]\cap L_j) = \cup_{i\in I} [p_i]_j$ and 
$Z\cap L_k=\cup_{i\in I} ([p_i]\cap L_k) = \cup_{i\in I} [p_i]_k$,
but $|[p_i]_j|=|[p_i]_k|$ by (1), so $|\cup_{i\in I} [p_i]_j|=|\cup_{i\in I} [p_i]_k|$ 
hence $|Z\cap L_j|=|Z\cap L_k|$.
\end{proof}

\begin{remark}\label{rem:two transversals}
A significant subclass of sets of skew lines $L_1,\ldots,L_s\subset\PP^3$ 
are those having two skew transversals, $T_1$ and $T_2$. Examples are easy to come by.
Given skew lines $T_1$ and $T_2$, note that any lines $L$ and $L'$ which meet both
$T_1$ and $T_2$ have to be skew unless they both meet $T_1$ at the same point or both meet
$T_2$ at the same point. Thus given skew lines $T_1$ and $T_2$, pick any finite subset
$A_1\subset T_1$ and any finite subset $A_2\subset T_2$ with $|A_1|=|A_2|$.
Then pick any bijection $b\colon A_1\to A_2$ and define the lines
$L_a$, $a\in A_1$, where $L_a$ is the line through the points $a$ and $b(a)$.
Then the lines $\mathcal L=\{L_a:a\in A_1\}$ are skew and $T_1$ and $T_2$ are transversals for $\mathcal L$. 
\qed\end{remark}

\begin{remark}\label{r.2 transversal-a}
Assume the ground field $\field$ is algebraically closed. 
Consider lines $L_i\subset\PP^3_\field$.

It is easy to see that any 
line $L_1$ and any two skew lines $L_1,L_2$ have infinitely many transversals defined over $\field$.

Three skew lines $L_1,L_2,L_3$ also have infinitely many transversals defined over $\field$ and the transversals are skew.
(This is because there is a unique smooth quadric $Q$ containing the three lines,
which are all in one ruling of $Q$. By B\'ezout's Theorem 
a line $T$ is a transversal for the three if and only if $T$ is in the other ruling.)

A set of four skew lines $\mathcal L=\{L_1,L_2,L_3,L_4\}$ always has either one, two or infinitely many transversals
defined over $\field$ and if there is more than one transversal, they are themselves skew.
(There is a unique smooth quadric $Q$ containing $L_1,L_2,L_3$ and these three lines
are all members of the same ruling on $Q$. Every transversal
for $\mathcal L$ is also a transversal for $L_1,L_2,L_3$, but a line is transversal for $L_1,L_2,L_3$
if and only if it is a line in the other ruling on $Q$, hence the transversals for $L_1,L_2,L_3$ are skew. 
If $L_4\subset Q$, then $\{L_1,L_2,L_3\}$ and $\mathcal L$ have the same transversals,
hence there are infinitely many. If $L_4\not\subset Q$, then either $L_4$ meets $Q$ in two points,
in which case $\mathcal L$ has exactly two transversals, or $L_4$ meets $Q$ in only one point,
in which case $\mathcal L$ has exactly one transversal; this occurs when $L_4$ is tangent to $Q$.)

Finally, if ${\mathcal L}=\{L_1,\ldots,L_s\}$ is a set of $s>4$ skew lines, then ${\mathcal L}$ has either 0, 1, 2 or infinitely many
transversals defined over $\field$ and for each $s>4$ each possibility occurs. 
(We know transversals for ${\mathcal L}$ are transversals for $\{L_1,\ldots,L_4\}$ so clearly
there are either infinitely many or 2 or fewer transversals. To get examples with infinitely many
transversals, just pick $s$ lines from the same ruling on any smooth quadric.
To get examples with exactly two
transversals, just pick four lines with exactly two transversals, and then pick any $s-4$ additional 
lines meeting both transversals but general among the lines meeting both transversals.
To get examples with exactly one
transversal, just pick four lines with exactly one transversal, and then pick any $s-4$ additional 
lines meeting this transversal but general among the lines meeting it.
And to get examples with no
transversals, just pick four lines with only one or two transversals, and then pick any $s-4$ additional 
general lines.)
\qed\end{remark}

\begin{definition}\label{StPosDef}
We will say that the lines
$L_1,L_2,L_3\subset\mathbb{P}^3$
are in \emph{standard position} and that
$T_1,T_2$ are the 
\emph{standard transversals} if: 

$L_1$ is defined by $y=z=0$,

$L_2$ is defined by $x=y, z=w$,

$L_3$ is defined by $x=w=0$,

$T_1$ is defined by $x=y=0$, and 

$T_2$ is defined by $z=w=0$. 

Note that $L_1,L_2,L_3$ are skew and each meets both
$T_1$ and $T_2$.
We will also say $L_1,L_2,L_3,L_4$ are
in \emph{standard position} if $L_4$
meets both $T_1$ and $T_2$ but does not meet
any of $L_1,L_2,L_3$ and is not on the unique
smooth quadric $Q\colon xz-yw=0$ containing $L_1,L_2,L_3$.
\end{definition}

\begin{remark}\label{StPosRem}
For any skew lines $L_1,L_2,L_3\subset \PP^3_{\field}$
having distinct transversals $T_1,T_2\subset \PP^3_{\field'}$
for an extension $\field\subseteq \field'$,
there is a choice of coordinates on $\PP^3_{\field'}$
such that $L_1,L_2,L_3$ are in standard position and $T_1$ and $T_2$ are the
standard transversals.
One simply chooses coordinates
such that $T_1$ is defined by $x=y=0$, 
$T_2$ by $z=w=0$, and such that:

$L_1$ meets $T_1$ at 0001 and $T_2$ at 1000;

$L_2$ meets $T_1$ at 0011 and $T_2$ at 1100; and

$L_3$ meets $T_1$ at 0010 and $T_2$ at 0100. 

The unique smooth quadric $Q$ containing the lines 
$L_1,L_2,L_3,T_1, T_2$ is defined by $xz-yw=0$. Then
$L_4\cap T_1$ is the point $(0:0:t:1)$ for some $t\neq0,1$ and 
$L_4\cap T_2$ is the point $(l:1:0:0)$ for some $l\neq 0,1$. Note that $L_4$ is the line
defined by $z=tw$ and $x=ly$. Substituting into
$xz-yw$ gives $(lt-1)yw$, so $L_4\not\subset Q$
is equivalent to $lt\neq1$.
\qed\end{remark}

\begin{proposition}\label{twoTrans}
Consider a set of $s\geq 4$ distinct skew lines  ${\mathcal L}=\{L_1,\ldots,L_s\}\subset\PP^3_\field$.
Assume there is no smooth quadric containing all of the lines $L_i$
but that there are two distinct lines $T_1, T_2$ transversal to ${\mathcal L}$ and also defined over $\field$. Then we have the following facts.
\begin{enumerate}
\item[(a)] Each orbit
under the action of $C_{\mathcal L}$ is one of the following:
\begin{itemize}
    \item $T_1\cap(L_1\cup\cdots\cup L_s)$, 
    \item $T_2\cap(L_1\cup\cdots\cup L_s)$, or
    \item $[p]$ for $p\in (L_1\cup\cdots\cup L_s)\setminus (T_1\cup T_2)$, $|[p]|>s$ 
and, if $p,q\in (L_1\cup\cdots\cup L_s)\setminus (T_1\cup T_2)$, then $[p]$ and $[q]$ are projectively equivalent.
\end{itemize}

\item[(b)] The group $G_i$ is a subgroup of the multiplicative group $\field^*$ hence abelian and,
for $p\in L_i\setminus (T_1\cup T_2)$, the action of 
$C_{\mathcal L}$ on $[p]$ is faithful (i.e., only the identity of $G_i$ takes $p$ to itself).

\item[(c)] 
Let $p\in L_i\setminus (T_1\cup T_2)$. 
If $[p]$ is finite, then $G_i$ is finite and cyclic.

\item[(d)] Each group $G_i$ is generated by elements of the form $f_{jil}f_{ijk}$,
hence $G_i$ is finite if and only if every element
$f_{jil}f_{ijk}$ has finite order. 

\end{enumerate}
\end{proposition}

\begin{proof}
As discussed in Remark \ref{StPosRem},
after a change of coordinates
we may assume $T_1,T_2$
are the standard transversals. 

(a) Note, for every $p\in L_1\cup\cdots\cup L_s$, that the orbit $[p]$ meets every line $L_i$. So, $|[p]|\ge s$.

By Corollary \ref{transCor},
$T_1\cap(L_1\cup\cdots\cup L_s)$ and
$T_2\cap(L_1\cup\cdots\cup L_s)$
are the only orbits with exactly $s$ elements. Every other orbit must be
$[p]$ for some $p\in (L_1\cup\cdots\cup L_s)\setminus (T_1\cup T_2)$. 

The map $f_e\colon \PP^3\to\PP^3$ defined for a nonzero $e\in \field$ by
$f_e((a:b:c:d))=(a:b:ec:ed)$ is the identity on $T_1$ and $T_2$ and hence
takes every line $L_i$ to itself. When $e\neq0,1$, the set of fixed points of $f_e$
is exactly $T_1\cup T_2$, so we see $f_e$ is then not the identity on any line $L_i$.
Moreover, if $f_e(p)=f_{e'}(p)$ for a point
$p\notin T_1\cup T_2$, then $e=e'$.
Thus as $e$ runs over all possible  
values of $\field$
different from 1, the images $f_e(p)$ of $p\in L_1\setminus(T_1\cup T_2)$ run over all 
points of $L_1\setminus(T_1\cup T_2)$. Therefore, 
for any two points $p,q\in L_1\setminus(T_1\cup T_2)$, there is an $e$ such that
$f_e(p)=q$. So for any two points $p,q\in (L_1\cup\cdots\cup L_s)\setminus(T_1\cup T_2)$, there are points
$p_1\in [p]\cap L_1$ and $q_1\in [q]\cap L_1$ and an $e$ with $f_e(p_1)=q_1$,
so projective equivalence of $[p]=[p_1]$ and $[q]=[q_1]$ follows
by Proposition \ref{FunctorCor}. 

(b) For $h\in G_i$, $T_2\cap L_i$ and $T_1\cap L_i$ are fixed points of $h$.
Regarding the fixed points as 0 and $\infty$ respectively
of $\PP^1$, $h$ is multiplication by
some element of $\field^*$ (i.e., some nonzero element of $\field$). Thus $G_i$
is isomorphic to a multiplicative subgroup of $\field^*$,
so $G_i$ is abelian.

Now say $g$ is an arrow in 
$C_{\mathcal L}$ with $g(p)=p$ for $p\in L_i\setminus(T_1\cup T_2)$. 
Then $g\in G_i$. Since $g(p)=p$ and $g$ is multiplication by some element of $\field^*$, 
this element must be 1 so $g$ is the identity.
Thus the action of $C_{\mathcal L}$ on $[p]$ 
is faithful. 

(c) We may assume $i=1$. Let $g\in G_1$.
Powers of $g$ applied to $p$
give only finitely many points, but this corresponds to multiplying
the element $c\in\field^*$ corresponding to the point $p$
by powers of the element $d\in\field^*$
corresponding to $g$, and so some power of $d$ is 1. 
Thus each $g$ has finite order,
so $G_1$, being finitely generated 
(by Proposition~\ref{Prop1}(a))
and abelian, is finite, and any finite multiplicative subgroup of a field is cyclic.

(d) By Proposition \ref{Prop1},
$G_i$ is generated by elements of the form $f_{jil}f_{ijk}$ and
$f_{kij}f_{jkl}f_{ijk}$. 
Let $i,j,k,l\in \{1,\ldots,s\}$ be distinct.

If $\mathcal{L'}=\{L_i,L_j,L_k,L_l\}$ is contained in a smooth quadric, then
$(G_\mathcal{L'})_i$ is the identity by Corollary \ref{transCor}, so 
$f_{jil}f_{ijk}$ and $f_{kij}f_{jkl}f_{ijk}$ are also the identity.

So assume $\mathcal{L'}=\{L_i,L_j,L_k,L_l\}$ is not contained in a smooth quadric.
After renumbering we may assume $i=1$ and
$\{L_j,L_k,L_l\}=\{L_2,L_3,L_4\}$, and after an appropriate
choice of coordinates we may
assume $L_2,L_3,L_4$ are in standard position with $T_1,T_2$ being the standard
transversals. We will regard
$T_2\cap L_1$ as being the point 0 on $L_1$ and
$T_1\cap L_1$ as being the point $\infty$ on $L_1$.

There are 12 possibilities for $f_{j1l}f_{1jk}$ with $j,k,l$ among 2, 3 and 4, namely
$f_{213}f_{123}$,
$f_{214}f_{124}$,
$f_{214}f_{123}$,
$f_{213}f_{124}$,
$f_{312}f_{132}$,
$f_{314}f_{134}$,
$f_{314}f_{132}$,
$f_{312}f_{134}$,
$f_{412}f_{142}$,
$f_{413}f_{143}$,
$f_{413}f_{142}$, and
$f_{413}f_{143}$.
But $f_{j1l}f_{1jk}$ is the identity 
whenever $k=l$, which
leaves six cases.

By direct computation we find that:

$f_{214}f_{123}$ is multiplication by $\frac{t-1}{t(1-l)}$;

$f_{213}f_{124}=f_{123}^{-1}f_{214}^{-1}=(f_{214}f_{123})^{-1}$ is multiplication by $\frac{t(1-l)}{t-1}$;

$f_{314}f_{132}$ is multiplication by $\frac{1}{lt}$;

$f_{312}f_{134}=f_{132}^{-1}f_{314}^{-1}=(f_{314}f_{132})^{-1}$ is multiplication by $lt$;

$f_{413}f_{142}$ is multiplication by $\frac{l-1}{l(1-t)}$; and

$f_{412}f_{143}=f_{142}^{-1}f_{413}^{-1}=(f_{413}f_{142})^{-1}$ is multiplication by $\frac{l(1-t)}{l-1}$.

There are six possibilities for $f_{k1j}f_{jkl}f_{1jk}$ which are not the identity
(since $j,k$ must be chosen from 2, 3 or 4, and if $l=1$ we get the identity).
Checking each one gives:

$f_{312}f_{234}f_{123}$ is multiplication by $\frac{l(t-1)}{1-l}$;

$f_{412}f_{243}f_{124}$ is multiplication by $lt$;

$f_{213}f_{324}f_{132}=f_{123}^{-1}f_{234}^{-1}f_{312}^{-1}=(f_{312}f_{234}f_{123})^{-1}$ is multiplication by $\frac{1-l}{l(t-1)}$;

$f_{413}f_{342}f_{134}$ is multiplication by $\frac{t(1-l)}{t-1}$;

$f_{214}f_{423}f_{142}=f_{124}^{-1}f_{243}^{-1}f_{412}^{-1}=(f_{412}f_{243}f_{124})^{-1}$ is multiplication by $\frac{1}{lt}$; and

$f_{314}f_{432}f_{143}=f_{134}^{-1}f_{342}^{-1}f_{413}^{-1}=(f_{413}f_{342}f_{134})^{-1}$ is multiplication by $\frac{t-1}{t(1-l)}$.

We see, as claimed, that the maps we get with three factors are the same as we get with two. 
\end{proof}

\begin{remark}\label{AlgorithmRem}
Working over an algebraically closed field $\field$,
our results allow us to find four lines with two distinct transversals whose group is a given finite cyclic group $G$.
Start with three lines $L_1,L_2,L_3$ in standard position and the two standard transversals, $T_1,T_2$.
Let $L_4\cap T_1$ be the point $(0:0:t:1)$ for $t\neq0,1$ and $L_4\cap T_2$ be $(l:1:0:0)$ for $l\neq0,1$.
If $|G|=1$, we must choose $l=1/t$ by Remark \ref{StPosRem} (since $G$ is trivial if and only if
all of the lines are in the same smooth quadric).

Say $|G|=m>2$. We saw in the proof of Proposition \ref{twoTrans} that $G_1$ is generated 
(as a multiplicative subgroup of $\field^*$) by
$\alpha=\frac{1}{lt}$, $\beta=\frac{t(1-l)}{t-1}$ and $\gamma=\frac{l(1-t)}{l-1}$.
Thus $\alpha\beta\gamma=1$, so in fact $G_1$ is generated by
$\alpha$ and $\beta$. But we can recover $(l,t)$ from $(\alpha,\beta)$, 
in particular 
$$t=\frac{\alpha\beta-1}{\alpha\beta-\alpha}, \hbox to.25in{\hfil}
l=\frac{1}{t\alpha}=\frac{\beta-1}{\alpha\beta-1}=\frac{\alpha\gamma-1}{\alpha\gamma-\alpha}.$$
The conditions $lt\neq1$, $l\neq1\neq t$ and $0\neq lt$
are collectively equivalent to $\alpha\beta\neq1$, $\alpha\neq1\neq \beta$ and
$0\neq\alpha\beta$. 
To get $G_1$ to be cyclic of order $m$, just pick
elements $\alpha,\beta\in K^*$ such that $\beta\neq 1/\alpha$ and
${\rm lcm}(|\alpha|,|\beta|)=m$. (Of course, if ${\rm char}(K)=p>0$,
this is possible if and only if $m$ is not divisible by $p$.)
\qed\end{remark}

To extend the previous
remark to more than 4 lines we recall the cross ratio
of four points on $\PP^1_{\overline{K}}$. The {\em cross ratio} of four points $(a_r:b_r)\in \PP^1_{\overline{K}}$, $r=1,2,3,4$,  is
$\chi_{1234}=\frac{(a_3b_1-a_1b_3)(a_4b_2-a_2b_4)}{(a_3b_2-a_2b_3)(a_4b_1-a_1b_4)}$; it is independent of the coordinates used
to represent the points (i.e., a linear change of coordinates does not change the cross ratio),
but if $\chi_{1234}=t$, then $t\in\overline{K}^*$, $t\neq 1$, and permuting the points gives the set
$\{t,\frac{1}{t},1-t,\frac{1}{1-t},\frac{t}{t-1},\frac{t-1}{t}\}$.
(See the end of the proof of Theorem \ref{Thm:HopfFibr} for an example of applying the next result to compute
$G_{\mathcal L}$.)

\begin{theorem}\label{CrossRatioRatios}
Let $\mathcal L=\{L_1'\ldots,L_s'\}$ be $s\geq4$ skew lines in $\PP^3_{\overline{\field}}$
with two distinct transversals $T_1', T_2'$. Let $q_{ij}=L_i'\cap T_j'$ and let
$\chi_{ij1kl}$ denote the cross ratios of the points, in order, 
$q_{il},q_{jl},q_{1l},q_{kl}$.
Then $G_1$, as a subgroup of $\overline{\field}^*$, is generated by the cross ratio ratios of the form 
$\chi_{ij1k1}/\chi_{ij1k2}$ for all distinct choices $q_{i1},q_{j1},q_{k1}$ and $q_{i2},q_{j2},q_{k2}$ 
of three of the points (other than $q_{1l}$) on each transversal.
\end{theorem}

\begin{proof}
By Proposition \ref{twoTrans}(d), $G_1$ is generated by elements of the form
$f_{i1k}f_{1ij}$. We may choose coordinates such that 
$L_1'=L_1, L_i'=L_3, L_j'=L_4, L_k'=L_2$
where $L_1,L_2,L_3$ are in standard position, $T_i'$ are the
standard transversals $T_i$, $L_4\cap T_1=(l:1:0:0)$ and $L_4\cap T_2=(0:0:t:1)$,
in which case (as in the proof of Proposition \ref{twoTrans}(d)), 
$f_{i1k}f_{1ij}(1:0:0:1)=(1:0:0:tl)$. 

Using just the $z$ and $w$ coordinates to parametrize $T_1$, the four points 
$q_{j1}=L_j'\cap T_1=L_4\cap T_1$, $q_{k1}=L_k'\cap T_1=L_2\cap T_1$, 
$q_{11}=L_1\cap T_1$ and $q_{i1}=L_i'\cap T_1=L_3\cap T_1$ become
$q_{j1}=(t:1)$, $q_{k1}=(1:1)$, $q_{11}=(0:1)$ and $q_{i1}=(1:0)$, so the cross ratio is
$\chi_{jk1i1}=\frac{(0-t)(1-0)}{(0-1)(1-0)}=t$.
Using just the $x$ and $y$ coordinates to parametrize $T_2$, the four points 
$q_{j2}=L_j'\cap T_2=L_4\cap T_2$, $q_{k2}=L_k'\cap T_2=L_2\cap T_2$, 
$q_{12}=L_1\cap T_2$ and $q_{i2}=L_i'\cap T_2=L_3\cap T_2$ become
$(l:1)$, $(1:1)$, $(1:0)$ and $(0:1)$, so the cross ratio is
$\chi_{jk1i2}=\frac{(1-0)(0-1)}{(1-0)(0-l)}=1/l$. Thus $f_{i1k}f_{1ij}\colon (a:0:0:b)\mapsto (a:0:0:btl)=(a:0:0:b\chi_{jk1i1}/\chi_{jk1i2})$
is multiplication by $\chi_{jk1i1}/\chi_{jk1i2}$.

Permuting the indices $i,j,k$ changes the cross ratios as follows: 

$\chi_{kj1i1}/\chi_{kj1i2}=(1/t)/l=\frac{1}{tl}$;

$\chi_{ik1j1}/\chi_{ik1j2}=(1-t)/(1-(1/l))=\frac{l(1-t)}{l-1}$;

$\chi_{ji1k1}/\chi_{ji1k2}=(t/(t-1))/((1/l)/((1/l)-1))=\frac{t(1-l)}{t-1}$;

$\chi_{ki1j1}/\chi_{ki1j2}=(1/(1-t))/(1/(1-(1/l)))=\frac{l-1}{l(1-t)}$;

$\chi_{ij1k1}/\chi_{ij1k2}=((t-1)/t)/(((1/l)-1)/(1/l))=\frac{t-1}{t(1-l)}$.

From the proof of Proposition \ref{twoTrans}(d), we see these cross ratio ratios
are exactly the generators of $G_1$. I.e., given a set of skew lines $L_1',\ldots,L_s'$
with two distinct transversals $T_1', T_2$, we get the points $q_{ij}=L_i'\cap T_j$.
Then $G_1$, as a subgroup of $\overline{\field}^*$, is generated by the cross ratio ratios of the form 
$\chi_{ij1k1}/\chi_{ij1k2}$ for all choices $q_{i1},q_{j1},q_{11},q_{k1}$ and $q_{i2},q_{j2},q_{12},q_{k2}$ 
of four of the points on each transversal.
\end{proof} 

As noted in Remark \ref{r.2 transversal-a},
it is also possible for ${\mathcal L}=\{L_1,\ldots,L_s\}\subset\PP^3_\field$, $s\geq4$, to have
a unique transversal. A particularly interesting case of this is when that transversal
in some sense counts double. We will say $T$ is a \emph{transversal of multiplicity 2} 
for ${\mathcal L}$ with respect to $\{i_1,i_2,i_3\}$ for three distinct indices $1\leq i_j\leq s$,
if we have: $T$ is transversal for ${\mathcal L}$; no single smooth quadric contains all of the lines $L_i$;
and the quadric $Q_{i_1i_2i_3}$ containing $L_{i_1},L_{i_2},L_{i_3}$ is such that
each line $L_i\in{\mathcal L}$ is either contained in $Q_{i_1i_2i_3}$ or is tangent to $Q_{i_1i_2i_3}$
at the point $L_i\cap T$. By the next lemma, 
$T$ being a \emph{transversal of multiplicity 2} for ${\mathcal L}$ is independent of 
the choice of the three indices, so we can simply say
$T$ is a \emph{transversal of multiplicity 2} for ${\mathcal L}$
if it is so with respect to any choice of $\{i_1,i_2,i_3\}$.

\begin{lemma}\label{oneTransLem}
Let ${\mathcal L}=\{L_1,\ldots,L_s\}\subset\PP^3_\field$, $s\geq4$, be skew.
Let $T$ be a transversal of multiplicity 2 for ${\mathcal L}$ with respect to $\{i_1,i_2,i_3\}$.
Then $T$ is a transversal of multiplicity 2 for ${\mathcal L}$ for every choice of three indices.
\end{lemma}

\begin{proof} 
Let $\{j_1,j_2,j_3\}$ be another choice of three indices.
If $\{j_1,j_2,j_3\}=\{i_1,i_2,i_3\}$, then clearly $T$ is a transversal of multiplicity 2 for ${\mathcal L}$
with respect to $\{j_1,j_2,j_3\}$. 

Now say $\{j_1,j_2,j_3\}$ and $\{i_1,i_2,i_3\}$ have exactly two indices in common, say $j_1=i_1, j_2=i_2$ 
but $j_3\neq i_3$. If $L_{j_3}\subset Q_{i_1i_2i_3}$, then $Q_{j_1j_2j_3}= Q_{i_1i_2i_3}$ and again
clearly $T$ is a transversal of multiplicity 2 for ${\mathcal L}$ with respect to $\{j_1,j_2,j_3\}$.
If however $L_{j_3}\not\subset Q_{i_1i_2i_3}$, then $Q_{j_1j_2j_3}\neq Q_{i_1i_2i_3}$, in which case
$Q_{j_1j_2j_3}\cap Q_{i_1i_2i_3}$ is a curve $C$ of degree 4 which as a divisor on either quadric is of type $(2,2)$
(meaning every line on the quadric meets $C$ twice). Since $L_{i_1}\cup L_{i_2}\cup T\subset Q_{j_1j_2j_3}\cap Q_{i_1i_2i_3}$,
we see as a divisor that $C=L_{i_1}+ L_{i_2}+T+L$ where $L$ is a line in the same ruling as $T$. Moreover, $L$
meets the four lines $L_{i_1}=L_{j_1}, L_{i_2}=L_{j_2}, L_{i_3}, L_{j_3}$, but $L_{i_3}$ meets $Q_{i_1i_2i_3}$
only at $T\cap L_{i_3}$ while $L_{j_3}$ meets $Q_{j_1j_2j_3}$ only at $T\cap L_{j_3}$, so $L$ is the line
through the points $T\cap L_{i_3}$ and $T\cap L_{j_3}$. I.e., $T=L$, so $Q_{j_1j_2j_3}\cap Q_{i_1i_2i_3}=L_{i_1}+ L_{i_2}+2T$.
In particular, $Q_{j_1j_2j_3}$ and $Q_{i_1i_2i_3}$ are tangent along $T$ and so every plane tangent to 
$Q_{i_1i_2i_3}$ at a point of $T$ is also tangent to $Q_{j_1j_2j_3}$ at that point.
Now recall that every line $L_j$ is either contained in $Q_{i_1i_2i_3}$ or tangent to $Q_{i_1i_2i_3}$ at $L_j\cap T$. 
Say $L_j$ is not contained in $Q_{j_1j_2j_3}$. If $L_j$ is tangent to $Q_{i_1i_2i_3}$ at $L_j\cap T$, then 
it is in the plane tangent to $Q_{j_1j_2j_3}$ at that point; i.e., $L_j$ is tangent to $Q_{j_1j_2j_3}$ at $L_j\cap T$.
Now say $L_j$ is contained in $Q_{i_1i_2i_3}$. Then the plane $A$ tangent to $Q_{i_1i_2i_3}$ at $L_j\cap T$
is spanned by $L_j$ and $T$. But $A\cap Q_{j_1j_2j_3}=L+T$ where $L$ is the line through 
the point $L_j\cap T$ in the ruling on $Q_{j_1j_2j_3}$ transverse to $T$. Thus 
$L_j\cap Q_{j_1j_2j_3} = (L_j\cap A)\cap Q_{j_1j_2j_3} = L_j\cap (A\cap Q_{j_1j_2j_3})=L_j\cap (L\cup T)$.
Since $L_j$ is not contained in $Q_{j_1j_2j_3}$, $L_j$ is neither $L$ nor $T$, and all three lines contain
the point $L_J\cap T$, so $L_j\cap (L\cup T)=L_j\cap T$; i.e., $L_j$ meets $Q_{j_1j_2j_3}$ in the single point $L_j\cap T$,
hence is tangent to $Q_{j_1j_2j_3}$ at $L_j\cap T$. Thus $T$ is a transversal of multiplicity 2 for ${\mathcal L}$
with respect to $\{j_1,j_2,j_3\}$. 

Next say $\{j_1,j_2,j_3\}$ and $\{i_1,i_2,i_3\}$ have exactly one index in common, say $j_1=i_1$ but $j_2\neq i_2$ 
and $j_3\neq i_3$. Given that $T$ is a transversal of multiplicity 2 for ${\mathcal L}$
with respect to $\{i_1=j_1,i_2,i_3\}$, we see from the preceding paragraph that $T$ is a transversal of multiplicity 2 for ${\mathcal L}$
with respect to $\{i_1=j_1,j_2,i_3\}$ and hence also for $\{i_1=j_1,j_2,j_3\}$.

Finally, if $\{j_1,j_2,j_3\}$ and $\{i_1,i_2,i_3\}$ have no indices in common, 
having it for $\{i_1,i_2,i_3\}$ implies we have it for $\{i_1,j_2,j_3\}$ and thus also
for $\{j_1,j_2,j_3\}$.
\end{proof} 

The case of sets ${\mathcal L}$ of 4 lines is key for understanding $G_{\mathcal L}$ in general,
since (by Proposition \ref{Prop1}(a)) generators of $G_{\mathcal L}$ come from subsets of four of the lines of ${\mathcal L}$.
By Remark \ref{r.2 transversal-a}, a set ${\mathcal L}$ of 4 disjoint lines can have either 1, 2 or infinitely many transversals.
When there are infinitely many, the lines are all contained in a smooth quadric so $G_{\mathcal L}$ is trivial
by Corollary \ref{transCor}. 
When there are exactly two then $G_{\mathcal L}$ is a subgroup of the multiplicative group $\field^*$ by Proposition \ref{twoTrans}(b).
The next result addresses the case of there being a unique transversal.
Note that ${\mathcal L}=\{L_1,L_2,L_3,L_4\}$ having a unique transversal $T$ means that $L_4$ is tangent to the quadric
containing $L_1,L_2,L_3$, and hence that $T$ is a transversal of multiplicity 2 for ${\mathcal L}$.
See Example \ref{Example:D4} for an explicit example of this.

\begin{proposition}\label{oneTransProp}
Consider a set of 4 distinct skew lines  ${\mathcal L}=\{L_1,\ldots,L_4\}\subset\PP^3_\field$.
Assume that $T$ is a transversal of multiplicity 2 for ${\mathcal L}$.
Let $q=(a:b:c:d)$ be the point $T\cap L_4$.
Then we have the following.
\begin{enumerate}
\item[(a)] $T$ (and hence $q$) is defined over $\field$ when ${\rm char}(\field)\neq 2$, or when ${\rm char}(\field)=2$ and $|\field|<\infty$.
\item[(b)] If $T$ is defined over $\field$, then $G_{\mathcal L}$ is a subgroup of the additive group $\field$, hence abelian and,
for $p\in L_i\setminus T$, the action of 
$C_{\mathcal L}$ on $[p]$ is faithful.
\item[(c)] Let $p\in \cup_i L_i\setminus T$. 
Then $G_{\mathcal L}$ and hence $[p]$ is infinite if and only if ${\rm char}(\field)=0$.
\end{enumerate}
\end{proposition}

\begin{proof} 
(a) As usual we may choose coordinates such that $L_1,L_2,L_3$ are in standard position
and $T_1,T_2$ are the standard transversals for these three lines. These are all contained in the
quadric $Q$ defined by $xz-yw$. 
The lines in the ruling on $Q$ transverse to $L_1,L_2,L_3$ are of the form
$vx-tw$, $vy-tz$ for scalars $t,v$, and since $q\in T$ we see for $T$ that $t=a,v=d$. 
(Note that $a$ and $d$ cannot both be 0, since then we would have $q\in L_3$ and hence $q\in L_3\cap L_4$.)
The lines in the same ruling as 
$L_1,L_2,L_3$ are, similarly, of the form $sx-uy$, $sw-uz$ for scalars $s,u$, so the line $L$ through $q$
in this ruling has $s=b,u=a$ or $s=c,u=d$ (one must take $s=c,u=d$ if $a=b=0$ and one must take
$s=b,u=a$ if $c=d=0$). 
The plane tangent to $Q$ at $q$ is the plane containing $T$ and $L$, hence is defined by
$cx-dy+az-bw$. The line $L_4$ meets $Q$ only at $q$ hence is tangent to $Q$ at $q$ so is contained in this plane.

If $d=0$, then $c=0$ (since $q\in Q$) and $T$ is defined by $z=w=0$ so $T=T_2$ is defined over $\field$. 

If $d=1$, then $q=(a:ac:c:1)$ and $c\neq0$ (since $q\not\in L_1$).
Now $T$ is $x-aw$, $y-az$, $L$ is $cx-y$, $cw-z$ and the tangent plane is
$cx-y+az-acw$ so $y=cx+az-acw$. 
Since $L_4$ is in the pencil of lines through $q$ in the tangent plane,
$L_4$ is defined by $cx-y+az-acw$ and some linear combination
of $x-aw$ and $z-cw$, but $L_4$ is not in $Q$ so is not either $T$ nor $L$,
so the linear combination can be taken to be $(x-aw)+\alpha(z-aw)$
(or $x+\alpha z-(a+\alpha c)w$) for some $\alpha\neq0$.
We now show that $T$ is defined over $\field$. 

Since $L_4$ is defined over $\field$,
some linearly independent linear combinations $\beta(cx-y+az-acw)+\delta(x+\alpha z-(a+\alpha c)w)$ and 
$\beta'(cx-y+az-acw)+\delta'(x+\alpha z-(a+\alpha c)w)$
have coefficients only in $\field$. Looking at the $y$ terms shows $\beta,\beta'\in\field$.
If $\beta,\beta'\neq0$ we can (by dividing through by $\beta$ and $\beta'$) reduce to the case that
$\beta=\beta'=1$. Subtracting $(cx-y+az-acw)+\delta'(x+\alpha z-(a+\alpha c)w)$ from
$(cx-y+az-acw)+\delta(x+\alpha z-(a+\alpha c)w)$ reduces to the case that $\beta=0$ and $\beta'=1$.
If $\beta=0$ then $\beta'\neq0$ and we again reduce to the case that $\beta=0$ and $\beta'=1$.
If $\beta'=0$ then $\beta\neq0$ and we reduce to the case that $\beta'=0$ and $\beta=1$,
which is symmetric to the case $\beta=0$ and $\beta'=1$.

So assume $\beta=0$ and $\beta'=1$, hence
$\beta(cx-y+az-acw)+\delta(x+\alpha z-(a+\alpha c)w)$ is $\delta(x+\alpha z-(a+\alpha c)w)$ and
$\beta'(cx-y+az-acw)+\delta'(x+\alpha z-(a+\alpha c)w)$ is $(\delta'+c)x-y+(a+\delta'\alpha)z-(ac+a\delta'+\delta'\alpha c)w$.
These are defined over $\field$
so $\delta\in\field$ hence $\alpha,a+\alpha c\in \field$ and
$\delta'+c,a+\delta'\alpha\in\field$, hence $\alpha(\delta'-c)=a+\delta'\alpha-(a+\alpha c)\in\field$,
so $\delta'+c,\delta'-c\in\field$. Thus $2\delta',2c\in\field$, so if ${\rm char}(\field)\neq2$ we have $\delta',c\in\field$.
From $\delta',\alpha,a+\delta'\alpha\in\field$ we now get $a\in\field$. Since $a\in\field$ we see
$T$ is defined over $\field$ when ${\rm char}(\field)\neq 2$.

Now assume ${\rm char}(\field)=2$ but $\field$ is finite. As above we reduce to the case that 
$\beta=0,\beta'=1$ and as before $\delta,\alpha,a+\alpha c,\delta'+c,a+\delta'\alpha\in \field$,
and also $ac+a\delta'+\delta'\alpha c\in\field$. Thus
$(\delta')^2\alpha=(ac+a\delta'+\delta'\alpha c)+(\delta'+c)(a+\delta'\alpha)\in\field$
hence $(\delta')^2\in\field$. But the map $f\colon \field\to\field$ given by $f(\lambda)=\lambda^2$
is injective since the characteristic is 2 hence surjective since $\field$ is finite, so 
$\delta'\in\field$. Now from $\alpha,\delta',a+\delta'\alpha\in\field$ we get $a\in\field$
(and hence $c\in\field$), so again $T$ is defined over $\field$.

(b) Choose coordinates such that $L_1,L_2,L_3$ are in standard position
and $T=T_1,T_2$ are the standard transversals for these three lines. These are all contained in the
quadric $Q$ defined by $xz-yw$. We will determine the group $G_1$ as a group of matrices
with respect to this choice of coordinates, thinking of a point $(a:0:0:b)\in L_1$ as a point $(a:b)\in\PP^1_\field$.
We regard the point $(0:0:0:1)$ given by $T\cap L_1$ as $\infty$; the remaining points of $L_1$ are of the form
$(1:0:0:b)$. We can assume $L_4$ meets $T$ at $(0:0:t:1)$ for some $t\in \field$, $t\neq0,1$,
and is defined by $ry+z-tw$, $r\neq0$, and $y-tx$.

As in the proof of Proposition \ref{twoTrans}(d), we explicitly determine the generators of $G_1$,
using Proposition \ref{Prop1}(a). Since $T$ is a transversal for all four lines, every element of $G_1$
maps $(0:0:0:1)$ to itself. One can determine the matrix representing an element of $G_1$
by computing the images also of $(1:0:0:0)$ and $(1:0:0:1)$. 

We now list the generators and their matrices (found by direct computation): 

$(f_{213}f_{124})^{-1}=f_{214}f_{123}$: 
$\displaystyle\begin{pmatrix}
1 & 0 \\
\frac{r}{1-t} & 1 \\
\end{pmatrix}$;
$(f_{312}f_{134})^{-1}=f_{314}f_{132}$: 
$\displaystyle\begin{pmatrix}
1 & 0 \\
r  & 1 \\
\end{pmatrix}$;

$(f_{412}f_{143})^{-1}=f_{413}f_{142}$: 
$\displaystyle\begin{pmatrix}
1 & 0 \\
\frac{rt}{t-1}  & 1 \\
\end{pmatrix}$;

$(f_{213}f_{324}f_{132})^{-1}=f_{312}f_{234}f_{123}$: 
$\displaystyle\begin{pmatrix}
1 & 0 \\
\frac{rt}{1-t}  & 1 \\
\end{pmatrix}$;
$(f_{214}f_{423}f_{142})^{-1}=f_{412}f_{243}f_{124}$: 
$\displaystyle\begin{pmatrix}
1 & 0 \\
-r  & 1 \\
\end{pmatrix}$;

$(f_{314}f_{432}f_{143})^{-1}=f_{413}f_{342}f_{134}$: 
$\displaystyle\begin{pmatrix}
1 & 0 \\
\frac{r}{t-1}  & 1 \\
\end{pmatrix}$.

Notice that $\begin{pmatrix}
1 & 0 \\
u  & 1 \\
\end{pmatrix}$ fixes the point $(0:1)$ (i.e., $\infty$) and
translates the point $(1:v)$ to $(1:v+u)$. In particular, these
maps all are translations hence $G_1$ can be regarded as
the additive subgroup of $\field$ generated by
$r, \frac{1}{t-1}, \frac{r}{t-1}, \frac{rt}{t-1}$. Thus $G_1$ acts 
faithfully on $L_1$ away from $T\cap L_1$, hence on $\in L_1\cup \cdots\cup L_4\setminus T$,
$C_{\mathcal L}$ also acts faithfully.

(c) By Proposition \ref{Prop1}(e), we can assume $T$ is defined over $\field$,
hence $G_{\mathcal L}$ is a nontrivial (since $Q$ does not contain $L_4$)
finitely generated additive subgroup of $\field$ by (b). This is infinite if ${\rm char}(\field)=0$
(since any nonzero element of $\field$ has infinite additive order)
and finite if ${\rm char}(\field)>0$
(since every element of $\field$ has finite additive order).
\end{proof}

\begin{corollary}\label{oneTransCor}
Let ${\mathcal L}=\{L_1,\ldots,L_s\}\subset\PP^3_\field$, $s\geq4$, be skew.
Assume there is a transversal  $T$ of multiplicity 2 for ${\mathcal L}$.
Then $G_{\mathcal L}$ is a subgroup of the additive group $\overline{\field}$, hence abelian.
In addition, $T\cap (L_1\cup\cdots\cup L_s)$ is the only orbit with $s$ elements; indeed,
for $p,q\in \cup_i L_i\setminus T$, the group $G_{\mathcal L}$ acts 
faithfully on each line so we have $s|G_{\mathcal L}|=|[p]|=|[q]|>s$
and in fact $[p]$ and $[q]$ are projectively equivalent.
Moreover, for $p\in \cup_i L_i\setminus T$, $[p]$ is infinite if and only if ${\rm char}(\field)=0$.
\end{corollary}

\begin{proof} 
By Proposition \ref{Prop1}(e), we can compute $G_{\mathcal L}$ working over the algebraic closure
$\overline{\field}$, where we can apply Proposition \ref{oneTransProp}(b),
showing that $G_\ell$ is a subgroup of $\overline{\field}$ for every subset $\ell\subseteq {\mathcal L}$
of 4 of the lines. But by Proposition \ref{Prop1}(a), the union of these $G_\ell$ contain generators for
$G_{\mathcal L}$, hence $G_{\mathcal L}$ is also a subgroup of $\overline{\field}$.
Since the lines in $\mathcal L$ are not all contained in a single quadric,
$G_{\mathcal L}$ has nontrivial elements, and any nonzero element of $\overline{\field}$
has infinite order when ${\rm char}(\field)=0$. 
However, when ${\rm char}(\field)=p>0$, $G_{\mathcal L}$ is a finitely generated subgroup of $\overline{\field}$,
hence a finite dimensional vector space over $\mathbb Z/p\mathbb Z$, hence finite,
acting as a group of translations on each line of ${\mathcal L}$ away from a single fixed point on that line,
hence $G_{\mathcal L}$ acts faithfully on line.

Note, for every $p\in L_1\cup\cdots\cup L_s$, that the orbit $[p]$ meets every line $L_i$. So, $|[p]|\geq s$.
By Corollary \ref{transCor}, $T\cap(L_1\cup\cdots\cup L_s)$ 
is an orbit with exactly $s$ elements and any other orbit with $s$ elements
must also lie on a transversal. But since $T$ has multiplicity 2, there is a 
quadric $Q$ containing at least 3 (but not all) of the lines $L_i$ and every line
not contained in $Q$ is tangent to $Q$. Every transversal for $\mathcal L$
is contained in $Q$, but if $L_i$ is one of the lines tangent to but not contained in $Q$,
then $L_i$ meets $Q$ in only one point (which is on $T$) so there can be no 
transversal other than $T$, hence no other orbits with $s$ elements.
Thus for $p,q\in \cup_i L_i\setminus T$ we have $|[p]|=|[q]|>s$.

To see that $[p]$ and $[q]$ are projectively equivalent,
we can by Proposition \ref{Prop1}(e) assume $\field$ is
algebraically closed. Also, we can after renumbering assume that $L_4$ is not on
the quadric $Q_{123}$ containing $L_1,L_2,L_3$ and hence
$L_4$ is tangent to $Q_{123}$ at a point of $T$.
Since $\field=\overline{\field}$, after a change of coordinates
we can assume that $L_1,L_2,L_3$ are in standard position
with $T=T_1$ and $T_2$ the standard transversals for $L_1,L_2,L_3$.
Thus $Q_{123}$ is defined by $xz-yw$.
Moreover, $L_4$ is tangent to $Q_{123}$ at a point
$v\in T$ but $Q_{123}$ does not contain $L_4$.
The points of $T$ are of the form $(0:0:a:b)$,
but $v$, being on $L_4$, cannot be on $L_3$ (so $b\neq0$,
so we may assume $b=1$), nor on $L_1$
(so $a\neq0$), nor on $L_2$, so $v=(0:0:a:1)$ for some
$a\neq0,1$. The tangent plane to $Q_{123}$ at $v$
is spanned by the two ruling lines in $Q_{123}$
through $v$, namely $T$ (defined by $x,y$) and
a line $L_v$ defined by $z=aw, y-ax$. The plane $\Gamma$ tangent to
$Q_{123}$ at $v$ is thus $y-ax$ and the lines in $\Gamma$
through $v$ are thus defined by $y=ax, z=aw-ux$ for some  
scalar $u$. In particular, $L_4$ is defined by
$y=ax, z=aw-ux$ for some $u\neq0$ (since
$L_4$ is not $L_v$ nor $T$, both of which
are contained in $Q_{123}$).

Consider the map $\Phi_g\in PGL_4(\field)$ given by the matrix 
$\begin{pmatrix}
1 & 0  & 0 & 0 \\
0 & 1  & 0 & 0 \\
0 & g  & 1 & 0 \\
g & 0  & 0 & 1 \\
\end{pmatrix}$.
One checks directly that $\Phi_g(L_i)=L_i$ for $i=1,2,3,4$ 
for every value of $g$. Moreover, $\Phi_g$ is the identity on $T$. 
As $g$ runs through $\field$, the image of $(1:0:0:1)\in L_1$
runs over every point of $L_1$ off $T$. In fact,  
$\Phi_g((1:0:0:1))=(1:0:0:1+g)$, so $\Phi_g$ restricted to $L_1$
is translation by $g$. Thus for any $p,q\in \cup_i L_i\setminus T$,
there are points of $p_1,q_1\in L_1$ with $[p_1]=[p]$ and $[q_1]=[q]$,
and a translation $\Phi_g$ taking $p_1$ to $q_1$.
Thus projective equivalence of $[p]=[p_1]$ and $[q]=[q_1]$ follows
by Proposition \ref{FunctorCor}.
\end{proof} 

We can now characterize commutativity of $G_{\mathcal L}$ in terms of transversals.

\begin{theorem}\label{TransThm}
Let ${\mathcal L}=\{L_1,\ldots,L_s\}\subset\PP^3_\field$, $s\geq3$, be skew.
Then $G_{\mathcal L}$ is abelian if and only if ${\mathcal L}$
has either 2 or more distinct transversals or it has a transversal of multiplicity 2.
\end{theorem}

\begin{proof} First say $s=3$. Then the result follows because
$|G_{\mathcal L}|=0$ and the lines are in the same ruling of a smooth quadric
and hence there are infinitely many transversals. So say $s>3$.

By Proposition \ref{twoTrans}, $G_{\mathcal L}$ is abelian if 
${\mathcal L}$ has 2 or more transversals, 
and by Corollary \ref{oneTransCor}, it is abelian if 
${\mathcal L}$ has a transversal of multiplicity 2.

Consider the converse, which we prove contrapositively. First suppose ${\mathcal L}$ has a unique transversal, $T$, but that it is not of multiplicity 2.
Then there are 4 distinct lines $L_i, L_j, L_k, L_l\in{\mathcal L}$ such that
$L_i$ is not contained in and not tangent to $Q_{jkl}$, where $Q_{jkl}$
is the quadric containing $L_j, L_k, L_l$. Since $L_i$ meets $Q_{jkl}$ in two points, 
there is a second transversal $T'$ for $L_i, L_j, L_k, L_l$, but $T'$ is not transversal for
${\mathcal L}$ so there is a line $L_m\in {\mathcal L}$ not meeting $T'$.

There are two cases: (case A) $L_m$ is or (case B) is not tangent to $Q_{jkl}$.

Case A: Say it is tangent (indicated by the black dot in Figure \ref{FigA}). 
Then the group for the lines $L_j,L_k,L_l,L_m$ give a nontrivial
element of $G_j$ which by Proposition \ref{oneTransProp}(b) we can regard as a 
translation $x\mapsto x+a$, $a\neq0$, on $\overline{\field}$, while the lines $L_i,L_j,L_k,L_l$ give a nontrivial
element of $G_j$ which by Proposition \ref{twoTrans}(b) we can regard as a
scaling $x\mapsto bx$, $b\neq 0,1$. These do not commute, so $G_j\cong G_{\mathcal L}$
is nonabelian.

\begin{figure}[h]
\begin{tikzpicture}[x=1cm,y=1cm]
\clip(-11.09,4) rectangle (0,6.5);
\draw [line width=1pt] (-7,6)-- (-7,4);
\draw [line width=1pt] (-6,6)-- (-6,4);
\draw [line width=1pt] (-5,6)-- (-5,4);
\draw [line width=1pt] (-4,6)-- (-4,4);
\draw [line width=1pt] (-8,5.44)-- (-2,5.44);
\draw [line width=1pt] (-3,6)-- (-3,4);
\draw [fill=white,color=white] (-3,4.58) circle (3.5pt);
\draw [line width=1pt] (-7.96,4.58)-- (-1.96,4.58);
\draw[color=black] (-6.97,6.32) node {$L_j$};
\draw[color=black] (-6.03,6.32) node {$L_k$};
\draw[color=black] (-5.01,6.32) node {$L_l$};
\draw[color=black] (-3.95,6.32) node {$L_i$};
\draw[color=black] (-2.97,6.32) node {$L_m$};
\draw[color=black] (-8.2,5.5) node {$T$};
\draw[color=black] (-8.2,4.6) node {$T'$};
\draw [fill=black,color=black] (-3,5.44) circle (3.5pt);
\end{tikzpicture}
\caption{Case A}
\label{FigA}
\end{figure}

Case B: Now say it is not tangent (see Figure \ref{FigB}). Then $L_m$ meets $Q_{jkl}$ in a second point, hence there
is a line $T''$ transversal to $L_j,L_k,L_l,L_m$ but disjoint from $L_i$ (since
$L_i\not\subset Q_{jkl}$). We have three points of $L_j$, namely 
$T\cap L_j$, $T'\cap L_j$ and $T''\cap L_j$. By Proposition \ref{twoTrans}(b),
we have a nontrivial element $f\in G_j$ coming from the lines
$L_i,L_j,L_k,L_l$ (this fixes the points $T\cap L_j$, $T'\cap L_j$) and another
nontrivial element $g\in G_j$ coming from the lines
$L_j,L_k,L_l,L_m$ (this fixes the points $T\cap L_j$, $T''\cap L_j$).
Choose coordinates on $L_j$ such that $T\cap L_j$ is $(0:1)$,
$T'\cap L_j$ is $(1:0)$ and $T''\cap L_j$ is $(1:1)$.
Then $f$ is represented for some $b\neq0,1$ by the matrix
$\begin{pmatrix}
b & 0 \\
0 & 1\\
\end{pmatrix}$ while 
the matrix for $g$ is
$\begin{pmatrix}
c & 0 \\
c-1 & 1\\
\end{pmatrix}$ for some $c\neq0,1$.
These do not commute.
This finishes the case that ${\mathcal L}$ has a unique transversal, $T$.

\begin{figure}[h]
\begin{tikzpicture}[x=1cm,y=1cm]
\clip(-11.09,3) rectangle (0,6.5);
\draw [line width=1pt] (-7,6)-- (-7,3);
\draw [line width=1pt] (-6,6)-- (-6,3);
\draw [line width=1pt] (-5,6)-- (-5,3);
\draw [line width=1pt] (-4,6)-- (-4,3);
\draw [fill=white,color=white] (-4,3.77) circle (3.5pt);
\draw [line width=1pt] (-8,5.44)-- (-2,5.44);
\draw [line width=1pt] (-3,6)-- (-3,3);
\draw [fill=white,color=white] (-3,4.58) circle (3.5pt);
\draw [line width=1pt] (-7.96,4.58)-- (-1.96,4.58);
\draw [line width=1pt] (-7.99,3.77)-- (-1.95,3.77);
\draw[color=black] (-6.97,6.32) node {$L_j$};
\draw[color=black] (-6.03,6.32) node {$L_k$};
\draw[color=black] (-5.01,6.32) node {$L_l$};
\draw[color=black] (-3.95,6.32) node {$L_i$};
\draw[color=black] (-2.97,6.32) node {$L_m$};
\draw[color=black] (-8.2,5.5) node {$T$};
\draw[color=black] (-8.2,4.6) node {$T'$};
\draw[color=black] (-8.2,3.8) node {$T''$};
\end{tikzpicture}
\caption{Case B}
\label{FigB}
\end{figure}

Now assume ${\mathcal L}$ has no transversals.
In particular, no smooth quadric contains all of the lines of ${\mathcal L}$.
Thus there are four lines, say $L_1,L_2,L_3,L_4$ after renumbering, not on any smooth quadric.
Let $Q_{123}$ be the quadric containing $L_1,L_2,L_3$.

First assume that no four of the lines in ${\mathcal L}$ have two distinct transversals (case C).
Thus $L_1,L_2,L_3,L_4$ has a unique transversal $T_1$ and $L_4$ is tangent to 
(as indicated  by the black dot in Figure \ref{FigC}) but not contained in $Q_{123}$. Since ${\mathcal L}$ has no transversals, there must be a line,
say $L_5\in {\mathcal L}$, that does not meet $T_1$. Since no four of the lines in ${\mathcal L}$ 
have two distinct transversals, the lines $L_1,L_2,L_3,L_5$ has a unique transversal $T_2$
so $L_5$ is also tangent to $Q_{123}$ (again indicated by the black dot). Choose coordinates on $L_1$ such that
$T_1\cap L_1$ is $(0:1)$ and $T_2\cap L_1$ is $(1:0)$. 
Then as in the proof of Proposition \ref{oneTransProp}(b),
the elements of $G_1$ coming from $L_1,L_2,L_3,L_4$ have a nontrivial element
$\begin{pmatrix}
1 & 0 \\
u  & 1 \\
\end{pmatrix}$
and the elements of $G_1$ coming from $L_1,L_2,L_3,L_5$ have a nontrivial element
$\begin{pmatrix}
1 & v \\
0  & 1 \\
\end{pmatrix}$, but these do not commute since $u,v\neq0$, hence $G_1$ is not abelian.

\begin{figure}[h]
\begin{tikzpicture}[x=1cm,y=1cm]
\clip(-11.09,3) rectangle (0,6.5);
\draw [line width=1pt] (-3,6)-- (-3,3);
\draw [fill=white,color=white] (-3,5.44) circle (3.5pt);
\draw [line width=1pt] (-7,6)-- (-7,3);
\draw [line width=1pt] (-6,6)-- (-6,3);
\draw [line width=1pt] (-5,6)-- (-5,3);
\draw [line width=1pt] (-4,6)-- (-4,3);
\draw [fill=white,color=white] (-4,4.58) circle (3.5pt);
\draw [fill=black,color=black] (-4,5.44) circle (3.5pt);
\draw [fill=black,color=black] (-3,4.58) circle (3.5pt);
\draw [line width=1pt] (-8,5.44)-- (-2,5.44);
\draw [line width=1pt] (-7.96,4.58)-- (-1.96,4.58);
\draw [line width=1pt] (-7.99,3.77)-- (-1.95,3.77);
\draw[color=black] (-6.97,6.32) node {$L_1$};
\draw[color=black] (-6.03,6.32) node {$L_2$};
\draw[color=black] (-5.01,6.32) node {$L_3$};
\draw[color=black] (-3.95,6.32) node {$L_4$};
\draw[color=black] (-2.97,6.32) node {$L_5$};
\draw[color=black] (-8.2,5.5) node {$T_1$};
\draw[color=black] (-8.2,4.6) node {$T'$};
\draw[color=black] (-8.2,3.8) node {$T''$};
\end{tikzpicture}
\caption{Case C}
\label{FigC}
\end{figure}

Now assume $L_1,L_2,L_3,L_4\in {\mathcal L}$ have two distinct transversals (case D), $T_1,T_2$.
Since ${\mathcal L}$ has no transversals, either (case (D1)) there is a line $L_5$ that misses both $T_i$, 
or (case (D2)) there are two lines $L_5,L_6$ where $L_5$ meets $T_1$ but not $T_2$
and $L_6$ meets $T_2$ but not $T_1$.

Case (D1): Either (case (i); see Figure \ref{FigD1i}) $L_1,L_2,L_3, L_5$ has two transversals, $S_1,S_2$ or (case (ii)) just one, $S_1$.
We may choose coordinates on $L_1$ such that $T_1$ and $T_2$ meet $L_1$
at $(0:1)$ and $(1:0)$, and $S_1\cap L_1$ is $(1:1)$. In case (i), $S_2\cap L_1$
is $(t:1)$ for some $t\neq0,1$. In case (i) $G_1$ has an element
$\begin{pmatrix}
c & 0 \\
0  & 1 \\
\end{pmatrix}$ with $c\neq0,1$ coming from $L_1,L_2,L_3,L_4$
and an element 
$\begin{pmatrix}
at-1 & t-at \\
a-1  & t-a \\
\end{pmatrix}$ with $c\neq0,1$ coming from $L_1,L_2,L_3,L_5$
(this is the matrix with eigenvectors $(t:1)$, whose eigenvalue is $a$
so we need $a\neq0$ to be invertible, and $(1:1)$, whose eigenvalue we may assume is 1,
since we are really working modulo scalar matrices).
These matrices do not commute so $G_1$ is nonabelian.

\begin{figure}[h]
\begin{tikzpicture}[x=1cm,y=1cm]
\clip(-11.09,2) rectangle (0,6.5);
\draw [line width=1pt] (-3,6)-- (-3,2.3);
\draw [fill=white,color=white] (-3,5.44) circle (3.5pt);
\draw [fill=white,color=white] (-3,4.58) circle (3.5pt);
\draw [line width=1pt] (-7,6)-- (-7,2.3);
\draw [line width=1pt] (-6,6)-- (-6,2.3);
\draw [line width=1pt] (-5,6)-- (-5,2.3);
\draw [line width=1pt] (-4,6)-- (-4,2.3);
\draw [fill=white,color=white] (-4,3.77) circle (3.5pt);
\draw [fill=white,color=white] (-4,2.87) circle (3.5pt);
\draw [line width=1pt] (-8,5.44)-- (-2,5.44);
\draw [line width=1pt] (-7.96,4.58)-- (-1.96,4.58);
\draw [line width=1pt] (-7.99,3.77)-- (-1.95,3.77);
\draw [line width=1pt] (-7.99,2.87)-- (-1.95,2.87);
\draw[color=black] (-6.97,6.32) node {$L_1$};
\draw[color=black] (-6.03,6.32) node {$L_2$};
\draw[color=black] (-5.01,6.32) node {$L_3$};
\draw[color=black] (-3.95,6.32) node {$L_4$};
\draw[color=black] (-2.97,6.32) node {$L_5$};
\draw[color=black] (-8.2,5.5) node {$T_1$};
\draw[color=black] (-8.2,4.6) node {$T_2$};
\draw[color=black] (-8.2,3.8) node {$S_1$};
\draw[color=black] (-8.2,2.87) node {$S_2$};
\end{tikzpicture}
\caption{Case D1i}
\label{FigD1i}
\end{figure}

In case (ii) (see Figure \ref{FigD1ii}), the first matrix stays the same but the second matrix 
has $(1:0)$ and $(1:1)$ as eigenvectors with eigenvalue 1 so it is
$\begin{pmatrix}
1 & a \\
0  & a+1 \\
\end{pmatrix}$ for some $a\neq0$ (since the matrix is nontrivial).
Again these do not commute so $G_1$ is nonabelian.

\begin{figure}[h]
\begin{tikzpicture}[x=1cm,y=1cm]
\clip(-11.09,3) rectangle (0,6.5);
\draw [line width=1pt] (-3,6)-- (-3,3);
\draw [fill=white,color=white] (-3,5.44) circle (3.5pt);
\draw [fill=white,color=white] (-3,4.58) circle (3.5pt);
\draw [line width=1pt] (-7,6)-- (-7,3);
\draw [line width=1pt] (-6,6)-- (-6,3);
\draw [line width=1pt] (-5,6)-- (-5,3);
\draw [line width=1pt] (-4,6)-- (-4,3);
\draw [fill=white,color=white] (-4,3.77) circle (3.5pt);
\draw [line width=1pt] (-8,5.44)-- (-2,5.44);
\draw [line width=1pt] (-7.96,4.58)-- (-1.96,4.58);
\draw [line width=1pt] (-7.99,3.77)-- (-1.95,3.77);
\draw[color=black] (-6.97,6.32) node {$L_1$};
\draw[color=black] (-6.03,6.32) node {$L_2$};
\draw[color=black] (-5.01,6.32) node {$L_3$};
\draw[color=black] (-3.95,6.32) node {$L_4$};
\draw[color=black] (-2.97,6.32) node {$L_5$};
\draw[color=black] (-8.2,5.5) node {$T_1$};
\draw[color=black] (-8.2,4.6) node {$T_2$};
\draw [fill=black,color=black] (-3,3.77) circle (3.5pt);
\draw[color=black] (-8.2,3.8) node {$S_1$};
\end{tikzpicture}
\caption{Case D1ii}
\label{FigD1ii}
\end{figure}

Case (D2): Here we have three subcases.
Case (i): Neither of $L_5,L_6$ are tangent to $Q_{123}$.
Case (ii): Only one (say $L_5$) is tangent to $Q_{123}$.
Case (iii): Both are tangent to $Q_{123}$.

Case 2(i): Here $L_1,L_2,L_3,L_5$ has a second transversal, $T_3$, which might or might not
meet $L_6$ (indicated by an open circle in Figure \ref{FigD2i}).
The lines $L_1,L_2,L_3,L_4,L_5,T_1,T_2,T_3$ now give the same situation
as, respectively, $L_j,L_k,L_l,L_i,T,T',T''$ do in case B (compare Figures \ref{FigB}
and \ref{FigD2i}, ignoring $L_6$).

\begin{figure}[h]
\begin{tikzpicture}[x=1cm,y=1cm]
\clip(-11.09,3) rectangle (0,6.5);
\draw [line width=1pt] (-3,6)-- (-3,3);
\draw [fill=white,color=white] (-3,4.58) circle (3.5pt);
\draw [line width=1pt] (-2,6)-- (-2,3);
\draw [fill=white,color=white] (-2,5.44) circle (3.5pt);
\draw [color=black, line width=1] (-2,3.77) circle (3.5pt);
\draw [line width=1pt] (-7,6)-- (-7,3);
\draw [line width=1pt] (-6,6)-- (-6,3);
\draw [line width=1pt] (-5,6)-- (-5,3);
\draw [line width=1pt] (-4,6)-- (-4,3);
\draw [fill=white,color=white] (-4,3.77) circle (3.5pt);
\draw [line width=1pt] (-8,5.44)-- (-1,5.44);
\draw [line width=1pt] (-7.96,4.58)-- (-1,4.58);
\draw [line width=1pt] (-7.99,3.77)-- (-1,3.77);
\draw[color=black] (-6.97,6.32) node {$L_1$};
\draw[color=black] (-6.03,6.32) node {$L_2$};
\draw[color=black] (-5.01,6.32) node {$L_3$};
\draw[color=black] (-3.95,6.32) node {$L_4$};
\draw[color=black] (-2.97,6.32) node {$L_5$};
\draw[color=black] (-2,6.32) node {$L_6$};
\draw[color=black] (-8.2,5.5) node {$T_1$};
\draw[color=black] (-8.2,4.6) node {$T_2$};
\draw[color=black] (-8.2,3.8) node {$T_3$};
\end{tikzpicture}
\caption{Case D2i}
\label{FigD2i}
\end{figure}

Case 2(ii): Here $T_1$ is transversal for $L_1,L_2,L_3,L_4,L_5$ (but not $L_6$)
with $L_5$ tangent to $Q_{123}$, and $T_2$ is transversal to $L_1,L_2,L_3,L_4,L_6$
 (but not $L_5$) and $L_6$ is not tangent to $Q_{123}$ (see Figure \ref{FigD2ii}).
 In this case $L_1,L_2,L_3,L_4,L_5,T_1,T_2$
give the same situation as do $L_j,L_k,L_l,L_i,L_m,T,T'$ in case A (compare Figures \ref{FigA}
and \ref{FigD2ii}, ignoring $L_6$).

\begin{figure}[h]
\begin{tikzpicture}[x=1cm,y=1cm]
\clip(-11.09,4) rectangle (0,6.5);
\draw [line width=1pt] (-3,6)-- (-3,4);
\draw [fill=white,color=white] (-3,4.58) circle (3.5pt);
\draw [fill=black,color=black] (-3,5.44) circle (3.5pt);
\draw [line width=1pt] (-2,6)-- (-2,4);
\draw [fill=white,color=white] (-2,5.44) circle (3.5pt);
\draw [line width=1pt] (-7,6)-- (-7,4);
\draw [line width=1pt] (-6,6)-- (-6,4);
\draw [line width=1pt] (-5,6)-- (-5,4);
\draw [line width=1pt] (-4,6)-- (-4,4);
\draw [line width=1pt] (-8,5.44)-- (-1,5.44);
\draw [line width=1pt] (-7.96,4.58)-- (-1,4.58);
\draw[color=black] (-6.97,6.32) node {$L_1$};
\draw[color=black] (-6.03,6.32) node {$L_2$};
\draw[color=black] (-5.01,6.32) node {$L_3$};
\draw[color=black] (-3.95,6.32) node {$L_4$};
\draw[color=black] (-2.97,6.32) node {$L_5$};
\draw[color=black] (-2,6.32) node {$L_6$};
\draw[color=black] (-8.2,5.5) node {$T_1$};
\draw[color=black] (-8.2,4.6) node {$T_2$};
\end{tikzpicture}
\caption{Case D2ii}
\label{FigD2ii}
\end{figure}

Case 2(iii): Here $T_1$ is transversal for $L_1,L_2,L_3,L_4,L_5$ (but not $L_6$)
with $L_5$ tangent to $Q_{123}$, and $T_2$ is transversal to $L_1,L_2,L_3,L_4,L_6$
 (but not $L_5$) and $L_6$ is tangent to $Q_{123}$.
 In this case $L_1,L_2,L_3,L_4,L_5,T_1,T_2$ again
gives the same situation as do $L_j,L_k,L_l,L_i,L_m,T,T'$ in case A  (compare Figures \ref{FigA}
and \ref{FigD2iii}, ignoring $L_6$).

\begin{figure}[h]
\begin{tikzpicture}[x=1cm,y=1cm]
\clip(-11.09,4) rectangle (0,6.5);
\draw [line width=1pt] (-3,6)-- (-3,4);
\draw [fill=white,color=white] (-3,4.58) circle (3.5pt);
\draw [fill=black,color=black] (-2,4.58) circle (3.5pt);
\draw [fill=black,color=black] (-3,5.44) circle (3.5pt);
\draw [line width=1pt] (-2,6)-- (-2,4);
\draw [fill=white,color=white] (-2,5.44) circle (3.5pt);
\draw [line width=1pt] (-7,6)-- (-7,4);
\draw [line width=1pt] (-6,6)-- (-6,4);
\draw [line width=1pt] (-5,6)-- (-5,4);
\draw [line width=1pt] (-4,6)-- (-4,4);
\draw [line width=1pt] (-8,5.44)-- (-1,5.44);
\draw [line width=1pt] (-7.96,4.58)-- (-1,4.58);
\draw[color=black] (-6.97,6.32) node {$L_1$};
\draw[color=black] (-6.03,6.32) node {$L_2$};
\draw[color=black] (-5.01,6.32) node {$L_3$};
\draw[color=black] (-3.95,6.32) node {$L_4$};
\draw[color=black] (-2.97,6.32) node {$L_5$};
\draw[color=black] (-2,6.32) node {$L_6$};
\draw[color=black] (-8.2,5.5) node {$T_1$};
\draw[color=black] (-8.2,4.6) node {$T_2$};
\end{tikzpicture}
\caption{Case D2iii}
\label{FigD2iii}
\end{figure}
\end{proof} 

\begin{remark}\label{InfRem}
Let ${\mathcal L}$ be a finite set of 3 or more skew lines in $\PP^3_\field$
where $\field$ is the algebraic closure of $\mathbb Q$.
In this remark we describe an algorithmic answer for the question of when
$G_{\mathcal L}$ is finite. If the generators for ${\mathcal L}$
given in Proposition \ref{Prop1}(a) commute then $G_{\mathcal L}$ is 
abelian, in which case it is finite if and only if they all have finite order. So suppose
they do not all commute with each other. It is not hard, given a finite set of
2 by 2 invertible matrices to tell if they are contained in a dihedral group
(note that half the elements of a dihedral group have order 2
and the other elements give a cyclic group). So also suppose
that the generators do not generate a dihedral group.
Then by the classification of
the finite subgroups of $SL_2(\mathbb C)$, the group $G_{\mathcal L}$ must have order 
24, 48 or 120 (and the order identifies the group up to isomorphism). So if one iteratively generates the elements of
$G_{\mathcal L}$ by taking all products of up to 2, then 3 then 4, etc.\ generators,
either one gets the whole group after a while and then one knows the order and hence the group,
or one gets more than 120 elements, and hence $G_{\mathcal L}$ is infinite. 
\qed\end{remark}

Note, up to projective equivalence, that there are infinitely many 
ways to choose a set $\mathcal L=\{L_1,L_2,L_3,L_4\}$ of four skew lines in $\PP^3_{\mathbb C}$ 
with $|G_{\mathcal L}|=1$. Just choose any four lines in the same ruling of a smooth quadric.
If $T$ is any transversal (hence $T$ is any line in the other ruling), 
the cross ratio of the intersection of the four lines with $T$ is the same
for all choices of $T$, but fixing $L_1,L_2,L_3$, infinitely
many cross ratios arise depending on the choice of $L_4$.
However, for each choice of $L_4$ there are at most 5 other 
choices of $L_4$ that have projectively equivalent cross ratios
and hence give a projectively equivalent set of four lines.
Things are quite different when the four lines do not all lie on the same quadric.

\begin{corollary}
Let $m$ be a fixed integer $m>2$. Then there are choices of
four skew lines $\mathcal L=\{L_1,L_2,L_3,L_4\}$, $L_i\subset \PP^3_{\mathbb C}$,
such that $|G_{\mathcal L}|=m$, in which case $G_{\mathcal L}$ is a cyclic group.
Moreover, up to projective equivalence, there are
only finitely many such choices. 
\end{corollary}

\begin{proof}
By Remark \ref{AlgorithmRem}, we can find choose four lines 
whose group is cyclic of order $m$, by choosing generators 
$\alpha,\beta$ of a multiplicative cyclic of order $m$.
But there are only finitely many
ways to choose $\alpha,\beta$ from a multiplicative cyclic of order $m$.
Up to projective equivalence, we may assume the lines $L_1,L_2,L_3$ 
are in standard position and the transversals are standard, and then
the choices of $\alpha$ and $\beta$
determine $L_4$, so there are only finitely choices of $L_4$.
\end{proof}

\subsection{Additional combinatorial questions}\label{Sec:CombQuestsl}

The groupoid gives new perspectives on some traditional combinatorial problems.
For example, a traditional area of study is that of spreads in $\PP^3_\field$
when $|\field|<\infty$.

A \emph{spread} is a set $\mathcal L=\{L_1,\ldots,L_s\}$ of skew lines in $\PP^3_\field$. 
It is \emph{maximal} if it is contained in no larger such set.
It is \emph{full} if every point of $\PP^3_\field$ is in one of the lines and \emph{partial} otherwise.

A traditional problem is to understand what values of $s$ arise for maximal spreads.
Since every line has $q+1$ points, where $q=|\field|$ and $|\PP^3_\field|=q^3+q^2+q+1$,
the maximum possible value for $s$ is $(q^3+q^2+q+1)/(q+1)=q^2+1$, which is what we get
for a full spread. 

Full spreads always exist by a construction 
\cite{BB}
which
\cite{kettinger2023}
notes is essentially the Hopf fibration.
Given any degree 2 extension $\field\subset J$ of fields, we have a canonical map
$h\colon \PP^3_\field=(\field^4)^*/\field^*=(J^2)^*/\field^*\to (J^2)^*/J^*=\PP^1_J$, where $^*$ means
the nonzero elements.
(When $\field= \mathbb R$ and $J=\mathbb C$, composing with the antipodal quotient
$S^3\to \PP^3_\field$  gives the original Hopf fibration $S^3\to S^1$.)
The fibers are the quotients $V^*/\field^*$ where $V\subset J^2$ is a $J$-vector subspace
of dimension 1; i.e., the fibers are lines in $\PP^3_\field$. These lines give a full spread
for $\PP^3_\field$. 

For our next result it is helpful to express $h$ explicitly when $|\field|<\infty$.
The case of characteristic 2 requires special handling (which you can see in 
\cite{Ganger}), 
so, to simplify our discussion, we avoid it here. Now, assuming $|\field|<\infty$ and ${\rm char}(\field)>2$, we can pick $\alpha\in J$
such that $J=\field[\alpha]$ and $\alpha^2\in\field$. Then, following 
\cite{Ganger},
the Hopf map $h\colon \PP^3_{\field}\to\PP^1_J$ is $(a:b:c:d)\mapsto (a+\alpha b,c+\alpha d)$.
Note $(a:b:c:d)$ and $(\alpha^2b:a:\alpha^2d:c)$ represent different points
of $\PP^3_{\field}$, but $h((a:b:c:d))=h((\alpha^2b:a:\alpha^2d:c))$, so
the fiber $L_p$ containing $p=(a:b:c:d)$ is the line spanned by the points 
$(a:b:c:d)$ and $(\alpha^2b:a:\alpha^2d:c)$.
Note that 
$$\alpha(a:b:c:d)+(\alpha^2b:a:\alpha^2d:c)=(a+\alpha b)(\alpha:1:0:0)+(c+\alpha d)(0:0:\alpha:1)$$
so $L_p$ meets the line $H_1\colon (x-\alpha y, z-\alpha w)$ spanned by 
$(\alpha:1:0:0)$ and $(0:0:\alpha:1)$. I.e., $H_1$ is a transversal for the fibers of
$h$. The Galois involution given by $\alpha\mapsto-\alpha$ gives the other transversal,
$H_2\colon (x+\alpha y, z+\alpha w)$. The transversals are defined over
$J$ but not over $\field$ and have no points of $\PP^3_\field$. 

\begin{theorem}\label{Thm:HopfFibr}
Let $\mathcal L$ be a full spread for $\PP^3_\field$, assuming 
$|\field|<\infty$ and ${\rm char}(\field)>2$. If $G_{\mathcal L}$ is abelian,
then $\mathcal L$ is projectively equivalent to the spread given by
the fibers of the Hopf map $h$. In particular, $G_{\mathcal L}$ 
is cyclic of order $q+1$. 
\end{theorem}

\begin{proof} Let $q=|\field|$ and let $\mathcal L=\{L_1,\ldots,L_{q^2+1}\}$.
If $G_{\mathcal L}$ is abelian, then there are two transversals (counted with multiplicity) by 
Theorem \ref{TransThm}. If the transversal $T$ has multiplicity 2, then $T$ is defined
over $\field$, hence has $q+1$ points, so any spread for which $T$ is a transversal has
at most $q+1$ lines. But a full spread has $q^2+1>q+1$ lines. 
Thus a full spread with an abelian group has
two distinct transversals, $T_1$ and $T_2$. If any point $T_1\cap L_i$ is defined
over $\field$, all points $T_1\cap L_j$ are, since they form an orbit of the groupoid.
But then $T_1$ is defined over $\field$ which we have just seen is impossible.
Thus no point $T_1\cap L_i$ (and similarly none for $T_2$) is defined over $\field$.
But $T_1$ is one of the two ruling lines on the quadric $Q_{123}$ determined by $L_1,L_2,L_3$
which meet $L_4$. Since $T_1\cap L_4$ is not defined over $\field$,
it is defined over the degree 2 extension $\field\subset J=K[\alpha]\subset \overline{\field}$
for some $\alpha^2\in\field$.
Over $J$, $T_1$ has $q^2+1$ points, which must be the points $T_1\cap L_j$.
Thus no point of $T_1$ over $J$ is defined over $\field$ (and the same for $T_2$).

Note that every point $(a:b:c:d)\in\PP^3_J$ which is not in $\PP^3_\field$ is on 
a unique line defined over $\field$. This is because we can write
$(a:b:c:d)=(u_1+\alpha U_1: u_2+\alpha U_2: u_3+\alpha U_3: u_4+\alpha U_4)$.
Since $(a:b:c:d)$ is not defined over $\field$, we see that 
$(u_1: u_2: u_3: u_4)$ and $(U_1: U_2: U_3: U_4)$ represent different points
of $\PP^3_\field$ and so define a line in $\PP^3_\field$ which over $J$ contains
$(a:b:c:d)$. But if $(a:b:c:d)$ were on two
lines defined over $\field$, it would itself be defined over $\field$.

So now we see if $L\subset\PP^3_J$ is a line no point of which is defined over $\field$,
then each point $p\in L$ of the $q^2+1$ points of $L$ defines a line $L_p$ defined over $\field$.
If $p_1,p_2$ are two different points of $L$ such that $L_{p_1}$ and $L_{p_2}$ meet,
the lines $L_{p_1}$ and $L_{p_2}$ determine a plane defined over $\field$ which must 
contain the Galois dual $L^*$ of $L$, hence $L$ and $L^*$ meet and the point where they meet
is dual to itself and so must be defined over $\field$, contrary to our assumption that
no point of $L$ is defined over $\field$. Thus a line $L$ with no points
defined over $\field$ determines a canonical full spread for $\PP^3_\field$.
We now show these spreads are projectively equivalent.

Pick two distinct points of the form $T_1\cap L_j$, say 
$u=(u_1+\alpha U_1: u_2+\alpha U_2: u_3+\alpha U_3: u_4+\alpha U_4)$ is $T_1\cap L_{j_1}$ and
$v=(v_1+\alpha V_1: v_2+\alpha V_2: v_3+\alpha V_3: v_4+\alpha V_4)$  is $T_1\cap L_{j_2}$.
The Galois automorphism induced by $\alpha\mapsto -\alpha$ fixes
the points of $\PP^3_\field$ and hence the lines of the spread,
so $u'=(u_1-\alpha U_1: u_2-\alpha U_2: u_3-\alpha U_3: u_4-\alpha U_4)$ is $T_2\cap L_{j_1}$ and
$v'=(v_1-\alpha V_1: v_2-\alpha V_2: v_3-\alpha V_3: v_4-\alpha V_4)$ is $T_2\cap L_{j_2}$. Thus 
$u_+=(u+u')/2,u_-=(u-u')/(2\alpha) \in L_{j_1}$ and 
$v_+=(v+v')/2,v_-=(v-v')/(2\alpha)\in L_{j_2}$ are defined over $\field$.
Because $u$ is not defined over $\field$, $u_+$ and $u_-$ are distinct points of $L_{j_1}$
and likewise $v_+$ and $v_-$ are distinct points of $L_{j_2}$.
Since $L_{j_1}$ and $L_{j_2}$ are distinct lines, the points $u_+,u_-,v_+$ and $v_-$
span $\PP^3_\field$. Thus the matrix 
$M=\begin{pmatrix}
u_1/\alpha^2 & U_1 & v_1/\alpha^2 & V_1 \\
u_2/\alpha^2 & U_2 & v_2/\alpha^2 & V_2 \\
u_3/\alpha^2 & U_3 & v_3/\alpha^2 & V_3 \\
u_4/\alpha^2 & U_4 & v_4/\alpha^2 & V_4 \\
\end{pmatrix}$
is nonsingular and, as an element of $PGL_4(\field)$, takes
$(\alpha:1:0:0)$ to $u$ and $(0:0:\alpha:1)$ to $v$;
i.e., it takes $T_1$ to $H_1$ and $T_2$ to $H_2$, and thus
the lines in $\mathcal L$ (each of which is defined by a point of 
$T_1$ and its Galois dual in $T_2$) to the fibers of $h$.
Thus $\mathcal L$ is projectively equivalent to
the spread coming from $h$. By Proposition \ref{twoTrans}(c),
the group is finite and cyclic; by 
\cite{Ganger}
it has order $q+1$
(Ganger shows the group for the Hopf spread is $J^*/\field^*$,
as a consequence of showing for any full spread that 
$G_1$ acts transitively on the points of $L_1$).
We can use
Theorem \ref{CrossRatioRatios} to recover this result.
A fiber of $h$ meets $H_1$ at a point $(\alpha a:a:\alpha b:b)$
where $(a:b)$ is a point of $\PP^1_J$, and that fiber meets
$H_2$ at $(-\alpha a^*:a^*:-\alpha b^*:b^*)$, where
the asterisk denotes the Galois conjugate (i.e.,
if $e,f\in\field$, then $(e+\alpha f^*=e-\alpha f$).
Thus four points on $H_1$ can be parameterized by
$(a_i:b_i)$, $1\leq i\leq 4$, and the corresponding points on $H_2$ are
$(a^*_i:b^*_i)$, $1\leq i\leq 4$, where we assume $(a_1:b_1)=(1:0)$
in order to compute generators for $G_1$. The cross ratio ratio is thus
$$\frac{\chi_1}{\chi_2}=\frac{\frac{(a_4b_2-a_2b_4)(a_3b_1-a_1b_3)}{(a_4b_1-a_1b_4)(a_3b_2-a_2b_3)}}
{\frac{(a^*_4b^*_2-a^*_2b^*_4)(a^*_3b^*_1-a^*_1b^*_3)}{(a^*_4b^*_1-a^*_1b^*_4)(a^*_3b^*_2-a^*_2b^*_3)}}=
\frac{\frac{(a_4b_2-a_2b_4)(0-b_3)}{(0-b_4)(a_3b_2-a_2b_3)}}
{\frac{(a^*_4b^*_2-a^*_2b^*_4)(0-b^*_3)}{(0-b^*_4)(a^*_3b^*_2-a^*_2b^*_3)}}=
\frac{\chi_1}{\chi_1^*}$$
where $\chi_1$ is the cross ratio of the four given points on $H_1$ and 
$\chi_2$ is the cross ratio of the four corresponding points on $H_2$.
Every nonzero element of $J$ occurs as the cross ration $\chi_1$ for some set of four points on $H_1$,
hence $G_1$ consists of the subgroup of $J^*$ of all elements of the form $c/c^*$, $c\in J^*$.
This is a finite group, hence cyclic. We have a homomorphism $J^*\to J^*$ given by $c\mapsto c/c^*$.
The kernel consists of the elements of $\field^*\subset J^*$. Thus $|G_1|=\frac{|J^*|}{|K*|}=\frac{p^2-1}{p-1}=p+1$.
\end{proof}

Typically there are spreads not projectively equivalent to the one given by the
Hopf fibration, so they must have nonabelian groups (see 
\cite{Ganger}
for an explicit example).
But all full spreads $\mathcal L$ for $\PP^3_\field$ whose group is abelian
are projectively equivalent and hence $\PP^3_\field$ is a single orbit of $C_\mathcal L$
(because that is true for the spread given by the fibers of the Hopf map $h$).
This prompts the following questions.

\begin{question}\label{quest:Hopfspread}
Is there a nonabelian group $G$ arising as the group for two or more full spreads
which are not projectively equivalent? Is the number of such spreads
related to some group theoretical property of $G$?
\end{question}

\begin{question}\label{quest:maxspreads}
What groups $G_\mathcal L$ arise for maximal partial spreads $\mathcal L$?
Is there always only a single orbit? If not, what do their orbits look like?
\end{question}

We can also ask about minimal spreads.

\begin{question}\label{quest:minspreads}
What minimal spreads $\mathcal L'$ are there in a maximal spread $\mathcal L=\{L_1\ldots, L_s\}$
such that $G_{\mathcal L'}=G_{\mathcal L}$? Note that $C_{\mathcal L'}$ and $C_{\mathcal L}$
have the same groups $G_{\mathcal L'}=G_{\mathcal L}$ and thus for any line $L_i\in \mathcal L'$,
the orbits on $L_i$ of $C_{\mathcal L'}$ and $C_{\mathcal L}$ are the same.
Suppose we obtain ${\mathcal L''}$ by removing one line
from ${\mathcal L'}$. How do the orbits change? Is it possible for $\mathcal L_1$
to be contained in two different maximal spreads $\mathcal L_2$ and $\mathcal L_3$
if $G_{\mathcal L_1}=G_{\mathcal L_2}$?
\end{question}

\section{An application to algebraic geometry}\label{sec: application to AG}

In this part of the paper we work over an algebraically closed field $\mathbb F$,
so when we write $\PP^3$ without qualification we mean $\PP^3_{\mathbb F}$.
We will keep $\field$ as denoting an arbitrary field  but with $\field\subseteq \mathbb F$.
The initial work on geproci sets, namely 
\cite{CM}
and 
\cite{POLITUS},
was over the complex numbers; as a result, much of the subsequent work, 
including 
\cite{POLITUS2, PSS, WZ}, also focused on the complex setting.
Now that 
\cite{kettinger2023, Ganger} 
have shown how natural and interesting
the geproci concept is over finite fields there is more incentive to
work in a characteristic free way.
Nonetheless, some results of 
\cite{CM,POLITUS}
which we suspect are true in general 
do not yet have characteristic free proofs, so some of our work here which
applies results from 
\cite{CM,POLITUS} 
have a characteristic 0 assumption (typically by assuming $\mathbb F=\mathbb C$), even though they might be true
in general.

So let $Z\subset \PP^3$ be a finite set of points and let $P\in\PP^3$ be a general point.
We denote by an overbar, $\overline{\hbox to.1in{\vbox to .105in{\vfil}\hfil}}: \PP^3\DashedArrow \PP^3$, the rational map
(defined away from $P$) given by projection from $P$ to a given plane $H\cong\PP^2$.
Thus $\overline{Z}_{P,H}$ (or for simplicity just $\overline{Z}$) is
the image of $Z$ under the projection. Following 
\cite{POLITUS},
we say that a set
$Z$ of $ab$ points is $(a,b)$-geproci if $\overline{Z}$
is the intersection of plane algebraic curves in $H$ of degrees $a$ and $b$ with $a \leq b$;
i.e., $\overline{Z}$ is a transverse complete intersection of type $(a,b)$.
We say that $Z$ is $\{a,b\}$-geproci if we drop the condition $a \leq b$.

Given any plane $H'$ and finite set $Z\subset H'$, the projection from a general point $P$ to $H$
restricts to an isomorphism $H'\to H$. Thus if $Z$ is 
a transverse complete intersection of type $(a,b)$ in $H'$, then $\overline Z$
is a transverse complete intersection of type $(a,b)$ in $H$.
Thus a degenerate set $Z$ is $(a,b)$-geproci if and only if it is a
transverse complete intersection of type $(a,b)$ in a plane containing it. 
When we say a geproci set is trivial, we just mean it is
degenerate (and hence a transverse complete intersection of 
type $(a,b)$ in some or equivalently every plane containing it). 

The real interest is in understanding nondegenerate geproci sets, but
the trivial case suggests an approach to doing so. Suppose $Z=C_1\cap C_2$ is a set of $ab$ points 
where $C_1$ and $C_2$ are space curves of degrees $a$ and $b$, respectively.
Because $P$ is general, we have $\overline Z=\overline{C_1\cap C_2}=\overline{C_1}\cap\overline{C_2}$,
so $Z$ is $\{a,b\}$-geproci. In the trivial case, we have $C_1,C_2\subset H'$ for some plane $H'$,
but when $Z$ is nondegenerate we can imagine that curves $C_1$ and $C_2$ 
might sometimes still exist with $Z=C_1\cap C_2$, in which case
we could get nontrivial geproci sets $Z$ this way.
More generally, we can look for $\{a,b\}$-geproci sets $Z$ occurring as 
subsets of space curves $C_1$ of degree $a$, even if there is
no space curve $C_2$ of degree $b$ with $Z=C_1\cap C_2$. 

This suggests
a rough classification of $(a,b)$-geproci sets $Z$ (so $a\leq b$): 

(1) $Z=C_1\cap C_2$ where 
$C_1,C_2$ are space curves with $\deg(C_1)=a$ and $\deg(C_2)=b$
(this includes the trivial case); 

(2) $Z$ is not of type (1) but there is a space curve $C_1$ with $Z\subset C_1$ and $\deg(C_1)\leq\max(a,b)$; and

(3) $Z$ is not of types (1) or (2).

All cases of types (2) and (3) are necessarily nontrivial.

For all known examples of nondegenerate $(a,b)$-geproci sets $Z$ of types (1) or (2),
each curve $C_1, C_2$ for type (1) or the curve $C_1$ for type (2),
is a union of skew lines (and necessarily the number of points of $Z$ on 
each line is the same for all of the lines since $\overline{Z}$ is an $(a,b)$ complete intersection).
This prompted the following terminology which we now recall.

\begin{definition}
Let $Z$ be $\{a,b\}$-geproci.
\begin{enumerate}
\item[(a)] We say $Z$ is $[a,b]$-geproci if, for some set $\mathcal L$ of $b$ skew lines,
$Z$ consists of $a$ points on each of the $b$ lines
and $\overline{Z}$ is the transverse intersection of a plane curve $C$ of degree $a$
with the image $\cup_{L\in \mathcal L}\overline{L}$, in which case we say
$Z$ is $[a,b]$-geproci with respect to $\mathcal L$.

\item[(b)] We say $Z$ is an $\{a,b\}$-\emph{grid} (which we write as $(a,b)$-grid to signify $a\leq b$) if $Z=C_1\cap C_2$,
where $C_1$ is a curve consisting of $a$ skew lines and $C_2$ is a curve  consisting of $b$ skew 
lines, such that each component of $C_1$ meets each component of $C_2$ in exactly one point.
\item[(c)] We say $Z$ is an $\{a,b\}$-\emph{half grid} (or $(a,b)$-half grid to indicate $a\leq b$) if it is not
a grid but is either $[a,b]$-geproci or $[b,a]$-geproci, in which case
we also say $Z$ is an $[a,b]$-\emph{half grid} or $[b,a]$-\emph{half grid} resp. 
(In particular, an $[a,b]$-geproci set $Z$ is either an $\{a,b\}$-grid or an $[a,b]$-half grid.)
\end{enumerate}
\end{definition}

Our first main algebraic geometric result (currently only over the complex numbers), 
shows that the curves $C_i$ in the rough classification into types (1), (2) or (3)
above, in the case of nondegenerate geproci sets $Z$, are unions of skew lines.
Thus a nondegenerate $\{a,b\}$-geproci set $Z$ is either an $\{a,b\}$-grid,
an $[a,b]$-half grid, a $[b,a]$-half grid or neither a grid nor a half grid. 

\begin{theorem}\label{curve theorem}
Assume $\mathbb F=\mathbb C$ and let $Z$ be a nondegenerate $(a,b)$-geproci set, hence $b\geq a$.
\begin{enumerate}
\item[(a)] If $Z$ has type (1), then each curve $C_1$ and $C_2$
is a union of skew lines and $Z$ is an $(a,b)$-grid.
Moreover, $C_1$ fails to be unique if and only if: $2=a=b$; or $2<a=b$ (in which case the pair $C_1,C_2$ is uniquely determined).
And $C_2$ fails to be unique if and only if $2=a\leq b$ or $2<a=b$ (but in the latter case, the pair $C_1,C_2$ is uniquely determined).
\item[(b)] If $Z$ has type (2), let $C$ be a curve of least degree $c$ which contains $Z$.
Then $4\leq c\in\{a,b\}$ and $C$ is a union of $c$ skew lines.
In this case, $Z$ is a $[d,c]$-half grid where $\{c,d\}=\{a,b\}$.
Moreover, the curve $C$ is unique except in at least some cases when $a<b=c$.
\end{enumerate}
\end{theorem}

We give the proof directly after Proposition \ref{curve prop}.

Regarding nonuniqueness, we note that the $[3,4]$-half grid 
$Z_{D_4}$, which is the projectivization of the roots of the $D_4$
root system (see Example \ref{Example:D4}), falls under Theorem \ref{curve theorem}(b).
In this case $C_1$ is a union of 4 skew lines
where each line contains three of the points 
\cite{POLITUS},
 but there are multiple choices of the 4 skew lines.
Also, by 
\cite[Theorem 4.6]{POLITUS},
there exists an 
$[a,b]$-half grid for every $a$ and $b$ with $4\leq b\leq a+1$.
Thus we can find examples of part (b) of Theorem \ref{curve theorem} 
for each of the cases 
$a<\deg(C_1)=b=a+1$,
$a=\deg(C_1)=b$ and
$\deg(C_1)=b<a$, and for the latter two the curve $C_1$ is unique.
(The idea, when $a\geq3$, is to construct a nongrid $[a,a+1]$-geproci set $Z$
on specially chosen skew lines $\mathcal{L}=\{L_1,\ldots,L_{a+1}\}$. 
For each subset $S\subseteq\mathcal{L}$, let $Z_S=Z\cap(\cup_{L\in S}L)$.
By 
\cite[Lemma 4.5]{POLITUS},
the proof of which is characteristic free,
$Z_S$ is $[a,|S|]$-geproci. By an appropriate choice of $S$ 
one can be sure that $Z_S$ is not a grid as long as $|S|>3$.)

There are thus only three types of nondegenerate geproci sets:
grids, half grids and nongrid non-half grids, and (by Theorem \ref{curve theorem}
and Lemma \ref{lem: c is a or b}) over the complex numbers 
a nongrid non-half grid $(a,b)$-geproci set is not contained in any curve
of degree $b$ or less (we have no reason to think this is 
different in positive characteristics).
Grids are easy to understand since they consist of either sets of $a$ collinear points, sets of
$2a$ points with $a$ points on each of two skew lines, or sets of $ab$ points on a smooth quadric
obtained as the points of intersection 
of $a$ lines in one ruling intersecting $b$ lines in the other ruling.
We know only three examples of nongrid non-half grids
in characteristic 0 (see 
\cite{POLITUS}
and Example \ref{Penrose et alia}).
Although there are many examples in positive characteristics
(based on the theory of spreads; see 
\cite{kettinger2023}
and \S\ref{Sec:CombQuestsl}), they are still rather mysterious.
Thus our main focus here is to understand half grids. 
By Proposition \ref{Prop:[a,b]geprociImpliesCollCompl}(a), a half grid on a set $\mathcal L$ of 3 or more lines
is a union of grids on subsets of 3 of the lines, and this is true in all characteristics. Even more,
by Proposition \ref{Prop:[a,b]geprociImpliesCollCompl} and Proposition \ref{p. chg and finite union orbits}, 
a half grid is a union of $C_{\mathcal L}$ orbits, so the combinatorics of skew lines
play an important role. 

The rest of the paper is motivated by the problem of understanding half grids.
There is not yet a complete classification of half grids,
but our results provide a method for constructing
all possible half grids whose group is abelian,
which includes the case of all half grids on 4 skew lines.
See for example Remark \ref{AlgorithmRem} for the case of half grids coming from
orbits on skew lines with two transversals. 

Proposition \ref{Prop:[a,b]geprociImpliesCollCompl} tells us that collinearly complete sets on 3 
skew lines are grids and thus geproci, and that $[a,b]$-geproci sets on any number $b$ of skew lines
are collinearly complete sets of those lines. 
Our next result gives a converse.
Since collinearly complete sets are unions of orbits by Proposition \ref{p. chg and finite union orbits},
this shows that being a geproci grid or half grid and being a finite union of orbits are two sides of the same coin. 

\begin{theorem}\label{CHGthm}
Assume $\mathbb F$ is algebraically closed of any characteristic, let $\mathcal L$ be a set of 
$b\geq3$ 
skew lines
in $\PP^3_{\mathbb F}$ and let $Z$ be a set of points with exactly $a$ points on each of the $b$ lines.
Then $Z$ is collinearly complete for $\mathcal L$ (or equivalently, $Z$ is a union of $C_{\mathcal L}$ orbits) 
if and only if either $Z$ is an $\{a,b\}$-grid or $Z$ is an $[a,b]$-half grid
with respect to $\mathcal L$.
\end{theorem}

We give the proof at the end of \S\ref{sec:ProofOfThmB} on page \pageref{proof_CHGthm}.

\subsection{A combinatorial property of $[a,b]$-geproci sets}\label{HalfGridsAreCollCompl}

The following result adapts the proof of 
\cite[Proposition 4.14]{POLITUS}
but does not make any assumptions about the characteristic.

\begin{proposition}\label{Prop:[a,b]geprociImpliesCollCompl}
Let $\mathcal L=\{L_1,\ldots,L_b\}$ be skew lines $L_i\subset\PP^3$ and
let $Z$ be $[a,b]$-geproci such that each line $L_i$ contains exactly $a$ points of $Z$.
Then $Z$ is $[a,b]$-geproci with respect to $\mathcal L$. Moreover:

\begin{itemize}
\item[(a)] If $b\geq 3$, then $Z$ is collinearly complete for $\mathcal L=\{L_1,\ldots,L_b\}$, and for
every subset $\{L_{i_1},L_{i_2},L_{i_3}\}$ of three of the lines,
$Z\cap(L_{i_1}\cup L_{i_2}\cup L_{i_3})$ is a $\{3,a\}$-grid (where $Q$
is the unique quadric containing the lines $L_{i_1}, L_{i_2}, L_{i_3}$ and
the grid lines are the ruling lines of $Q$ which meet $Z$).
\item[(b)] If $a\leq 2$, then $Z$ is a grid.
\item[(c)] If $b\leq 3$, then $Z$ is a grid.
\item[(d)] In particular, if $Z$ is a half grid (i.e., not a grid), then $a\geq 3$, $b\geq4$
and $Z$ is collinearly complete for $\mathcal L$.
\end{itemize}
\end{proposition}

\begin{proof} 
Say $Z$ is $[a,b]$-geproci with respect to $\mathcal L'=\{L_1',\ldots,L_b'\}$.
If $a>b$, then since $L_i$ contains $a$ points of $Z$, it contains more than one point
of some $L_j'$, hence $L_i=L_j'$. Thus $\mathcal L\subset \mathcal L'$, hence
$\mathcal L = \mathcal L'$.
If $a=b$, then there is a linear pencil $\Lambda$ of curves $C\subset H$ of degree $a$ containing $\overline Z$
and $\overline Z=C_1\cap C_2$ for any two $C_1,C_2\in\Lambda$ distinct curves. 
Since $\overline Z\subset \cup_i\overline{L_i}$, we see $\cup_i\overline{L_i}\in\Lambda$. Hence
either $\cup_i\overline{L_i}=\cup_i\overline{L_i'}$ (in which case 
$Z$ is $[a,b]$-geproci with respect to $\mathcal L$) or
$\overline Z=(\cup_i\overline{L_i})\cap(\cup_i\overline{L_i'})$ (in which case 
$Z$ is again $[a,b]$-geproci with respect to $\mathcal L$).
If $a<b$, then since $\overline Z\subset \cup_i\overline{L_i}$
and $\overline Z=C\cap(\cup_i\overline{L_i'})$ for some curve $C\subset H$
of degree $a$, either $C\subset \cup_i\overline{L_i}$ or
$\overline Z=C\cap(\cup_i\overline{L_i})$. In the latter case
$Z$ is $[a,b]$-geproci with respect to $\mathcal L$. In the former case,
$C$ is a union of lines, say $C=\cup_{i=1}^a \overline{L_i}$.
Thus $C$ is the image of $\cup_{i=1}^a L_i$. But $a<b$, so 
$Z\not\subset\cup_{i=1}^a L_i$. And points not in
$\cup_{i=1}^a L_i$ do not map to $\cup_{i=1}^a \overline{L_i}$
under a general projection; this contradicts
$\overline Z\subset C=\cup_{i=1}^a \overline{L_i}$.

(a) Note that $Z$ is collinearly complete for $\mathcal L=\{L_1,\ldots,L_b\}$
if it is so for each subset 
$$\{L_{i_1},L_{i_2},L_{i_3}\}\subset \mathcal L$$
of 3 of the lines. And for every nonempty subset $S\subseteq \mathcal L$,
the subset $Z'=\cup_{L\in S}Z\cap L$ of $Z$ is $[a,|S|]$-geproci 
with $a$ points on each line $L\in S$ by 
\cite[Lemma 4.5]{POLITUS}.
Thus it is enough to consider the case that $b=3$, in which case
there is a unique smooth quadric $Q$ containing $L_1,L_2,L_3$.
Every transversal $T$ for $\mathcal L=\{L_1,L_2,L_3\}$ also lies on $Q$.
Being collinearly complete means $Z$ is the $[a,3]$-grid 
determined by the curves $C_1$ and $C_2$,
where $C_1$ consists of the transversals through points of $Z$ and 
$C_2=L_1\cup L_2\cup L_3$.

If $a=1$, then the points of the image $\overline{Z}$ of $Z$ under projection from a general point
to a plane are collinear, hence the points of $Z$ are collinear, so they lie on
a transversal $T\subset Q$ for $\mathcal L$, so 
$Z$ is collinearly complete for $\mathcal L$.
Thus now we can assume $a\geq2$.

Assume $Z$ is not a grid. Thus some transversal $T$ for $\mathcal L$
contains a point $P\in Z$ but does not contain 3 points of $Z$. 
For specificity we may assume $P\in L_1$ but $T\cap L_2$ is empty.
Thus the plane $H'$ spanned by $L_3$ and $P$ contains 
$a+1$ points of $Z$. Let $Z'$ be the remaining $2a-1$ 
points of $Z$, so $Z'$ consists of the $a$ points on 
$L_2$ and the $a-1$ points on $L_1$ other than $P$. 
Let $\pi$ be the projection from a general point of $H'$ to a general plane $\Pi$. 
Then $\pi(Z)$ consists of $a+1$ points on $\pi(L_3)$ and
$a$ each on $\pi(L_1)$ and $\pi(L_2)$
(with $\pi(P)$ on both $\pi(L_1)$ and $\pi(L_3)$). Using B\'ezout's Theorem 
successively on these sets of collinear points, we see that any curve of degree $a$ 
containing $\pi(Z)$ must contain $\pi(L_1)\cup\pi(L_2)\cup\pi(L_3)$. 
But this is impossible if $a=2$, so say $a\geq3$. We have that the 
the dimension of the space of forms of degree $a$ vanishing on
$\pi(Z)$ is the dimension $D=\binom{a-3+2}{2}$ of the space of all forms of degree $a-3$.
However the dimension $d$ of the space of forms of degree $a$ vanishing on $\overline{Z}$
(the image of $Z$ under a general projection)
must, by the geproci property, satisfy $d>\binom{a-3+2}{2}=D$, 
but by semicontinuity we must have $D\geq d$.
This contradiction implies $Z$ is a grid.

(b) If $a=1$ then $Z$ is a set of $b$ collinear points, hence $Z$ is a $(1,b)$-grid.
So say $a=2$. If $b=1$ or $b=2$, $Z$ is clearly a grid, so assume $b\geq 3$.
By (a) and Proposition \ref{p. chg and finite union orbits}, $Z$ is a union of orbits,
and by Proposition \ref{Prop:not 2} every orbit not in a transversal for the $b$ lines has
more than $|Z|=2b$ points. Thus $Z$ must consist of $a$ points on each of two transversals,
hence $Z$ is a grid. 

(c) If $b=1$, then the points of $Z$ all lie on $L_1$ hence are collinear and $Z$ is a grid.
If $b=2$, then $Z$ is a grid (just connect up each point on $L_1$ with a different
point on $L_2$ to get a set of grid lines transverse to $L_1,L_2$.)
And the case $b=3$ is covered by (a).

(d) This is now immediate. 
\end{proof}

\begin{corollary}\label{grids and groups Cor}
Let $a,b\geq3$. For any field $\field$, let $Z\subset\PP^3_\field\subset\PP^3_{\overline{\field}}$ be 
$[a,b]$-geproci and let $\mathcal L=\{L_1,\ldots,L_b\}$ be a set of $b$ 
skew lines with $a$ points of $Z$ on each line.
Then $Z$ is an $\{a,b\}$-grid if and only if $|G_{\mathcal L}|=1$.
\end{corollary}

\begin{proof}
Assume $Z$ is an $\{a,b\}$-grid. Then there are curves $C_a$ and $C_b$ where $C_a\cap C_b=Z$,
and where $C_a$ consists of $a$ skew lines and $C_b$ consists of $b$ skew lines.
Since $a\geq 3$, we can pick three of the lines of $C_a$; these are contained in a 
unique smooth quadric $Q$. The lines of $C_b$ are transversal for these three lines,
so $C_b\subset Q$, and the lines of $C_a$ are transversal for the lines of $C_b$ so also
$C_a\subset Q$. Thus every point $z\in Z$ is on two lines in $Q$, but these are the only
lines in $Q$ through $Z$. Hence the lines in $C_a\cup C_b$ are exactly
the ruling lines which meet $Z$.

Any line containing three or more points of $Z$ must likewise be
contained in $Q$ (so must be a component of $C_a$ or $C_b$).
In particular, we have $L_i\subset Q$ for all $i$, hence $|G_{\mathcal L}|=1$
by Corollary \ref{transCor}.

Now assume $|G_{\mathcal L}|=1$, hence, by Corollary \ref{transCor}, there is a quadric $Q$
with $L_i\subset Q$ for all $i$. By Proposition \ref{Prop:[a,b]geprociImpliesCollCompl}(a), 
$Z$ is $[a,b]$-geproci with respect to $\mathcal L$.
Take a point $z\in Z$ and three lines $L_1,L_2,L_3\in \mathcal L$ with $z\in L_1$.
By Proposition \ref{Prop:[a,b]geprociImpliesCollCompl}(a), there is a transversal $T$
for $L_1,L_2,L_3$ through $z$ and $T\cap L_i$ is in $Z$ also for $i=2,3$, hence $T\subset Q$.
Likewise, for every other line $L$ in $\mathcal L$ the set $Z\cap(L_1\cup L_2\cup L)$ is a $\{3,a\}$-grid,
which thus has a grid line $T'$ other than $L_1$ through $z$, hence $T'\subset Q$.
But there is only one line on $Q$ through $z$ other than $L_1$, so $T'=T$.
I.e., $T\cap L$ is also in $Z$. Hence $Z$ is the $\{a,b\}$-grid
whose grid lines are the ruling lines on $Q$ through points of $Z$.
\end{proof}

\begin{example}\label{Example:D4}
Here, following 
\cite[Example 4]{kettinger2023}
but with details, 
we give a characteristic free example of a $[3,4]$-half grid to show that 
Proposition \ref{Prop:[a,b]geprociImpliesCollCompl}(d) 
is sharp in every characteristic. We denote the set of 12 points by $Z_{D_4}$
since in characteristic 0 it is projectively equivalent to 
the projectivization of the 12 roots of the $D_4$ root system. But the 12 actual roots only give 6 points mod 2,
so to get one representation that works in all characteristics we use a representation
that in characteristic 0 is only projectively equivalent to the actual $Z_{D_4}$. 
The representation we use is shown in Figure \ref{3ptPerspective},
taken from 
\cite{POLITUS}.
We adapt an argument from 
\cite{POLITUS}
to show
it is a $[3,4]$-half grid in every characteristic.

Let $\field$ be $\mathbb Q$ or a field of prime order and let $\mathbb F=\overline{\field}$.
The three dashed lines shown are skew, as are the three dotted-dashed lines,
and each dashed line meets each dotted-dashed line in a single point. Those 9 points give a $(3,3)$-grid.
let $P=(a:b:c:d)\in\PP^3_{\mathbb F}$ be a general point.
The three dashed lines are: $x,y-w$; $z,x-w$; and $y,z-w$.
The plane spanned by each line and $P$ is: $(b-d)x-a(y-w)$; $(a-d)z-c(x-w)$; and $(c-d)y-b(z-w)$.
The three dotted-dashed lines are: $x,z-w$; $y,x-w$; and $z,y-w$.
The plane spanned by each line and $P$ is: $(c-d)x-a(z-w)$; $(a-d)y-b(x-w)$; and $(b-d)z-c(y-w)$.
The main diagonal through the points $(0:0:0:1), (1:1:1:1), (1:1:1:2)$ is $x-y,y-z$ and the plane spanned
by it and $P$ is $(b-c)(x-y)-(a-b)(y-z)$. 

The union of the three dashed lines with the diagonal (or alternatively the three dotted-dashed lines
and the diagonal) is a quartic curve that contains all 12 points.
The cone over either curve (with vertex $P$) is a quartic surface which contains all 12 points.
Thus $S_1=((b-d)x-a(y-w))((a-d)z-c(x-w))((c-d)y-b(z-w))((b-c)(x-y)-(a-b)(y-z))$
and $S_2=((c-d)x-a(z-w))((a-d)y-b(x-w))((b-d)z-c(y-w))((b-c)(x-y)-(a-b)(y-z))$
are the equations of these quartic cones.

The cubic cones defined by $G_1=((b-d)x-a(y-w))((a-d)z-c(x-w))((c-d)y-b(z-w))$
and $G_2=((c-d)x-a(z-w))((a-d)y-b(x-w))((b-d)z-c(y-w))$ (i.e., the cone over the curve consisting
of the three dashed lines and the cone over the curve consisting of the three dotted-dashed lines)
each contain the 9 points not on the diagonal. The cubic form $C=G_1-G_2$
vanishes at all 12 points but is not identically 0. Thus $S_1$ (or $S_2$ if you like)
and $C$ both vanish on the 12 lines where each line is spanned by $P$ and one of the 12
points of $Z_{D_4}$. Intersecting this with a plane not containing $P$
shows that the projection of the 12 points to the plane from $P$ is a transverse intersection of 
a cubic curve and a quartic curve (since by direct computation 
none of the planes in the quartic cone is a component of the cubic cone
and both cones contain 12 different lines through $P$)
and hence $Z_{D_4}$ as given is $[3,4]$-geproci in
every characteristic. We will see in the next paragraph that it is not a grid, so
it must be a $[3,4]$-half grid.

We now show that $Z_{D_4}$ is a $[3,4]$-half grid on 4 lines and we check when the fourth line
is tangent to the quadric containing the other three lines: these 4 lines have
two transversals except when the characteristic $p$
is 3, in which case there is only one transversal (but with multiplicity 2).
Consider the lines represented by the dashed lines in 
Figure \ref{3ptPerspective}, namely 
$L_1\colon x, y-w; L_2\colon z, x-w; L_3\colon y, z-w$. The line given by the main diagonal
is $L_4\colon x-y, y-z$. Lines $L_1,L_2,L_3$ determine the quadric
$Q\colon  xy+xz+yz-xw-yw-zw+w^2$. Of the 12 points of $Z_{D_4}$, 
this quadric contains exactly the 9 points not on the diagonal
of the cube; thus $Z_{D_4}$ cannot be a grid.
Indeed, line $L_4$ meets $Q$
in the points $(a:a:a:d)$ where $3a^2-3ad+d^2=0$.
There is no such point when $a=0$ so we can assume $a=1$.
Then we get $d^2-3d+3=0$; this has a unique root if and only if $p=3$.
In that case $L_4$ meets $Q$ at $(1:1:1:0)$ and the unique transversal
is $T\colon x-y+w,x+y+z$. In fact, when $p=3$,
given any 4 skew lines with each line containing 3 of the points of $Z_{D_4}$,
one can check by direct computation that the quadric defined by 3 of the lines is tangent to the fourth line.
Thus, when $p=3$, for any set of 4 skew lines with respect to which $Z$ is a $[3,4]$-half grid,
there is a unique transversal.
(When $p$ is 0 or a prime other than 3, there are two transversals, 
but sometimes they are defined over the prime field and sometimes not.
They are not defined over $\mathbb Q$, nor over $\mathbb Z/p\mathbb Z$ when $p=2$.
When $p>3$, we can rewrite $d^2-3d+3=0$ as $(2d-3)^2=-3$, hence for $\alpha=(-1+(2d-3))/2$ and 
$\beta=(-1-(2d-3))/2$ we have $x^3-1=(x-1)(x-\alpha)(x-\beta)=0$.
The roots (and hence also $d$) are integers mod $p$ exactly when $|\field^*|=p-1$ has order divisible by 3 
\cite{L}.
So the transversals are defined over $\mathbb Z/p\mathbb Z$ exactly when $p\equiv 1\mod 3$.)
\qed\end{example}

\begin{figure}[ht!]
\begin{tikzpicture}[scale=1.1,x=1.0cm,y=1.0cm]
\clip(-4.496680397937174,0.75) rectangle (5.928312362917013,6.5);
\fill[line width=2.pt,color=black,fill=white] (-1.176350188308613,1.7059124529228467) -- (-0.8953798813051128,2.7315518487365438) -- (0.,3.230935550935551) -- (0.,2.) -- cycle;
\fill[line width=2.pt,color=black,fill=lightgray,fill opacity=0.30000000149011612] (0.,3.230935550935551) -- (0.6656059123409749,2.7359609200174266) -- (0.8751724137931034,1.7082758620689655) -- (0.,2.) -- cycle;
\fill[line width=2.pt,color=black,fill=gray,fill opacity=0.60000000149011612] (-1.176350188308613,1.7059124529228467) -- (-0.23559785204630612,1.546900694008833) -- (0.8751724137931034,1.7082758620689655) -- (0.,2.) -- cycle;
\draw [line width=0.4pt] (-4.,1.)-- (0.,2.);
\draw [line width=0.4pt] (0.,2.)-- (0.,6.);
\draw [line width=0.4pt] (0.,2.)-- (3.,1.);
\draw [line width=0.4pt,dash pattern=on 2pt off 2pt] (-4.,1.)-- (0.8751724137931034,1.7082758620689655);
\draw [line width=0.4pt,dash pattern=on 1pt off 2pt on 4pt off 4pt] (3.,1.)-- (-1.176350188308613,1.7059124529228467);
\draw [line width=0.4pt,dash pattern=on 1pt off 2pt on 4pt off 4pt] (-4.,1.)-- (0.,3.230935550935551);
\draw [line width=0.4pt,dash pattern=on 2pt off 2pt] (0.,3.230935550935551)-- (3.,1.);
\draw [line width=0.4pt,dash pattern=on 2pt off 2pt] (-1.176350188308613,1.7059124529228467)-- (0.,6.);
\draw [line width=0.4pt,dash pattern=on 1pt off 2pt on 4pt off 4pt] (0.,6.)-- (0.8751724137931034,1.7082758620689655);
\draw [line width=0.4pt,color=black] (-1.176350188308613,1.7059124529228467)-- (-0.8953798813051128,2.7315518487365438);
\draw [line width=0.4pt,color=black] (-0.8953798813051128,2.7315518487365438)-- (0.,3.230935550935551);
\draw [line width=0.4pt,color=black] (0.,3.230935550935551)-- (0.,2.);
\draw [line width=0.4pt,color=black] (0.,2.)-- (-1.176350188308613,1.7059124529228467);
\draw [line width=0.4pt,color=black] (0.,3.230935550935551)-- (0.6656059123409749,2.7359609200174266);
\draw [line width=0.4pt,color=black] (0.6656059123409749,2.7359609200174266)-- (0.8751724137931034,1.7082758620689655);
\draw [line width=0.4pt,color=black] (0.8751724137931034,1.7082758620689655)-- (0.,2.);
\draw [line width=0.4pt,color=black] (0.,2.)-- (0.,3.230935550935551);
\draw [line width=0.4pt,color=black] (-1.176350188308613,1.7059124529228467)-- (-0.23559785204630612,1.546900694008833);
\draw [line width=0.4pt,color=black] (-0.23559785204630612,1.546900694008833)-- (0.8751724137931034,1.7082758620689655);
\draw [line width=0.4pt,color=black] (0.8751724137931034,1.7082758620689655)-- (0.,2.);
\draw [line width=0.4pt,color=black] (0.,2.)-- (-1.176350188308613,1.7059124529228467);
\begin{scriptsize}
\draw [fill=black] (-4.,1.) circle (2.5pt);
\draw[color=black] (-4.137197888942202,1.3898073887022795) node {0100};
\draw [fill=black] (0.,2.) circle (2.5pt);
\draw[color=black] (-0.35,2.2345912848404637) node {0001};
\draw [fill=black] (3.,1.) circle (2.5pt);
\draw[color=black] (3.0165040400577405,1.3538591378027822) node {1000};
\draw [fill=black] (0.,6.) circle (2.5pt);
\draw[color=black] (0.4462041007436909,6.045105880187166) node {0010};
\draw [fill=black] (-1.176350188308613,1.7059124529228467) circle (2.5pt);
\draw[color=black] (-1.65,1.8391605249459944) node {0101};
\draw [fill=black] (0.8751724137931034,1.7082758620689655) circle (2.5pt);
\draw[color=black] (1.25,1.947005277644486) node {1001};
\draw [fill=black] (-0.23559785204630612,1.546900694008833) circle (2.5pt);
\draw[color=black] (-0.25,1.2) node {1101};
\draw [fill=black] (0.,3.230935550935551) circle (2.5pt);
\draw[color=black] (-0.35,3.4) node {0011};
\draw [fill=black] (-0.8953798813051128,2.7315518487365438) circle (2.5pt);
\draw[color=black] (-1.4,2.881659801031413) node {0111};
\draw [fill=black] (0.6656059123409749,2.7359609200174266) circle (2.5pt);
\draw[color=black] (1.05,2.91760805193091) node {1011};
\end{scriptsize}
\end{tikzpicture}
\caption[Representing the $D_4$ configuration as a cube in 3-point perspective.]{The $D_4$ configuration 
represented by a unit cube in 3-point perspective.
(Not visible: the back vertex point $(1:1:1:1)$, the center point $(1:1:1:2)$
(which except in characteristic 2 can be written $(\frac{1}{2}:\frac{1}{2}:\frac{1}{2}:1)$), 
the orthogonal lines through the point $(1:1:1:1)$ along the three back edges,
and the four main diagonals through opposite vertices of the cube and through $(1:1:1:2)$. 
Of the 10 visible points, the 9 points other than $(0:0:0:1)$
give a $(3,3)$-grid; the three grid lines in one ruling are shown with small dashes, each of the other three grid lines
is shown with dashes and dots.)}
\label{3ptPerspective}
\end{figure}

\begin{remark} 
Over $\mathbb C$, every nontrivial $\{3,b\}$-geproci set other than $Z_{D_4}$ is a grid 
\cite{POLITUS}.
However, if $\field=\mathbb Z/2\mathbb Z$, then $Z=\PP^3_K\subset\PP^3_{\overline K}$ is 
$(3,5)$-geproci 
\cite{kettinger2023},
indeed a $[3,5]$-half grid 
\cite{Ganger}.
We suspect but do not know that the only nontrivial nongrid $\{3,b\}$-geproci sets
are $Z_{D_4}$ (in any characteristic) and $Z=\PP^3_\field$ when $\field=\mathbb Z/2\mathbb Z$.
\qed\end{remark}

\subsection{Space curves containing geproci sets}\label{s: proof 1.1} 

In this section we will prove Theorem \ref{curve theorem}. 
We assume throughout this section that $\mathbb F=\mathbb C$
since our proof cites \cite{POLITUS} which assumes the ground field is
$\mathbb C$, but the first place in our proof where the characteristic might matter is 
Proposition \ref{curve prop}.

Recall a \emph{secant line} of $Z$ is a line joining two points of $Z$, regardless of whether that line contains any other points of $Z$.

The proof of Theorem \ref{curve theorem} involves several steps
stated as propositions below. The lemmas are tools used along the way. 

\begin{lemma}\label{lem: c is a or b}
Let $c$ be the minimal degree among all curves containing a given $(a,b)$-geproci set $Z$.
If $c\leq b$, then $c=a$ or $c=b.$
\end{lemma}

\begin{proof}
Let $Z$ be an $(a,b)$-geproci set so $a\le b$. 
If $Z$ is contained in a curve $C$ of degree $c<b$ then $Z$ is 
contained in a curve of degree $a$. Indeed, since $\overline{Z}$ is an $(a,b)$ complete intersection,
we have $\overline{Z}=A\cap B\subset H$ where $A\subset H$ is a curve of degree $a$ and 
$B\subset H$ is a curve of degree $b$, and the ideal of $\overline{Z}$ in the coordinate ring of $H$ is generated
by the form $F_A$ defining $A$ and the form $F_B$ defining $B$. Since $\overline{C}\subset H$ is
a curve of degree $c<b$ containing $\overline{Z}$, we see $a\leq c<b$ 
so $F_A$ divides the form $F_{\overline{C}}$ defining $\overline{C}$;
i.e., $A$ is contained in $\overline{C}$, so $C$ contains a curve $A'$ of degree $a$ which contains $Z$.
By minimality, $C=A'$ and $c=a$. 
\end{proof}

\begin{lemma} \label{Lem5}
Let $C$ be a reduced, irreducible curve of degree $b$ and suppose that  $Z \subset C$ is $\{a,b\}$-geproci. 
Let $\pi$ be projection from a point $R$ not on $C$ to a plane $H'$ not containing $R$, 
such that $\pi(C)$ has degree $b$ and $|\pi(Z)| = ab$.  
Then $\pi(Z)=F\cap\pi(C)$ for a curve $F\subset H'$ of degree $a$ so $\pi(Z)$ consists of smooth points of $\pi(C)$
(and hence $Z$ consists of smooth points of $C$)
and the elements defining $F$ and $\overline{C}$ in the ideal of $\pi(Z)$ on $H'$
generate the ideal.
\end{lemma}

\begin{proof}
First consider the case that $R$ is the general point $P$ and $H'=H$, so $\overline{Z}=\pi(Z)$. 
Since $C$ is irreducible, also $\overline{C}$ is irreducible. 
Since $\overline{Z}$ is a complete intersection of type $(a,b)$, we can find a curve 
$F\subset H$ of degree $a$ containing $\overline{Z}$ and not containing $\overline{C}$.
If $Z$ contained a point $P_1$ lying on the singular locus 
of $C$ then the projection $\overline{P_1}$ would lie on the singular locus of $\overline{C}$. 
Then the complete intersection of $\overline{C}$ and $F$ will not be reduced so consists of fewer than $ab$ points, 
but must contain $\overline{Z}$, which has $ab$ points. Thus the points of $\overline{Z}$ are smooth points of $\overline{C}$. 

Now specialize $R$ to any point off $C$ such that $\pi(C)$ has degree $b$ and $|\pi(Z)| = ab$. 
If $a<b$, then the dimension of the linear system of curves of degree $a$ containing $\pi(Z)$ is, by semicontinuity,
at least as big as the linear system of curves of degree $a$ containing $\overline{Z}$, so there are curves $F$ of degree $a$
which do not contain $\pi(C)$. If $a\geq b$, the dimension of the linear system of 
curves of degree $a$ containing containing $\pi(C)$ is $\binom{a-b+2}{2}$ which is the same as the linear system
of curves of degree $a$ containing containing $\overline{C}$, and the linear system of curves of degree $a$
containing $\overline{Z}$ is bigger than that, and by semicontinuity so is the dimension of the linear system of curves of degree $a$
containing $\pi{Z}$. Thus again there is a curve $F$ of degree $a$ containing $\pi(Z)$ but not containing $\pi(C)$.
Since $F$ and $\pi(C)$ have no components in common but have $ab$ distinct points in common, by B\'ezout's Theorem
$\pi(Z)=F\cap\pi(C)$. Thus the elements defining $F$ and $\overline{C}$ in the ideal of $\pi(Z)$ on $H'$
generate the ideal. If $z\in Z$ were a singular point of $C$, then $\pi(z)$ would be a singular point
of $\pi(C)$ and then $F\cap\pi(C)$ would not consist of $ab$ distinct points, so we see that the points of
$\pi(Z)$ are smooth points of $\pi(C)$.
\end{proof}

The first step considers when $Z$ lies on an irreducible curve $C$. 

\begin{proposition} \label{thm: pts on irred curve}
Let $C \subset \PP^3$ be a reduced, irreducible, nondegenerate curve of degree $b$ and let $Z \subset C$ be a reduced set of points. 
Then $Z$ is not $\{a,b\}$-geproci  for any $a$.  
\end{proposition}

\begin{proof}
Assume $Z = \{ P_1,\dots,P_{ab} \}\subset C$ is $\{a,b\}$-geproci. 
Let $D$ be the union of the secant lines of $Z$. For any point $R \in \mathbb P^3$ 
we will denote by $\pi_R$ the projection from~$R$ to some plane not containing $R$.

Take a general plane $H'$ through $P_1$. Since $P_1$ is a smooth point of $C$ by Lemma \ref{Lem5},
$H'$ meets the curve $C$ transversely in 
a set $Y$ of $b$ noncollinear points. Let $Q$ be a point of $Y$ different from $P_1$,
hence $Q$ is a general point of $C$. Consider the line $\lambda$ joining $P_1$ 
and $Q$ and fix a general point $\Lambda$ in $\lambda$.  

We claim that $\pi_\Lambda$ satisfies the assumptions in Lemma \ref{Lem5}. 
Notice that $\Lambda$ does not lie on $D$. 
Thus $\pi_\Lambda(Z)$ is a set of $ab$ distinct points. 
Now we confirm that $\pi_\Lambda(C)$ is a curve of degree $b$. 
This follows if we prove that the map $\pi_\Lambda$ restricted to $C$ is birational (generically one to one). 
It suffices to find a line through $\Lambda$ containing only one point of $C$. 
The existence of the line is clear since there is a point $Q' \in Y$ outside the line $\lambda$ and, 
since $\Lambda$ is general in $\lambda$, the line through $Q'$ and $\Lambda$ misses the remaining points 
of $Y$ and hence the remaining points of $C$ (since $\lambda$ is contained in $H'$).

But $\pi_\Lambda(C)$ is singular at $\pi_\Lambda(P_1)=\pi_\Lambda(Q)$, contrary to
Lemma \ref{Lem5} which asserts that the points of $\pi(Z)$ are smooth points of $\pi(C)$.
Thus $Z$ cannot be $\{a,b\}$-geproci.     
\end{proof}

Now we start considering reducible curves $C$. 

\begin{lemma} \label{key lemma}
Let $Z$ be $\{a,b\}$-geproci and assume $C$ is a curve of least degree $c$
containing $Z$ and that $c\leq\max(a,b)$.
Let $C_1, \dots, C_k$ be the irreducible components of $C$.
Set $c_i = \deg C_i$ (so $\sum c_i = c$) and $Z_i = Z \cap C_i$. Then $C$ is reduced and:

\begin{itemize}

\item[(a)] $Z$ lies in the smooth part of $C$; in particular, $Z_i$ lies in the smooth part of $C_i$ and 
\[
Z \cap (C_i \cap C_j)  = \emptyset \ \ \ \hbox{ for all } 1 \leq i < j \leq k;
\]

\item[(b)] $Z_i$ is $(a,c_i)$-geproci;

\item[(c)] $\overline{Z}$ is the complete intersection of $\overline{C}$ with a curve $F\subset H$ of degree $a$.

\end{itemize}
\end{lemma}

\begin{proof}
Minimality of $c$ means $C$ is reduced.
By Lemma \ref{lem: c is a or b} and symmetry we may assume $c=b$ (we do not assume either $a\leq b$ or $b\leq a$).
Now, $\overline{Z}$ is the reduced complete intersection of a curve $D_a$ and 
a curve $D_b$ of degrees $a$ and $b$ respectively, and $\overline{C}$ is 
a curve of degree $b$ containing $\overline{Z}$. 

If $D_a$ and $\overline{C}$ have no common component 
(as is the case if $\overline{C}=D_b$) then the result is immediate: 
$\overline{Z}$ is the reduced complete intersection of $D_a$ and $\overline{C}$, so 
$\overline{Z_i}$ is the reduced complete intersection of $D_a$ and $\overline{C_i}$, 
and $D_a$ avoids the singular locus of $\overline{C}$. This gives all of the claims in the lemma. 

The only danger is that $D_a$ contains one of the irreducible curves $\overline{C_i}$ 
as a component even though $D_b$ has no component in common with $D_a$. 
First assume $b < a$. Since $\overline{Z}$ is a complete intersection of 
type $(b,a)$ with $b<a$, there is a unique curve in $H$ of degree $b$ containing $\overline{Z}$,
hence $\overline{C}=D_b$ has no component in common with $F=D_a$.

Now assume $b=a$. Since $\overline{Z}$ is a complete intersection of 
type $(a,b)$ with $a=b$, the curves in $H$ of degree $a=b$ containing $\overline{Z}$
comprise a linear pencil and any two distinct members of the pencil
have intersection $\overline{Z}$. But $\overline{C}$, $D_a$ and $D_b$
all are members of the pencil. If all three are different or if $\overline{C}=D_b$, 
then $\overline{C}$ and $F=D_a$ have no common components. 
If $\overline{C}=D_a$, then $\overline{C}$ and $F=D_b$ have no common components.

Finally assume $b>a$. Since $\overline{Z}$ is a complete intersection of 
type $(b,a)$ with $b>a$, there is a unique curve in $H$ of degree $a$ containing $\overline{Z}$,
and every curve in $H$ of degree $b$ either contains $D_a$ or
meets $D_a$ exactly at $\overline{Z}$.
If $D_a\subset \overline{C}$, then let $C'$ be the components of $C$ which
project to $D_a$. Since the projection is general and $D_a$ contains $\overline{Z}$,
we must have $Z\subset C'$ so by minimality we have $c=a<b$, which is a contradiction.
I.e., $D_a\subset \overline{C}$ cannot happen.
\end{proof}

The next proposition is an immediate corollary of Lemma \ref{key lemma}(b) and Proposition \ref{thm: pts on irred curve}.

\begin{proposition} \label{pts on curve}
Let $Z$ be an $\{a,b\}$-geproci set and assume $C$ is a curve of least degree $c$
containing $Z$ and that $c\leq\max(a,b)$. 
Let $C_1, \dots, C_k$ be the irreducible components of $C$. 
Then $C$ is reduced and each $C_i$ is a plane curve.  
\end{proposition}

We next consider the situation where $Z$ is nondegenerate and lies on a curve consisting of plane curves. 

\begin{proposition} \label{union of plane curves}
Let $Z$ be a nondegenerate $\{a,b\}$-geproci set and assume $C$ is a curve of least degree $c$
containing $Z$ and that $c\leq\max(a,b)$. 
Let $C = C_1 \cup \dots \cup C_k$, $k\geq 2$, be the irreducible components of $C$
and assume each is a plane curve. Then $\deg(C_i)=1$ for all $i$.
\end{proposition}

\begin{proof}
We have $a,b>1$ since $Z$ is nondegenerate.
As usual we may assume $c=b$.
By Lemma \ref{key lemma}, $Z$ must lie on the smooth part of $C$ (so each point of $Z$ lies 
on a unique $C_i$) and $Z_i=C_i\cap Z$ 
is a planar complete intersection of type $(a,b_i)$ where $b_i=\deg(C_i)$
(since $C_i$ is planar and $Z_i$ is $\{a,b_i\}$-geproci by the lemma).

Say $b_1>1$. Then there is a unique plane $H_1$ containing $C_1$
and, since $a>1$, $Z_1$ spans that plane. 
Since $Z$ is nondegenerate, there is a component, say $C_2$,
with $Z_2$ not contained in $H_1$. 

Now $\overline{Z}$ is the reduced complete intersection in $H$
of $\overline{C}$ with some curve $D_a\subset H$ of degree $a$. 
Successively deleting sets $Z_i$ from $Z$ preserves the 
geproci property by 
\cite[Lemma 4.5]{POLITUS}. 
Thus we reduce to the case $C = C_1 \cup C_2$ of two curves 
of degrees $b_1>1, b_2\geq 1$ respectively and $Z=Z_1\cup Z_2$
is nondegenerate.

Specialize the projection. 
Let $P_1 \in Z_1$; if $b_2>1$, we may assume $P_1$ is not in 
the plane spanned by $C_2$.
Let: $Q$ be a general point of $C_2$,
$\Lambda$ a general point of the line $\lambda$ spanned by $P_1$ and $Q$, and
$\pi$ the projection from $\Lambda$ to the plane $H_1$.

Then $\pi(C_1)=C_1$ and $\pi(C_2)$ are distinct irreducible 
curves of degrees $b_1, b_2$ respectively,  and $\pi(Z_1)=Z_1$ and $\pi(Z_2)$ 
are disjoint sets of $ab_1$ and $ab_2$ points respectively and $P_1=\pi(P_1) \in \pi(C_2)$
is the only point of $Z_1$ on $\pi(C_2)$. 
Thus $\pi(Z)$ contains $ab_2+1$ points of $\pi(C_2)$, 
so $\pi(C_2)$ is a fixed component of $[I_{\pi(Z)}]_a$ on $H_1$.  
The other $ab_1-1$ points of $\pi(Z)$ lie on $\pi(C_1)=C_1$.
But every curve $A\subset H_1$ defined by a form in $[I_{\pi(Z)}]_a$ contains
$\pi(C_2)$, so $A$ meets $C_1$ in at least $(ab_1-1) + (b_1b_2)$ points,
and since $b_1>1$ we have $(ab_1-1) + (b_1b_2)>ab_1$ so
$C_1$ is also a fixed component of $[I_{\pi(Z)}]_a$.
This means $[I_{\pi (Z)}]_a$ consists of all forms of degree $a-b_1-b_2$
times the forms defining $C_1$ and $\pi(C_2)$, hence
\[
\dim [I_{\pi (Z)}]_a = \binom{a-b_1-b_2+2}{2}.
\]  
This contradicts semicontinuity since 
\[
\dim [I_{\overline{Z}}]_a = \binom{a-b_1-b_2+2}{2}+1.
\eqno\qed
\]  
\phantom\qedhere
\end{proof}

\begin{proposition}\label{curve prop}
Let $Z$ be a nondegenerate $\{a,b\}$-geproci set and assume $C$ is a curve of least degree $c$
containing $Z$ and that $c\leq\max(a,b)$. Then $c=a$ or $c=b$, and $C$ is a 
union of $c$ disjoint lines each of which contains the same number, $ab/c$, of points of $Z$. 
\end{proposition}

\begin{proof}
If $a=1$ or $b=1$, then the points of $Z$ are collinear hence $Z$ is degenerate.
Thus we may assume $a>1$, $b>1$.
If $a=2$ or $b=2$ (say $b=2, a\geq 2$), then $Z$ consists of $a$ points on each of 2 skew lines
by 
\cite[Proposition 4.8]{POLITUS}; 
thus $c=2$. If $a>2=b$, then the only curve of degree 2
containing $Z$ consists of the two skew lines, each containing $a$ points of $Z$.
If $a=2=b$, then $Z$ consists of 4 linearly general points, so the only curves of degree
2 which contain $Z$ are pairs of skew lines, each containing 2 points of $Z$.

Now consider the case that $a,b\geq 3$.
By Lemma \ref{lem: c is a or b}, $c=a$ or $c=b$; we may assume $c=b$.
By Proposition \ref{union of plane curves}, $C$ is a union of lines, say $L_1\cup\cdots\cup L_b$.
By Lemma \ref{key lemma}, the points of $Z$ are smooth points of $C$.
But $\overline{Z}\subset H$ is a complete intersection $\overline{Z}=D_1\cap D_2$ of curves $D_1,D_2\subset H$
with $\deg(D_1)\leq \deg(D_2)$ and $\{\deg(D_1), \deg(D_2)\}=\{a,b\}$.
So if $\deg(\overline{C})=\deg(D_1)<\deg(D_2)$, then $\overline{C}=D_1$, so
$\overline{C}\cap D_2=\overline{Z}$ and every line in $\overline{C}$
contains $a$ points of $\overline{Z}$ hence every line in $C$ contains $a$ points of $Z$.
If $\deg(\overline{C})=\deg(D_1)=\deg(D_2)$, then $\overline{C}$, $D_1$ and $D_2$ 
are all members of the pencil of curves in $H$ defining $\overline{Z}$ so
$\overline{C}\cap D_i=\overline{Z}$ either for $i=1$ or $i=2$ or both and again
every line of $C$ contains $a$ points of $Z$.
If $\deg(\overline{C})>\deg(D_1)$, we have by Lemma \ref{lem: c is a or b}
that $\deg(\overline{C})=\deg(D_2)>\deg(D_1)$ and by minimality 
$D_1$ cannot be contained in $\overline{C}$, so (since 
$\overline{Z}$ is an $(a,b)$ complete intersection) 
we have $\overline{C}\cap D_1=\overline{Z}$.
Thus every line of $\overline{C}$ contains $a$ points of $\overline{Z}$
so every line of $C$ contains $a$ points of $Z$.

Thus each line of $C$ contains exactly $a$ points of $Z$. 
Since $Z$ is nondegenerate, we can choose three of the lines, $L_1,L_2,L_3$ say,
which are not coplanar. 
By 
\cite[Lemma 4.5]{POLITUS},
we can remove the lines $L_4,\ldots,L_b$ one at a time
to reduce to the case that $Z'=Z\cap (L_1\cup L_2\cup L_3)$ is $\{3,a\}$-geproci.
Now by 
\cite[Proposition 4.14]{POLITUS}
a $\{3,a\}$-geproci set is a grid, 
so in particular these three lines are skew to each other. Likewise, $L_1,L_2,L_i$ are skew for all $i>2$,
and now $L_2,L_3,L_j$ are skew for all $j>3$, etc. Thus all of the lines are skew to each other.
Thus $C$ consists of $b$ skew lines and $Z$ consists of $a$ points on each line.
\end{proof}

Now we prove Theorem \ref{curve theorem}. Recall that overbar denotes
the image of projection from a general point to a given plane $H$.

\begin{proof}[Proof of Theorem \ref{curve theorem}]
Since $Z$ is nondegenerate, we have $2\leq a\leq b$. 

(a) By Proposition \ref{curve prop}, $C_1$ is a union of $a$ skew lines
with $b$ points on each line. Let $L$ be any of these lines and let $C$ be any component of $C_2$.
Since $\overline{Z}=\overline{C_1}\cap\overline{C_2}$,
we see $\overline{L}\cap\overline{C}$ consists of $\deg(C)$ points of $\overline{Z}$
and hence $L\cap C$ consists of $\deg(C)$ points of $Z$.
Thus there is a plane containing $L$ which contains these $\deg(C)$ points
plus at least one other point of $C$, hence the plane contains $C$. But $a\geq 2$ so $C$ is also contained 
in a plane containing another line $L'$ of $C_1$, hence $C$ is itself a line,
transversal to the lines in $C_1$. Thus $C_2$ is a union of $b$ lines,
each of which is transversal to the lines in $C_1$. Thus $Z$ is an $(a,b)$-grid.

If $2=a=b$, we can think of $Z$ as the four points of a tetrahedron
so there are 3 ways to choose $C_1$ and for each there are two ways to choose $C_2$.
So neither $C_1$ nor $C_2$ is unique.

If $2=a<b$, then the only lines through $b$ points of $Z$ are the two lines
comprising $C_1$, so $C_1$ is uniquely determined. But there are $b!$ ways to
pair up the $b$ points on one component of $C_1$ with the $b$ points
on the other component. Each pair in such a pairing determines a line, the resulting $b$ lines 
for a given pairing give a choice of $C_2$ and distinct pairings give different curves 
$C_2$. So in this case $C_1$ is unique but $C_2$ is not.

If $2<a=b$, then again there is a unique smooth quadric $Q$ containing $C_1$,
and the lines comprising $C_1,C_2$ are the ruling lines of $Q$ which meet $Z$,
so the partition of the lines into the two curves $C_1$ and $C_2$ of degree $a=b$ is unique
but which one we call $C_1$ and which we call $C_2$ is arbitrary.

If $2<a<b$, then there is a unique smooth quadric $Q$ containing $C_1$,
and the lines comprising $C_1,C_2$ are the ruling lines of $Q$ which meet $Z$.
Thus $C_1$ and $C_2$ are uniquely determined.

(b) By \cite[Proposition 4.8]{POLITUS},
a nondegenerate $(a,b)$-geproci set is an $(a,b)$-grid
if $a=2,3$ which means $Z$ would have type (1).
Thus we must have $4\leq a\leq b$.

Let $C$ be a space curve of minimum degree $c$ containing $Z$.
Since $Z$ is of type (1), there is a curve $C_1$ containing $Z$ with $\deg(C_1)\leq b$, so
$c\leq b$ and by Lemma \ref{lem: c is a or b} we have either $c=a$ or $c=b$.
By Proposition \ref{curve prop}, $C$ is a union of skew lines with the same 
number $d$ of points of $Z$ on each line, so $cd=|Z|=ab$ hence $\{c,d\}=\{a,b\}$ 
and $Z$ is $[d,c]$-geproci. Since $Z$ is not of type (1), $Z$ is not an $(a,b)$-grid, so
$Z$ is a $[d,c]$-half grid.

If $a<b$ and $\deg(C)=a$, then $\overline{C}$ is uniquely determined by
$\overline{Z}$ since it is defined by the unique
generator of the ideal of $\overline{Z}$ of smaller degree $a$. 
And $\overline{C}$ determines $C$: if $\ell$ is a line in $\overline{C}$,
then $\ell\cap \overline{Z}$ is a set of collinear points, hence the  
points of $Z$ projecting to these collinear points are collinear
and thus lie on a unique line $L$ which projects to $\ell$. 
Thus the components of $\overline{C}$ determine the
components of $C$, so $C$ is unique.

If $a=b$, there is a pencil of curves of degree $a$ which contain $\overline{Z}$,
one of which is $\overline{C}$.
If there were a second curve $Z\subset C'$ of degree $a$ it would also be a union of lines,
and $\overline{C}$ and $\overline{C'}$ would both be in this pencil,
so $\overline{Z}=\overline{C}\cap \overline{C'}$ and hence
$Z=C\cap C'$, which makes $Z$ of type (1). Thus $C$ is unique
in this case.

The remaining case is $a<b=\deg(C)$. 
Example \ref{Example:D4} explicitly gives a quartic containing
the $[3,4]$-half grid $Z_{D_4}$. But the symmetries of the cube
shown in Figure \ref{3ptPerspective} take that set of four lines to
others so $C$ need not be unique in this case. 
(Lack of uniqueness also happens in positive characteristics.
The sets $Z_{D_4}$ in positive characteristics give some examples.
For others, consider the points $Z=\PP^3_F$ over the field $F$ of $q$ elements,
where $q$ is any power of the characteristic. 
This is $[q+1,q^2+1]$-geproci in $\PP^3_{\overline{F}}$;
it is a half grid with respect to the $q^2+1$ skew lines
coming from the Hopf fibration (see \S\ref{Sec:CombQuestsl}). 
Lines through two or more points of $Z$ are defined over
$F$ and so contain exactly $q+1$ points of $Z$, thus the minimum degree for a 
curve which is a union of skew lines and which contains $Z$ is $q^2+1$
so $c=q^2+1$. But $Z$ is also a $[q+1,q^2+1]$-half grid 
with respect to many other choices of $q^2+1$ skew lines: just apply linear automorphisms
of $\PP^3_F$ to get others.)
\end{proof}

\begin{example}\label{Penrose et alia}
Although most of the known nontrivial nongrid $\{a,b\}$-geproci sets over $\mathbb C$ 
are half grids, and so the minimum degree of a curve on which 
the half grid lies is at most $\max(a,b)$, not every nontrivial nongrid $\{a,b\}$-geproci set is a half grid.
For example, 
\cite{POLITUS}
gives
examples of geproci sets over the complex numbers of 40, 60 and 120 points 
which are nontrivial and neither grids nor half grids.
Additional examples in positive characteristics are given in 
\cite{kettinger2023}.
However, the three nontrivial complex nongrid non-half grids are unions of half grids.
For example, the 40 point example $Z_{40}$, dubbed the Penrose configuration in 
\cite{POLITUS}, 
is $(5,8)$-geproci. Thus we conclude from the work in this section 
that $Z_{40}$ does not lie on any curve of degree 8 or less.
But $Z_{40}$ consists of 10 sets of 4 collinear points, and these 10 sets can be partitioned into
two sets of 5 sets of 4 collinear points each, each of which is a $[4,5]$-half grid. By Proposition \ref{Prop:[a,b]geprociImpliesCollCompl}, 
every $[a,3]$-geproci set is a grid, so every half grid is a union of grids
(in every characteristic), and (over $\mathbb C$ at least)
the known examples of nontrivial nongrid non-half grids are unions of half grids (and thus also of grids).
So understanding half grids may be key to understanding nontrivial nongrid non-half grids.
\qed \end{example}

\begin{definition}
    Let $X\subset \PP^2$ be a 0-dimensional scheme with Hilbert function $H(X,t) = \dim [R/I_X]_t$, where $R=\field[x,y,z]$ is the homogeneous coordinate ring for $\PP^2$. Let $L$ be a general linear form. Then the {\it $h$-vector} of $X$, denoted $h_X(t)$, is the Hilbert function of the artinian reduction of $X$, i.e. of $R/(I_X,L)$, and we have
    \[
    h_X(t) = \Delta H(X,t) = H(X,t) - H(X,t-1).
    \]
\end{definition}

\begin{example}\label{ClassicalExamples}
Our computations suggest further connections between the  
geproci sets of type (3) 
and half grids.  
\begin{enumerate}
  \item {\bf The $H_4$ configuration}. Let $\varepsilon$ be a primitive root of $t^5-1=0.$ Consider the following 12 lines
    \begin{itemize}
        \item  $L_j$ defined by $\varepsilon^jx - y= \varepsilon^jz - w=0$ for $j=0,\ldots, 4$ (on the quadric  $xw-yz=0$);
        \item  $L_j'$ defined by  $\varepsilon^j\eta 
        x - y= \varepsilon^jz - \eta w=0$ for $j=0,\ldots, 4,$ with $\eta=\varepsilon^4+\varepsilon-1$ (on the quadric $\eta^2xw-yz=0$);
        \item $\ell_1$ and $\ell_2$  defined by $y=z=0$ and $x=w=0$ respectively.
    \end{itemize}
Let $\mathcal L_1=\{L_j\}\cup \{\ell_1\}$, $\mathcal L_2=\{L_j'\}\cup \{\ell_2\}$ and $\mathcal L=\mathcal L_1\cup \mathcal L_2$. 
    
We note that the lines $T_1$ and $T_2$ defined by $x=y=0$ and $z=w=0$ are transversals for $\mathcal L$. The $C_{\mathcal L}-$orbit of a point $p$ not in $T_1$ and $T_2$ is finite and $|[p]_{\mathcal L}|=120$. The set of points $Z=[p]_{\mathcal L}$ has $h$-vector  
$$h_Z=(1,  3,  6,  10,  15,  21,  18,  14,  14,  14,  3,  1).$$
Set $p_1=(1:0:0:-1), p_2=(1:0:0:1)\in \ell_1 $, $q_1=(0:1:-1:0), q_2=(0:1:1:0)\in \ell_2 $ then we note that
$Z_i=[p_i]_{\mathcal L_i}$, and 
$Z_i'=[q_i]_{\mathcal L_i}$, for $i=1,2,$ are all projectively equivalent to the $(5,6)$-geproci set in the standard construction \cite{POLITUS} 
(see \S \ref{stConstr})
and $Z_1\cup Z_2 \cup Z_1'\cup Z_2'$ is projectively equivalent to $Z$. Moreover $Z_1\cup Z_1'$ and $Z_2\cup Z_2'$ are  projectively equivalent to the $H_4$ configuration.

\item {\bf The $F_4$ configuration.} Let $\varepsilon$ be a primitive root of $t^3-1=0.$ The following construction 
is similar to that in item $(1)$, starting from the following 8 lines:
    \begin{itemize}
        \item  $L_j$ defined by $\varepsilon^jx - y= \varepsilon^jz - w=0$ for $j=0,1, 2$ (on the quadric  $xw-yz=0$);
        \item  $L_j'$ defined by  $\varepsilon^j\eta 
        x - y= \varepsilon^jz - \eta w=0$ for $j=0,1, 2,$ with $\eta=\varepsilon^2+\varepsilon-1=-2$ (on the quadric $\eta^2xw-yz=0$);
        \item $\ell_1$ and $\ell_2$  defined by $y=z=0$ and $x=w=0$ respectively.
    \end{itemize}
Let $\mathcal L_1=\{L_j\}\cup \{\ell_1\}$, $\mathcal L_2=\{L_j'\}\cup \{\ell_2\}$ and $\mathcal L=\mathcal L_1\cup \mathcal L_2$. 
    
We note that the lines $T_1$ and $T_2$ defined by $x=y=0$ and $z=w=0$ are transversals for $\mathcal L$. The $C_{\mathcal L}-$orbit of a point $p$ not in $T_1$ and $T_2$ is finite and $|[p]_{\mathcal L}|=48$. The set of points $Z=[p]_{\mathcal L}$ has $h$-vector  
$$h_Z=(1,  3,  6,  10,  13,  11,  3,  1).$$
We checked that no five of these lines are contained in the same quadric.

Set $p_1=(1:0:0:-1), p_2=(1:0:0:1)\in \ell_1 $, $q_1=(0:1:-1:0), q_2=(0:1:1:0)\in \ell_2 $ then we note that
$Z_i=[p_i]_{\mathcal L_i}$, and 
$Z_i'=[q_i]_{\mathcal L_i}$, for $i=1,2,$ are all projectively equivalent to the $D_4$ configuration and $Z_1\cup Z_2 \cup Z_1'\cup Z_2'$ is projectively equivalent to $Z$. Moreover $Z_1\cup Z_1'$ and $Z_2\cup Z_2'$ are  projectively equivalent to the $F_4$ configuration. The lines in $\mathcal L$ are then 3-secant to the $F_4$ configuration.  The $F_4$ configuration can be covered by 8 skew 3-point lines in 64 different ways. But only for 8 of them does the groupoid $C_{\mathcal L}$ have a finite orbit producing a $(6,8)$-geproci set as above.

    \item The idea of looking at the secant lines can be also applied to {\bf the Penrose configuration} $Z_{40}$. It has 90 4-point lines and 459 sets of 10 skew lines covering it. Only for 27 of them the groupoid $C_{\mathcal L}$ has a finite orbit producing a (8,10)-geproci set $Z_{80}$ with $h$-vector
$$h_{Z_{80}}=(1,  3,  6,  10,  15,  21,  13,  11).$$

One such example is the following set of 10 lines $\mathcal L=\mathcal L_1\cup \mathcal L_2$   
$$ \begin{array}{lllll}
 \mathcal L_1= &    w=y +\varepsilon^2z=0,&  y +z -\varepsilon^2w=  x -\varepsilon^2z -\varepsilon^2w=0, \\
  &     x=y -\varepsilon^2z=  0,& y -\varepsilon z -w=  x -\varepsilon^2 z -w=0,\\
  &     & y -\varepsilon z -\varepsilon^2w=  x -z -w=0,\\
      \\

 \mathcal L_2= &             z=  x +\varepsilon w=0,& y -\varepsilon z +\varepsilon w=  x +\varepsilon z -w=0, \\ &  y=  x -\varepsilon w=0,  & y +z +\varepsilon w=  x -z -\varepsilon^2 w=0, \\ &&  y +z -w=  x +\varepsilon z -\varepsilon^2 w=0\\
\end{array}$$
where $\varepsilon$ is a root of $t^2-t+1=0.$
We checked that $\mathcal L$ has no transversals and at most 4 of the lines in $\mathcal L$ are contained on the same quadric. 
Also, for suitable points $p_1,p_2,q_1,q_2\in Z_{80}$,  we note that
$Z_i=[p_i]_{\mathcal L_i}$, and 
$Z_i'=[q_i]_{\mathcal L_i}$, for $i=1,2,$ are all projectively equivalent to the half Penrose configuration and $Z_1\cup Z_2 \cup Z_1'\cup Z_2'$ is projectively equivalent to $Z_{80}$. Moreover $Z_1\cup Z_1'$ and $Z_2\cup Z_2'$ are projectively equivalent to the $Z_{40}$ configuration. 
\end{enumerate}
\end{example}

\subsection{The standard construction}\label{stConstr}
The standard construction of \cite{POLITUS} gives an infinite class of examples
of nondegenerate nongrid half grids over the complex numbers.
In this section we extend this construction to finite characteristics and then classify the finite abelian groups
which arise as groups $G_{\mathcal L}$ for finite sets of skew lines ${\mathcal L}$.
We also answer a question raised in \cite{Fields}. 

Given a finite multiplicative subgroup $G\subset \mathbb C^*$, the standard construction of \cite{POLITUS} is 
a construction of two geproci sets $Z_0$ and $Z_\infty$ of $|G|$ points on each of $|G|+1$ skew lines, 
and, when $(-1)^{|G|}=1$, it also gives a related geproci set $Z_{0\infty}$
of $|G|$ points on each of $|G|+2$ skew lines. We show here that this construction generalizes almost verbatim
to the case of a finite multiplicative subgroup $G\subset {\mathbb F}^*$ for any field $\mathbb F$
(although for simplicity we assume here that $\mathbb F$ is algebraically closed).
The main point to note is that if $\mathbb F$ has characteristic 0, then 
for each positive integer $m$ there is a unique multiplicative subgroup $G\subset \mathbb F^*$
of order $m$, and it is cyclic, generated by a primitive $m$th root of 1.
If $\mathbb F$ has positive characteristic $p$, there is a multiplicative subgroup $G\subset\mathbb F^*$
of order $m$ if and only if $p$ does not divide $m$, but in that case $G$ is unique, cyclic  
and generated by a primitive $m$th root of 1. When $m=1,2$ the geproci sets are grids,
but for $m>2$, the sets $Z_0$ and $Z_\infty$ are $(m,m+1)$-half grids whose half grid lines have two transversals
and the set $Z_{0\infty}$ is an $(m,m+2)$-half grid whose half grid lines have two transversals.
In each case the half grid points comprise a single orbit with respect to the half grid lines.

When $\mathbb F$ has positive characteristic $p$ and $A\subset\mathbb F$ is any finite
additive subgroup of $\mathbb F$ (hence $A$ is a finite dimensional $\mathbb Z/p\mathbb Z$-vector space
and hence has order $p^s$ for some $s\geq0$), 
a modification of the construction gives a geproci set $Z$ of $|A|$ points on each of $|A|+1$ skew lines.
When $|A|<3$, the set $Z$ is a grid, but for $|A|\geq3$ the set $Z$ is an $(|A|,|A|+1)$-half grid whose half grid lines
have a single transversal of multiplicity 2 and the set $Z$ is a single orbit with respect to the
half grid lines.

We begin with the case of a finite multiplicative subgroup $G\subset\mathbb F$.
The standard construction gives a set $S$ of $m$ points on each of $m+2$ skew lines
$\lambda_0, \lambda_\infty, V_0,\ldots,V_{m-1}$, as shown in Figure \ref{StConstFig2Trans}.
There are three cases: 
$\mathcal L=\mathcal L_0=\{\lambda_0, V_0,\ldots,V_{m-1}\}$ and $Z=Z_0=S\cap(\cup_{L\in\mathcal L_0}L)$;
$\mathcal L=\mathcal L_\infty=\{\lambda_\infty, V_0,\ldots,V_{m-1}\}$ and $Z=Z_\infty=S\cap(\cup_{L\in\mathcal L_\infty}L)$; and,
when $(-1)^m=1$ (hence when either $m$ is even or $p=2$),
$\mathcal L=\mathcal L_{0\infty}=\{\lambda_0, \lambda_\infty, V_0,\ldots,V_{m-1}\}$ and $Z=Z_{0\infty}=S$.
In each case we show $Z$ is a single groupoid orbit and an $(m,|\mathcal L|)$-half grid with respect to the lines $\mathcal L$.
(The original result, given in \cite{POLITUS}, was for $p=0$ and $\mathbb F=\mathbb C$, in which case
by \cite{Ganger} $Z$ is a single orbit with respect to $\mathcal L$ and the lines in $\mathcal L$
are fibers of the Hopf fibration.)

Here we define the line
$\ell_0$ by $y,z$; 
$\ell_\infty$ by $x,w$; 
$V_i$ for $0\leq i < m$ by $u^ix-y,u^iw-z$;
$\lambda_0$ by $x,z$; 
and $\lambda_\infty$ by $y,w$.
Note that $\ell_0, V_0$ and $\ell_\infty$ are the lines $L_1,L_2$ and $L_3$ in standard position (see Definition \ref{StPosDef}) and that
the standard transversals $T_1\colon x,y$ and $T_2\colon z,w$ are transversals for $\mathcal L_{0\infty}$.

Consider also the lines $H_j$, $0\leq j < m$, defined by $u^jx-w,u^jy-z$, and let
$q_{ij}=(1:u^i:u^{i+j}:u^j)$ be the point $V_i\cap H_j$. The lines $T_1, T_2$ and the lines $H_j$ 
all are lines in the same ruling of the quadric $Q\colon xz-yw$, shown as horizontal lines in Figure \ref{StConstFig2Trans}.
The lines $V_i$ are all in the other ruling on $Q$, shown as vertical lines in the figure.
The lines $\lambda_i$ are not contained in $Q$.
The points $q_{ij}$ form an $(m, m)$-grid whose grid lines are the lines $V_i$ and $H_j$.
The set $S$ consists of the points $q_{ij}$, $0\leq i,j<m$, together with the points $(0:1:0:-u^i)\in\lambda_0$, $0\leq i<m$,
and the points $(1:0:-u^i:0)\in\lambda_\infty$, $0\leq i<m$.

The proof given in \cite{POLITUS,Ganger} also works in positive characteristics.
For convenience we give it here. Take $P=(a:b:c:d)$ as the general point
of projection, where we regard $a,b,c,d$ as variables. The plane $V^*_i$ spanned by $V_i$ and $P$ is defined by
$F_i=(u^id-c)(u^ix-y)-(u^ia-b)(u^iw-z)$, and 
the plane $H^*_j$ spanned by $H_j$ and $P$ is defined by
$G_j=(u^jb-c)(u^jx-w)-(u^ja-d)(u^jy-z)$. Let $F=\prod_{i=0}^{m-1}F_i$ and let 
$G=\prod_{j=0}^{m-1}G_j$. Also let $\lambda_i^*$ be the plane spanned and
by $\lambda_i$ and $P$, so $\lambda_0^*$ is defined by $\Lambda_0=cx-az$ and
$\lambda_\infty^*$ is defined by $\Lambda_\infty=dy-bw$.
Then $M_0=F\Lambda_0$ is the cone over the lines in $\mathcal L_0$, 
$M_\infty=F\Lambda_\infty$ the cone for the lines in $\mathcal L_\infty$
and $M_{0\infty}=F\Lambda_0\Lambda_\infty$ the cone for the lines in $\mathcal L_{0\infty}$.

Note $F$ and $G$ define a pencil of cones of degree $m$ with vertex $P$ and
which all vanish at all of the points $q_{ij}$. For each point $q$ of $\PP^3_{\mathbb F}$ not
on any line through $P$ and a point $q_{ij}$, there is a unique element of the pencil vanishing at $q$.
This unique member turns out to be the same member 
$H_0$ of the pencil for all $q\in Z_0\cap\lambda_0$ (but $H_0$ does not vanish on all of $\lambda_0$).
In fact, $H_0=F-G$ when $(-1)^m=1$ and $H_0=F+G$ when $(-1)^m=-1$. 
Moreover, let $H_{0\infty}=H_\infty=F-G$; then $H_\infty$ vanishes at the points of $Z_\infty\cap\lambda_\infty$
(but not on all of $\lambda_\infty$); in this case $(-1)^m$ does not matter. 

Given this, $M_0$ and $H_0$ have no components in common (since $H_0$ does not vanish on any of the lines
in $\mathcal L_0$). Thus $M_0$ and $H_0$ define a complete intersection curve of degree $(m+1)m$;
but this complete intersection curve contains the lines through $P$ and each of the $(m+1)m$ points of $Z_0$,
hence this complete intersection curve is this set of $(m+1)m$ lines, so the projection $\overline{Z_0}$ of $Z_0$ from $P$ to a plane $H$ 
is the complete intersection defined on $H$ by $H_0|_H$ and $M_0|_H$. I.e., $Z_0$ is an
$(m,m+1)$-half grid with respect to the lines $\mathcal L_0$. Similarly, $Z_\infty$ is an
$(m,m+1)$-half grid with respect to the lines $\mathcal L_\infty$ and, when $(-1)^m=1$,
$Z_{0\infty}$ is an $(m,m+2)$-half grid with respect to the lines $\mathcal L_{0\infty}$.

To see that $H_i=0$ does not contain $\lambda_i$, $i=0,\infty$, just plug $0100$ and $1000$ into $H_0$ and $H_\infty$
respectively. For $H_0$ at $0100$ we get $\prod_i (-(u^id-c)) \pm \prod_j (-(u^ja-d)u^j)$ which is clearly not 0 (it is the sum or difference
of two nonzero polynomials, one divisible by $c-d$, the other not). Similarly, $H_\infty$ at $1000$ is not 0.
To see that $H_0$ vanishes at the points of $Z_0\cap\lambda_0$ plug $(0,1,0,-u^k)$ into $H_0$.
We get 
$$F(0,1,0,-u^k)=\prod_i((u^id-c)(-1)+(u^ia-b)(u^{i+k}))=\prod_i(u^{2i+k}a-u^{i+k}b+c-u^id)$$
and
$$G(0,1,0,-u^k)=\prod_j((u^jb-c)(u^k)-(u^ja-d)(u^j))=(-1)^m(u^{km})\prod_j(u^{2j-k}a-u^jb+c-u^{j-k}d)$$
$$=(-1)^m\prod_i(u^{2i+k}a-u^{i+k}b+c-u^id)$$ 
where the last equality comes from substituting $i+k$ in for $j$, which merely permutes the factors,
and using $u^{km}=1$. Thus we see that $F(0,1,0,-u^k)-G(0,1,0,-u^k)=0$ if $(-1)^m=1$ 
and $F(0,1,0,-u^k)+G(0,1,0,-u^k)=0$ if $(-1)^m=-1$.

A similar argument works for $Z_\infty$, using the points $(1:0:-u^k:0)\in\lambda_\infty$: 
$$F(1,0,-u^k,0)=\prod_i((u^id-c)u^i-(u^ia-b)u^k)=\prod_i(-u^i)\prod_i(u^ka-u^{k-i}b+c-u^id)$$
and
$$G(1,0,-u^k,0)=\prod_j((u^jb-c)u^j-(u^ja-d)u^k))=\prod_j(-u^j)\prod_j(u^ka-u^jb+c-u^{k-j}d)$$
$$=\prod_j(-u^j)\prod_i(u^ka-u^{k-i}b+c-u^id)$$
where for the last equality for $G$ we make the substitution $j=k-i$.
So we see that $F(1,0,-u^k,0)-G(1,0,-u^k,0)=0$.

Since $Z_0$, $Z_\infty$ and $Z_{0\infty}$ are half grids with respect to their respective half grid lines $\mathcal L$,
each is a union of $C_{\mathcal L}$ orbits. To confirm each is a single orbit, consider the groupoid map $f_i\colon V_0\to \lambda_0$, $0<i<m$,
where $f_i(q)$ is the point where the plane spanned by $q\in V_0$ and $V_i$ meets $\lambda_0$. 
Thus the line $L_{iq}$ spanned by $q$ and $f_i(q)$ is a transversal for $\lambda_0,V_0$ and $V_i$.
If $f_i(q)=f_j(q)$, then $L_{iq}=L_{jq}$ is a transversal for $\lambda_0,V_0,V_i,V_j$. If $i\neq j$, then
$L_{iq}=L_{jq}$ meets $Q$ at points of all three lines $V_0,V_i,V_j$, hence $L_{iq}\subset Q$.
Thus $\lambda_0$ meets $Q$ at three points, $\lambda_0\cap T_1,\lambda_0\cap T_2, \lambda_0\cap L_{iq}$,
so $\lambda_0\subset Q$, which is not true. Therefore, $f_i(q)=f_j(q)$ implies $i=j$, hence 
$f_1(q),\ldots,f_{m-1}(q)$ are distinct, so the groupoid orbit of a point $q\in V_0$ contains at least $m-1$ points
on each line. The only way $Z_0$, $Z_\infty$ or (when $m$ is even) $Z_{0\infty}$ could fail to be a single orbit
is if they included an orbit consisting of 1 point on each line, but such an orbit must be on a transversal for the
corresponding set of lines
and none of the points of $Z_0$, $Z_\infty$ or $Z_{0\infty}$ lie on either transversal.
Since there are two transversals, the group of the groupoid in each case is a cyclic group of order $m$.

\begin{figure}[h]
\begin{tikzpicture}[x=1cm,y=1cm]
\clip(-11.09,3) rectangle (0,6.5);
\draw [line width=1pt, dotted] (-7.55,3.3) -- (-4.3,6);
\draw [fill=white,color=white] (-5.9,4.7) circle (12pt);
\draw [line width=1pt, dashed] (-7.7,6.05) -- (-4.45,3.3);
\draw [fill=white,color=white] (-6,4.6) circle (5pt);
\draw [line width=1pt] (-3,6)-- (-3,3);
\draw [line width=1pt] (-7,6)-- (-7,3);
\draw [line width=1pt] (-5,6)-- (-5,3);
\draw [line width=1pt] (-8,5.44)-- (-2,5.44);
\draw [line width=1pt] (-7.96,4.58)-- (-2,4.58);
\draw [line width=1pt] (-7.99,3.77)-- (-2,3.77);
\draw [fill=black] (-7,5.45) circle (2.5pt);
\draw [fill=black] (-7,4.6) circle (2.5pt);
\draw [fill=black] (-7,3.78) circle (2.5pt);
\draw [fill=black] (-5,5.45) circle (2.5pt);
\draw [fill=black] (-5,4.6) circle (2.5pt);
\draw [fill=black] (-5,3.78) circle (2.5pt);
\draw [fill=black] (-3,5.45) circle (2.5pt);
\draw [fill=black] (-3,4.6) circle (2.5pt);
\draw [fill=black] (-3,3.78) circle (2.5pt);
\begin{scriptsize}
\draw[color=black] (-4.2,3.3) node {$\lambda_0$};
\draw[color=black] (-7.8,3.3) node {$\lambda_\infty$};
\draw[color=black] (-3.4,4) node {$00u^i1$};
\draw[color=black] (-3.75,4.8) node {$1u^iu^{i+j}u^j$};
\draw[color=black] (-5.4,4.8) node {$100u^j$};
\draw[color=black] (-7.5,4.8) node {$01u^j0$};
\draw[color=black] (-7.45,5.25) node {$0100$};
\draw[color=black] (-5.4,5.65) node {$1000$};
\draw[color=black] (-3.5,5.65) node {$1u^i00$};
\draw[color=black] (-7.45,4) node {$0010$};
\draw[color=black] (-5.4,3.55) node {$0001$};
\draw[color=black] (-5.01,6.32) node {$\ell_0=L_1$};
\draw[color=black] (-6.97,6.32) node {$\ell_\infty=L_3$};
\draw[color=black] (-2.4,6.32) node {$V_i, 0\leq i<m$};
\draw[color=black] (-8.2,5.5) node {$T_2$};
\draw[color=black] (-1,4.6) node {$H_j, 0\leq j<m$};
\draw[color=black] (-8.2,3.8) node {$T_1$};
\end{scriptsize}
\end{tikzpicture}
\caption{The standard construction on $Q \colon xz-yw$ in $\PP^3_{\mathbb F}$, given $u\in {\mathbb F}$,
an $m$th primitive root of 1 for any $m>2$ not divisible by ${\rm char}(\mathbb F)$.
The lines $\lambda_0,\lambda_\infty, V_0,\ldots,V_{m-1}$ have two transversals, $T_1,T_2$.
Moreover, $G_{\mathcal L}$ is the multiplicative cyclic group of order $m$ for 
$\mathcal L=\{\lambda_0,V_0,\ldots,V_{m-1}\}$, $\mathcal L=\{\lambda_\infty, V_0,\ldots,V_{m-1}\}$ and,
when either $m$ is even or ${\rm char}(\mathbb F)=2$, also for $\mathcal L=\{\lambda_0,\lambda_\infty,V_0,\ldots,V_{m-1}\}$,
and in each case orbits are half grids (so geproci) with $m$ points on each line.}
\label{StConstFig2Trans}
\end{figure}

We now extend the standard construction to the case of a single transversal of multiplicity 2.
Pick a finite additive subgroup $A\subset \mathbb F$ of order at least 3.
Since $0<|A|<\infty$ we see ${\rm char}(\mathbb F)=p>0$, hence $|A|=p^s$ for some $s$, since additive subgroups
are $\mathbb Z/p\mathbb Z$ vector spaces.  
The lines $V_i$ defined by $ix-y,iw-z$, $i\in A$, and $H_j$ defined by $x-jw,y-jz$, $j\in A$,
are ruling lines on the quadric $Q \colon xz-yw$. Let $q_{ij}=(j:ij:i:1)$ denote the point $V_i\cap H_j$.
Also, let $\lambda_r$ be the line defined by $w$ and $x-rz$ for any $0\neq r\in\mathbb F$
(note $\lambda_r$ is tangent to $Q$ at $(0:1:0:0)$ but not contained in $Q$).
However, the matrix
$$\begin{pmatrix}
1 & 0 & 0 & 0\\
0 & 1 & 0 & 0\\
0 & 0 & r & 0\\
0 & 0 & 0 & r\\
\end{pmatrix}$$
induces an automorphism of $\PP^3_{\mathbb F}$
which is the identity on $T_1$ and $T_2$ and takes $\lambda_1$ to $\lambda_r$,
so up to projective equivalence we can assume $r=1$. We will do this and denote
$\lambda_1$ by $\lambda$.

Let $\mathcal L=\{\lambda, V_i, i\in A\}$ and let $Z$ be the set consisting of the $|A|^2$ points $q_{ij}$ together with the $|A|$ points
$(1:i:1:0)\in\lambda$, $i\in A$ (see Figure \ref{StConstFig1Trans}). 
An argument similar to that for the standard construction above shows $Z$ is an 
$(|A|,|A|+1)$-half grid (hence geproci) consisting of a single $C_{\mathcal L}$ orbit.

The points $q_{ij}$ comprise an $(|A|,|A|)$-grid whose grid lines are $V_i$ and $H_j$, $i,j\in A$.
The plane $V^*_i$ spanned by $V_i$ and $P$ is now defined by
$F_i=(c-id)(ix-y)-(ia-b)(z-iw)$, and 
the plane $H^*_j$ spanned by $H_j$ and $P$ is now defined by
$G_j=(b-jc)(x-jw)-(a-jd)(y-jz)$. Let $F=\prod_{i\in A}F_i$ and let 
$G=\prod_{j\in A}G_j$. Also let $\lambda*$ be the plane spanned and
by $\lambda$ and $P$, so $\lambda*$ is defined by $\Lambda=d(x-z)-(a-c)w$.
Then $M=F\Lambda$ defines the cone, with vertex $P$, over the lines in $\mathcal L$ and so has degree $|A|+1$.
And $F$ and $G$ define degree $|A|$ cones, both with vertex $P$, over the lines $V_i$ and $H_j$ respectively.
For each point $(1:i:1:0)\in\lambda$, $i\in A$, a unique member of the pencil defined by $F$ and $G$ vanishes at that point,
but in fact it's the same member for each $i\in A$, namely $F-G$, but $F-G$ does not vanish on all of $\lambda$.
We confirm these claims now.

To see that $F-G$ does not vanish on all of $\lambda$, plug in $(0:1:0:0)$:
$(F-G)(0,1,0,0)= (\Pi_i (id-c) -\Pi_j(a-jd))$ is clearly nonzero since
it is the difference of two nonzero forms in different variables.
Now we check that $F-G$ vanishes at the points $(1:i:1:0)\in\lambda$:
$$(F-G)(1,i,1,0)= \Pi_j((c-jd)(j-i)-(ja-b))-\Pi_j((b-jc)-(a-jd)(i-j))=$$
$$\Pi_j(-ja+b+(j-i)c-(j-i)jd)-\Pi_j((j-i)a+b-jc+j(i-j)d)=$$
$$\Pi_j((j-i)a+b-jc+j(i-j)d)-\Pi_j((j-i)a+b-jc+j(i-j)d)=0$$
where the second to last equality comes from substituting $i-j$ in for $j$ in 
$$\Pi_j(-ja+b+(j-i)c-(j-i)jd).$$
Thus $F-G,M$ meets the plane $H$ in a complete intersection of type $(|A|,|A|+1)$
so $Z$ is a $(|A|,|A|+1)$-geproci. It is not a grid since $\lambda$ is not contained 
in $Q$, so $Z$ is an $(|A|,|A|+1)$-half grid with respect to the lines $\mathcal L$.
The same argument as before (using the goupoid maps
$f_i\colon V_0\to\lambda$, $i\in A, i\neq0$) shows $Z$ is a single orbit.
By Corollary \ref{oneTransCor}, $G_{\mathcal L}$ is an additive subgroup of $\mathbb F$ of order
$|A|$, and thus $G_{\mathcal L}$ is isomorphic to $A$ (both being $\mathbb Z/p\mathbb Z$ vector spaces
with the same number of elements).

\begin{figure}[h]
\begin{tikzpicture}[x=1cm,y=1cm]
\clip(-11.09,3) rectangle (0,6.5);
\draw [line width=1pt] (-7.96,4.58)-- (-2,4.58);
\draw [line width=1pt] (-7.99,3.77)-- (-2,3.77);
\draw [fill=white,color=white] (-6.6,4.58) circle (5pt);
\draw [fill=white,color=white] (-6.23,3.77) circle (5pt);
\draw [line width=1pt, dotted] (-7.2,5.9) -- (-6,3.3);
\draw [line width=1pt] (-3,6)-- (-3,3);
\draw [line width=1pt] (-7,6)-- (-7,3);
\draw [line width=1pt] (-5,6)-- (-5,3);
\draw [line width=1pt] (-8,5.44)-- (-2,5.44);
\draw [fill=black] (-7,5.45) circle (2.5pt);
\draw [fill=black] (-7,4.6) circle (2.5pt);
\draw [fill=black] (-7,3.78) circle (2.5pt);
\draw [fill=black] (-5,5.45) circle (2.5pt);
\draw [fill=black] (-5,4.6) circle (2.5pt);
\draw [fill=black] (-5,3.78) circle (2.5pt);
\draw [fill=black] (-3,5.45) circle (2.5pt);
\draw [fill=black] (-3,4.6) circle (2.5pt);
\draw [fill=black] (-3,3.78) circle (2.5pt);
\begin{scriptsize}
\draw[color=black] (-5.9,3.2) node {$\lambda$};
\draw[color=black] (-3.4,4) node {$00i1$};
\draw[color=black] (-3.7,4.8) node {$(j\!:\!ij\!:\!i\!:\!1)$};
\draw[color=black] (-5.4,4.8) node {$j001$};
\draw[color=black] (-7.4,4.8) node {$0j10$};
\draw[color=black] (-7.4,5.25) node {$0100$};
\draw[color=black] (-5.4,5.65) node {$1000$};
\draw[color=black] (-3.5,5.65) node {$1i00$};
\draw[color=black] (-7.4,4) node {$0010$};
\draw[color=black] (-5.4,4) node {$0001$};
\draw[color=black] (-4.6,6.32) node {$V_0=L_1$};
\draw[color=black] (-6.97,6.32) node {$\ell_\infty=L_3$};
\draw[color=black] (-3.7,6.32) node {$\ldots$};
\draw[color=black] (-2.6,6.32) node {$V_i, i\in A$};
\draw[color=black] (-8.2,5.5) node {$T_2$};
\draw[color=black] (-1.3,4.6) node {$H_j, j\in A$};
\draw[color=black] (-1.8,4.3) node {$\vdots$};
\draw[color=black] (-1.4,3.8) node {$H_0=T_1$};
\end{scriptsize}
\end{tikzpicture}
\caption{A generalized standard construction on $Q \colon xz-yw$ in $\PP^3_{\mathbb F}$, given a finite additive
subgroup $0\neq A\subset {\mathbb F}$ (so ${\rm char}(\mathbb F)=p>0$).
The lines $V_i, i\in A$, are in one ruling on $Q$, while $T_1, T_2$ and $H_j, j\in A$, are in the other ruling and $\lambda$
is tangent to $Q$ at the point $(0:1:0:0)$ but not contained in $Q$. 
Here $G_{\mathcal L}=A$ for $\mathcal L=\{\lambda,V_i,i\in A\}$ and orbits are half grids (so geproci) with
$|A|$ points on each line.}
\label{StConstFig1Trans}
\end{figure}

\begin{theorem}\label{AbelianGrpClassifThm}
Consider an algebraically closed field $\mathbb F$ of characteristic $p\geq 0$.
Let $G$ be a finite abelian group. Then $G\cong G_{\mathcal L}$ for some set $\mathcal L$ 
of 3 or more skew lines in $\PP^3_{\mathbb F}$ if and only if 
$G$ is isomorphic to an additive subgroup of ${\mathbb F}$ or to a 
multiplicative subgroup of ${\mathbb F}^*$. In particular, if and only if
$G$ is a multiplicative cyclic group of order $m$ not a multiple of $p$
or $G=0$ or $p>0$ and $G$ is a finite dimensional $\mathbb Z/p\mathbb Z$ vector space.
\end{theorem}

\begin{proof}
By Theorem \ref{TransThm}, Proposition \ref{twoTrans} and Corollary \ref{oneTransCor},
if $G_{\mathcal L}$ is finite abelian, then it is an additive subgroup of ${\mathbb F}$ or a 
multiplicative subgroup of ${\mathbb F}^*$. Conversely, every finite subgroup of
${\mathbb F}$ or of ${\mathbb F}^*$ arises as $G_{\mathcal L}$ by the standard
constructions given above. The finite subgroups of ${\mathbb F}^*$ of order $m$ are the
solution sets of $x^m-1=0$, hence cyclic of order $m$ not divisible by $p$,
and finite nontrivial subgroups of ${\mathbb F}$ occur only when $p>0$,
in which case they are finite dimensional $\mathbb Z/p\mathbb Z$ vector spaces,
and dimension occurs since ${\mathbb F}$ is infinite dimensional over $\mathbb Z/p\mathbb Z$.
\end{proof}

\begin{remark}
Question 5.3 of \cite{Fields} asked if every complex $[m,s]$-half grid 
is projectively equivalent to one obtained by taking the half grid points on 
some choice of $s$ of the lines
of the $[m,r]$-half grid given by the standard construction 
(where $r$ is $m+1$ or $m+2$ if $m$ is even but just $m+1$ if $m$ is odd). 
Using the methods of this paper we can show the answer is yes for $[m,4]$-half grids. 
(For reasons of space we do not include a proof of that here, but we point out that this means
to classify $[m,4]$-half grids one merely needs to check which subsets of three of the lines $V_i$ together with $\lambda$ 
shown in Figure \ref{StConstFig1Trans} are projectively equivalent.)
However, in general the answer is no:
the $[6,8]$-half grid given in Example \ref{ClassicalExamples}(2)
is defined over $\mathbb C$ but no smooth quadric contains more than four of the half grid lines,
whereas in the standard construction there is a smooth quadric that contains
all but at most two of the lines. 
\end{remark}

\subsection{Do combinatorics determine geometry?}\label{sec:ProofOfThmB}
In this section we will prove Theorem \ref{CHGthm}, the main part of which
shows that finite sets which are collinearly complete on 3 or more skew lines
of $\PP^3$ over an algebraically closed field $\mathbb F$ of arbitrary
characteristic are geproci. Theorem \ref{CHGthm} is actually if and only if, 
but the converse part, which is much easier, is given by Proposition \ref{p. chg and finite union orbits} 
and Proposition \ref{Prop:[a,b]geprociImpliesCollCompl}.

The proof depends on applying properties of $h$-vectors of projections, so we first need some
preliminaries. 

\begin{lemma} \label{dgo}
    Let $X$ be a 0-dimensional complete intersection scheme in $\PP^2$ of type $(a,b)$ linking 
    a subscheme $X_1$ to a residual subscheme $X_2$. The $h$-vector of $X$ is related to the $h$-vectors of $X_1$ and $X_2$ by the following formula.
    \[
    h_X(t) - h_{X_1}(t) = h_{X_2}(a+b-2-t).
    \]
\end{lemma}

\begin{proof}
    See \cite{DGO} Theorem 3 or \cite{JuanBook} Corollary 5.2.19. Notice that $a+b-2$ is the degree of the last non-zero entry of $h_X$.
\end{proof}

For a homogeneous ideal $I \subset k[x,y,z]$ and a positive integer $m$ we denote by $I_{\leq m}$ the ideal generated by the homogeneous components of $I$ of degree $\leq m$. In particular, the scheme defined by $I_{\leq m}$ is the (scheme-theoretic) base locus of the component $[I]_m$. The following result is presented only in the form needed in this paper. See the papers cited in the proof for (much) more general versions. 

\begin{lemma} \label{davis}
    Let $X \subset \PP^2$ be a reduced set of points. Assume for some $s\leq d$ that $h_X(d) = h_X(d+1) = s$. Then
    the base locus of $(I_X)_{\leq d}$ contains a reduced curve $C$ 
    defined by a homogeneous polynomial $f$ of degree $s$, and $f$ is a factor of every element in $[I_X]_{\leq d}$.
    Moreover, let $X_1$ be the set of points of $X$ lying on $C$ and $X_2$ the set of points of $X$ not lying on $C$.
    Then:

    \begin{itemize}
    \item[(i)] $(I_X)_{\leq d}$ is the saturated ideal of $C \cup X_2$;

        \item[(ii)] $(I_{X_1})_{\leq d} = I_C = (f)$;

        \item[(iii)] for any $t$, 
        \[
        h_{X_1}(t) = \left \{
        \begin{array}{ll}
        \Delta H(C,t)  & \hbox{for $t \leq d+1$}, \\
         h_X (t) & {\hbox{for $t \geq d$}}
         \end{array}
        \right.
        \]
        (in particular, $h_{X_1}(t) = s$ for $s \leq t \leq d+1$); and
        \item[(iv)] for $s \leq t \leq d$ we have $h_{X_2}(t-s) = h_X(t) - s$.
    \end{itemize}
\end{lemma}

\begin{proof}
Most of this can already be found in \cite{Davis}, but we cite \cite{BGM} since \cite{BGM} is more accessible
and its exposition aligns better with our situation. Indeed, \cite[Theorem~3.6]{BGM}
gives everything but part (iv). From (i) we see that $f$ is a GCD for $(I_X)_{\leq d}$,
so we can apply \cite[Theorem 2.4 (b)]{BGM} with $k=s$ and $r=1$, which gives (v).
\end{proof}

\begin{example}
Let $X' \subset \PP^2$ be a reduced set of points with $h$-vector $(1,2,3,4,5,5,4,3)$. (This can be produced as the residual of a set of three non-collinear points inside a complete intersection of type $(5,6)$.) Let $X''$ be a complete intersection of type $(6,7)$ containing $X'$. Then $X''$ links $X'$ to a residual set of points $X$. Since $X''$ has $h$-vector $(1,2,3,4,5,6,6,5,4,3,2,1)$, the $h$-vector of $X$ is then computed using Lemma \ref{dgo} to be $(1,2,3,4,2,2,1)$. By Lemma \ref{davis}, we see that $X=X_1\cup X_2$ is a disjoint union such that $X_1$ lies on a conic and has $h$-vector $(1,2,2,2,2,2,1)$. Thus
$|X_1|=1+2+2+2+2+2+1=12$ and $X_2$ consists of 3 points off the conic.
(This could also be deduced from simple geometric considerations, but it illustrates the utility of Lemmas \ref{dgo} and \ref{davis}.)
\end{example}

\begin{theorem}\label{lemmaforThmB}
    Let $X=X_1\cup X_2 \subseteq \mathbb P^2$ be a set of $a(b+1)$ reduced points. 
    Assume
    \begin{enumerate}
     \item $b\ge 3$;
        \item $X_1$ is a complete intersection of type $(a,b)$ consisting of $a$ points 
        on each of $b$ distinct lines $\ell_1, \ldots, \ell_b$;
        \item $X_2$ is a set of a collinear points on a line $\ell'$ where $\ell'$ is distinct from the lines $\ell_i$;
        \item  $X$ is contained in the smooth locus of $\ell'\cup \ell_1\cup \cdots\cup\ell_b$ 
        and exactly two of the lines meet at any point;
        \item $X\setminus \ell_i$ is a complete intersection of type $(a,b)$ for any $\ell_i$.
    \end{enumerate}
 Then $X$ is a complete intersection of type $(a,b+1)$.   
\end{theorem}

\begin{proof}
Let $\gamma_i$ be the curve defined by the general element in the linear system $[I_{X_i}]_a$, for $i=1,2$. 
Note that $\ell'$ is not among the lines $\ell_i$ by (3), and $\gamma_2\cap\ell'=X_2\subsetneq \ell'$ by generality of $\gamma_2$, so
$\ell'$ is not a component of $\ell_1\cup\cdots\cup \ell_b\cup \gamma_2$.
Also, $X_1$ is the complete intersection of $\ell_1\cup\cdots\cup \ell_b$ and $\gamma_1$ so
$X_1=(\ell_1\cup\cdots\cup \ell_b)\cap \gamma_1$, hence
$\gamma_1$ and $\ell_1\cup\cdots\cup \ell_b$ have no common components
and every component of $\gamma_1$ contains at least one point of $X_1$.
We also see that $\gamma_1$ does not contain $\ell'$ since otherwise 
$X_1$ would contain singular points of $\ell'\cup \ell_1\cup \cdots\cup\ell_b$, contradicting (4).
Moreover, $[I_{X_2}]_a$ contains elements that do not vanish at any point of $X_1$
(for example, take the curve whose components consist of a general line through each point of $X_2$), hence
we can assume $\gamma_2$ does not contain any point of $X_1$ by generality of $\gamma_2$.
Thus $\gamma_1$ and $\gamma_2$ have no common components and 
$\ell_1\cup \cdots\cup\ell_b$ and $\gamma_2$ have no common components.
In particular, the curves $C_1=\ell'\cup \gamma_1$ and $C_2=\ell_1\cup\cdots\cup\ell_b\cup \gamma_2$ 
are reduced and have no common components. 
Thus $V=C_1\cap C_2$ is a complete intersection of type $(a+1,a+b)$ containing $X$. 
Note that as sets $V=B_1\cup B_2 \cup X_1\cup X_2$, where
$B_1=\ell'\cap(\ell_1\cup \cdots\cup\ell_b)$, $B_2=\gamma_1\cap\gamma_2$
(recall $X_1=\gamma_1\cap \ell_1\cup\cdots\cup\ell_b$ and $X_2=\ell'\cap\gamma_2$).

It is not necessarily true that $V$ is reduced. 
In particular, $\deg V=\deg(C_1)\deg(C_2)=(a+1)(a+b)=b+a^2+ab+a=\deg(B_1)+\deg(B_2)+\deg(X_1)+\deg(X_2)$;
thus $V$ being reduced is equivalent to $B_1,B_2,X_1,X_2$ each being reduced and being pairwise disjoint.
Our assumptions imply $B_1,X_1,X_2$ are reduced and pairwise disjoint. Also,
$B_1\cap B_2=(B_1\cap\ell')\cap B_2=B_1\cap(\ell'\cap B_2)\subseteq B_1\cap(\ell'\cap\gamma_2)=B_1\cap X_2
=\varnothing$ and $B_2\cap X_1\subseteq \gamma_2\cap X_1=\varnothing$ (since a general line through 
a point of $X_2$ misses every point of $X_1$), but we cannot assume $B_2$ is reduced or that $B_2$ is 
disjoint from $X_2$. (For example, if $a<b$, then there is a unique curve
$\gamma_1$ which, if the theorem we're proving is true, must contain $X_2$ and thus $X_2\subset B_2$, so $V$
would not be reduced at any point of $X_2$. It is also possible a priori that $\gamma_1$ is singular
at a point of $X_2$, or, if the characteristic is positive, that
some point of $X_2$ is on every tangent line of some component of $\gamma_1$, in which case
$\gamma_1$ and $\gamma_2$ do not meet transversely and hence $B_2$ would not be reduced.)

In any case, let $Y$ be the residual of $X$ in $V$ (so as a set $Y=B_1\cup B_2$, but
$Y$ need not be reduced and $X$ and $Y$ need not be disjoint).
We have the following relevant facts:

\begin{enumerate}[label=$(\roman*)$]
\item $Y\cap \ell_i$ is the single reduced point $q_i\in B_1$, where $\ell'\cap\ell_i=\{q_i\}$,
and $V\cap \ell_i$ is the reduced set $(X_1\cap\ell_i)\cup\{q_i\}$.
This is because $Y\cap \ell_i=(B_1\cap\ell_i)\cup(B_2\cap\ell_i)\subseteq (B_1\cap\ell_i)\cup(B_2\cap X_1)
=(B_1\cap\ell_i)\cup\varnothing=\ell'\cap\ell_i\subseteq Y\cap\ell_i$.
Similarly, $V\cap \ell_i=(B_1\cap\ell_i)\cup(X_1\cap\ell_i)=\{q_i\}\cup(X_1\cap\ell_i)$
is reduced.

\item $Y$ is contained in a unique curve of degree $a$ so $h_Y(a)=a$. To see this, by (i) note
$Y\setminus\ell_i=Y\setminus(\ell_i\cap Y)=Y\setminus\{q_i\}$ is linked by the complete intersection 
$V \setminus \ell_i=V \setminus (\ell_i\cap V)=V\setminus((X_1\cap\ell_i)\cup\{q_i\})$ to $X \setminus \ell_i=X \setminus(\ell_i\cap X_1)$,  
and the latter is a complete intersection of type 
$(a,b)$ so, by Lemma \ref{dgo}, for each $i=1,\ldots, b$ we have
the following relation among the $h$-vectors of $X\setminus \ell_i,V\setminus \ell_i$ 
    and $Y\setminus \ell_i$:  
$$h_{Y\setminus \ell_i}(a)=h_{V\setminus \ell_i}(a)-h_{X\setminus \ell_i}(2a+b-2-a)=$$
$$h_{V\setminus \ell_i}(a)-h_{X\setminus \ell_i}(a+b-2)=a+1-1=a.$$
Thus $Y\setminus \ell_i$ lies on a unique curve of degree $a.$ 
Letting $i$ vary, any two such such curves intersect in a zero-dimensional scheme of degree at least $a^2+b-2$ 
which includes $b-2$ (which is positive since $b\geq3$) distinct points of $B_1\subset\ell'$. 
Since $B_1\cap B_2=\varnothing$, any two such curves share at least one point 
outside the base locus $B_2$ of the pencil defined by $\gamma_1$ and $\gamma_2$. 
Thus, the two curves must coincide, giving the unique curve of degree $a$ containing $Y$. 
\item Applying $(ii)$ gives $h_X(a+b-1)=h_V(2a+b-1-(a+b-1))-h_Y(2a+b-1-(a+b-1))$
$=h_V(a)-h_Y(a)=a+1-a=1$. 
\item $h_X(j)=j+1$ for $j<\min\{a,b+1\}$: this follows since $I_X\subseteq I_{X_1}$ has no minimal generators of such degrees. 
\item If $0<h_X(i)\le i$ and $h_X(i)<\min\{a,b+1\}$, then $h_X(i)>h_X(i+1)$ 
(namely, the $h$-vector of $X$ is of so-called ``decreasing type"). To see this, by way of contradiction let $d$ be the least integer
for which for some $i$ we have $0<h_X(i)\le i$ but $d=h_X(i)\leq h_X(i+1)$. 
By Macaulay’s Theorem for O-sequences of codimension 2 we have $h_X(i+1)\leq h_X(i)^{<i>}$. Since 
$$h_X(i)=d=\binom{i}{i}+\binom{i-1}{i-1}+\cdots +\binom{i-d+1}{i-d+1}$$
we get
$$h_X(i)\leq h_X(i+1)\leq h_X(i)^{<i>}=\binom{i+1}{i+1}+\binom{i}{i}+\cdots +\binom{i-d+2}{i-d+2}=h_X(i)$$
so we actually must have $h_X(i) = h_X(i+1)=d$. By the minimality of $d$, 
there cannot be another flat in $h_X$ after this one ends. But by $(iii)$ we know $h_X$ has at least 
$a+b$ non-zero entries and from $d=h_X(i)<\min\{a,b+1\}$ we see $a+b>2d$.  
Thus, the following  is a lower bound for the $h$-vector of the points of $X$ on $C$:
		 $$(1, 2,3,  \ldots,d-1 , d,d, \ldots,d, d-1 \ldots, 3,2,1).$$
By Lemma \ref{davis}, $C$ contains a subset of $X$ of (at least)
		 $$2(1+2+ \cdots + d)+ d(a+b  - 2d)=2\binom{d+1}{2}+ d(a+b-2d)=d(a+b-d+1)$$
points. Now observe the curve $C$ has degree $d$ and is of the form either 
$C=\ell_{i_1}\cup\cdots\cup \ell_{i_j}\cup C'$ or $C=\ell'\cup\ell_{i_2}\cup\cdots\cup \ell_{i_j}\cup C'$
where $C'$ does not contain any of the lines $\ell_1,\ldots, \ell_b,\ell'$ as a component.
Thus $C$ contains at most $N=ja+(b+1-j)(d-j)$ distinct points of $X$ (i.e., $a$ on each of the $j$ lines 
and  $(b+1-j)(d-j)$ due to the intersection of $C'$ with the remaining lines). 
This gives a contradiction since
$d(a+b-d+1)=j(a+b+1-d)+(d-j)(b+1+a-d)>ja+(d-j)(b+1)\ge N$.
\end{enumerate}

To prove $X$ is a complete intersection we need to show the existence of a curve containing $X$ of degree $a$ but not containing
any of the lines $\ell_1,\ldots, \ell_b, \ell'$ as a component (since this curve of degree $a$ together with
the curve $\ell_1\cup\cdots\cup\ell_b, \ell'$ defines a complete intersection of type $(a,b+1)$ containing $X$
which is itself a reduced set of $a(b+1)$ points and thus equals the complete intersection).
We consider the following cases:
\begin{itemize}
 \item $a\ge b+2$: In this case, the curve $\ell'\cup\ell_1\cup\cdots\cup\ell_b$ 
 is contained in every curve containing $X$ of degree less than $a$
 (since $\ell'\cup\ell_1\cup\cdots\cup\ell_b$ does contain $X$ and each line contains $a>b+1$ of the points so 
 by B\'ezout's Theorem must be a component of every curve of degree at most $a-1$ containing $X$). Thus
 $$h_X(i)=\begin{cases}i+1 & \text{for}\ 0\le i<b+1 \text{ by (iv)};\\
b+1 & \text{for}\ b+1\leq i\leq a-1 \text{ by B\'ezout's Theorem};\\
1 & \text{for}\ i=a+b-1 \text{ by (iii); and}\\
0 & \text{for}\ i>a+b-1 \text{ by (v)}\\
\end{cases}$$ 
so
   \[
   ab +a =  |X| = h_X(0) + h_X(1) + \cdots + h_X(a+b-1)
    \]
and
$$\left(\sum_{i=0}^{b+1} h_X(i)+ h_X(a+b-1)\right)+\sum_{i=b+2}^{a-1} h_X(i)= \binom{b+3}{2}+(b+1)(a-b-2).$$

Now consider degrees $a\leq i\leq a+b-2$. From $(v)$, if $h_X(i)<b+1$ then $h_X(i)>h_X(i+1)$.
Therefore, we get 
$$(*)\,\,\, h_X(a+b-2)\ge 2,\ h_X(a+b-3)\ge 3,\ \ldots,\ h_X(a)\ge b.$$ 
Thus, 
$$\sum_{i=a}^{a+b-2} h_X(i)\ge \binom{b+1}{2}-1.$$
Hence, collecting all the information about the $h_X$, we have
  $$\begin{array}{rl}
      |X| \ge & \displaystyle \binom{b+3}{2}+ (b+1)(a-b-2)+ \binom{b+1}{2}-1\\[15pt]
       =&\displaystyle \dfrac{b^2+5b+6+b^2+b- 2b^2-6b-4}{2}-1+(b+1)a \\
       =&(b+1)a\\
       =&|X|.
  \end{array}$$
 Thus the inequalities in $(*)$ are equalities; in particular, $h_X(a)=b=h_X(a-1)-1$. Hence $I_X$ has a minimal generator of degree $a$. 
By B\'ezout, this generator defines a curve having none of the lines $\ell_1, \ldots, \ell_b,\ell'$ as a component
(since removing a line which is component gives a curve of degree $a-1$ containing $a$ points on each of the other lines).

 \item $a=b+1$: In this case we have
 $$h_X(i)= \begin{cases}
     i+1& \text{if }\ 0\le i\le b \text{ by (iv);}\\
     1 & \text{if }\ i=2b \text{ by (iii)};\\
     0 & \text{for}\ i>2b \text{ by (v)}.\\
 \end{cases}$$ 
Now consider degrees $b+1\leq i\leq 2b-1$. From $(v)$, 
if $h_X(i)<b+1$ then $h_X(i)> h_X(i+1)$.
Therefore, we get   $$h_X(2b-1)\ge 2,\ h_X(2b-2)\ge 3,\ \ldots,\ h_X(b+1)\ge b.$$ 
  Thus,
  $$\begin{array}{rl}
      |X|=&\sum h_X(i)=\displaystyle \sum_{i=0}^{b} h_X(i)+  \sum_{i=b+1}^{2b} h_X(i) \\
       \ge & \displaystyle \binom{b+2}{2}+ \binom{b+1}{2}\\[15pt]
       =&\displaystyle \dfrac{b^2+3b+2+b^2+b}{2} = b^2+2b+1\\
       =&(b+1)^2=|X|.
  \end{array}$$
 Then, all the inequalities above are equalities. In particular  $h_X(b+1)= b$. 
 Therefore, $I_X$ has two minimal generators of degree $b+1$ and we
 finish by arguing as before.
 
 \item $a\le b$: In this case we have 
 $$h_X(i)= \begin{cases}
     i+1& \text{if }\ 0\le i\le a-1 \text{ by (iv);}\\
      1 & \text{if }\ i=a+b-1 \text{ by (iii); and}.\\
      0 & \text{for}\ i>a+b-1 \text{ by (v)}.\\
 \end{cases}$$

Now consider degrees $a\leq i\leq a+b-2$. From $(v)$, 
if $h_X(i)<a+1$ then $h_X(i)> h_X(i+1).$
Therefore, we get  $$h_X(a+b-2)\ge 2,\ h_X(a+b-3)\ge 3,\ \ldots,\ h_X(b+1)\ge a-1, h_X(b)\ge a, $$
$$h_X(b-1)\ge a, \ldots, h_X(a)\ge a.$$ 
Thus, 
  $$\begin{array}{rl}
      |X|=&\sum h_X(i)=\displaystyle \sum_{i=0}^{a-1} h_X(i) +  \sum_{i=a}^{b} h_X(i)+\sum_{i=b+1}^{a+b-1} h_X(i)  \\
       \ge & \displaystyle \binom{a+1}{2}+ a(b-a+1)+ \binom{a}{2}\\[15pt]
       =&\displaystyle \dfrac{a^2+a+a^2-a}{2}+a(b+1)-a^2= \\
       =&a(b+1)\\
       =&|X|.
  \end{array}$$
Again, all the inequalities above are indeed equalities so $h_X(a)= a$ hence $I_X$ has a minimal generator of degree $a$
and we finish as before.
\qedhere
\end{itemize}
\end{proof}

We can now prove one of our main results, which shows there is a tight connection between
the combinatorics of skew lines and the algebraic geometric notion of geproci sets:

\begin{proof}[Proof of Theorem \ref{CHGthm}]\label{proof_CHGthm}

If either $Z$ is an $\{a,b\}$-grid or $Z$ is an $[a,b]$-half grid with respect to $\mathcal L$,
then $Z$ is $[a,b]$-geproci, hence collinearly complete with respect to $\mathcal L$ by Proposition \ref{Prop:[a,b]geprociImpliesCollCompl}(a),
which by Proposition \ref{p. chg and finite union orbits} means $Z$ is a finite union of finite $C_{\mathcal L}$ orbits.

For the converse, assume $Z$ is collinearly complete for $\mathcal L$. Enumerate the lines as $\mathcal L=\{L_1,\ldots,L_b\}$.
Let $\mathcal L_i=\{L_1,\ldots,L_i\}$, let $Z_i=Z\cap(L_1\cup\cdots\cup L_i)$ and let $\overline{Z_i}$ be the general projection of
$Z_i$ to a plane. When $i=3$, $Z_i$ is a grid, hence $[a,3]$-geproci, by Proposition \ref{Prop:[a,b]geprociImpliesCollCompl}(a).
This proves the result when $b=3$. 

Now assume $b>3$ and apply induction. Assume $3\leq i<b$; by induction we may assume
that the set of $ai$ points $Z'$ of $Z$ on any choice of $i$ of the lines in $\mathcal L$ is $[a,i]$-geproci.
In particular, $X_1=\overline{Z_i}$ is a complete intersection of type $(a,i)$.
Let $X_2\subset\ell'$ be the projected image of the points of $Z$ on $L_{i+1}$, where the projection of
$L_{i+1}$ is $\ell'$. The fact that the projection is general ensures that 
the hypotheses of Theorem \ref{lemmaforThmB} are satisfied, so we conclude $\overline{Z_{i+1}}=X=X_1\cup X_2$ is a 
complete intersection of type $(a,i+1)$. 
It follows by induction that $Z$ is $[a,b]$-geproci and thus either an $\{a,b\}$-grid or $Z$ is an $[a,b]$-half grid.
\end{proof}

\subsection{Classifying single orbit $[m,4]$-half grids}\label{Sec:ClassifyingCombinatorial}

First, by analogy, consider the problem of classifying $[m,3]$-grids $Z$ for some $m>3$. 
Since $m>3$, the $3m$ points uniquely determine the 3 skew lines, each of which contains 
$m$ of the points. There is only one projective equivalence class of 3 skew lines
so to classify $[m,3]$-grids $Z$ we can fix any convenient choice of the 3 skew lines.
Each $[m,3]$-grid $Z$ on those 3 lines is a union of $m$ orbits, where each orbit is 
the intersection of the 3 lines with a transversal. The 3 lines determine a unique smooth quadric $Q$.
The 3 lines all come from one ruling on $Q$, all transversals come from the other ruling.
Thus there is a bijection between $[m,3]$-grids $Z$ and a choice of $m$ points on
any one of the 3 lines. I.e., the projective equivalence classes of $[m,3]$-grids
is the same as the projective equivalence classes of $m$ points on $\PP^1$
(which by the way has positive dimensional moduli; cf. Remark \ref{PosDimModuli}).

In the same way, classifying $[m,4]$-half grids for $m\geq4$ amounts to classifying
sets of 4 skew lines, classifying single orbit geproci sets on those lines
and then classifying unions of single orbits having a total of $m$ points per line.
Unions of single orbits introduce positive dimensional moduli, just as in the
case of $[m,3]$-grids, by the relative disposition of the orbits (see Remark \ref{PosDimModuli}). 
But unlike the case of $[m,3]$-grids where single orbits on the 3 lines are all projectively 
equivalent (being 3 collinear points), single orbit $[m,4]$-half grids are not all projectively equivalent.
However, there are only finitely many projective equivalence classes of single orbit $[m,4]$-half grids.
So what we propose to do here is count the number of classes.

We focus on the case of there being two distinct transversals. 
Our approach 
in this case
is as follows. Finite nontrivial orbits on four skew lines arise exactly 
when there are two distinct transversals and the lines are not all on a single quadric.
Up to projective equivalence we can assume three of the lines and the two transversals
are in standard position. By Remark \ref{AlgorithmRem}, there are only finitely many choices of the fourth line so that
single orbits off the transversals have $m$ points on each line. We begin by finding the number
of these choices for the fourth line (in a characteristic free way)
in Proposition \ref{LineCount}. 

To explain this in more detail, fix $m>4$, and consider 
distinct skew lines $L_1,L_2,L_3$ and two distinct transversals $T_1,T_2$.
We would first like to count how many choices we have for a line $L_4$
meeting both $T_1$ and $T_2$ such that $G_1$ has order $m$.
By Remark \ref{AlgorithmRem},
counting the number of lines $L_4$ which have $|G_1|=m$
amounts to counting the pairs $(\alpha,\beta)\in{\field}^2$ such that
$\alpha\beta\neq1$, $\alpha\neq1\neq\beta$, $0\neq\alpha\beta$,
where $\alpha$ and $\beta$ generate a multiplicative group of order $m$.
Thus $\alpha^m=\beta^m=1$ so $\alpha$ and $\beta$
are $m$th roots of 1 which generate the full group of 
$m$th roots of 1.

One approach is to work in the ring $\mathbb Z_m=\mathbb Z/m\mathbb Z$
and find the number $n_m$ of pairs of nonzero elements $(i,j)$ with $i+j\neq0$ such that
$\langle i,j\rangle=\mathbb Z_m$, where $\langle i,j\rangle$ is the additive subgroup generated by $i$ and $j$. 
(The condition for $\langle i,j\rangle =\mathbb Z_m$ is ${\rm gcd}({\rm gcd}(i,m),{\rm gcd}(j,m))={\rm gcd}(i,j,m)=1$.)
In the next result we determine $n_m$ in terms of Euler's totient function $\phi(m)$
and the prime factorization $m=p_1^{e_1}\cdots p_r^{e_r}$ of $m$.

\begin{proposition}\label{LineCount}
Let $p={\rm char}(\overline\field)$ and let $m=p_1^{e_1}\cdots p_r^{e_r}$ be the prime factorization for an integer $m>1$
(so the $p_k$ are distinct primes and $e_k\geq1$ for each $k$).
Consider three distinct skew lines $L_1,L_2,L_3$ and two distinct transversals $T_1,T_2$.
If $p>0$ and divides $m$, let $n_m=0$, but if $p=0$ or if $p>0$ and does not divide $m$, let
$$n_m=\phi(m)(\phi(m)-1)+2\phi(m)(m-1-\phi(m))+
\phi(m)\sum_S\big(\big(\prod_{k\in S}p_k^{e_k-1}\big)\big(\prod_{k\not\in S}p_k^{e_k}-\prod_{k\not\in S}\phi(p_k^{e_k})\big)\big).$$
Then there are exactly $n_m$ lines $L_4$ meeting both $T_1$ and $T_2$ such that $G_1$ has order $m$
(where the sum is over all nonempty proper subsets $S\subsetneq\{1,\dots,r\}$
and thus has $2^r-2$ terms).
In particular, if $m=p$ is prime, then the displayed equation gives  $n_m=(m-1)(m-2)$, and if
$m=p^e$ for a prime $p$ with $e\geq1$, then it gives 
$n_m=(m-\frac{m}{p})(m-\frac{m}{p}-1)+2(m-\frac{m}{p})(\frac{m}{p}-1)$.
\end{proposition}

\begin{proof}
If $p>0$ and $p|m$, then $\overline\field^*$ has no subgroup of order $m$, so $n_m=0$.
So assume $p=0$ or that $p>0$ but does not divide $m$; in this case
$\overline\field^*$ has a unique subgroup of order $m$, in which case
the expression for $n_m$ has three main terms.
The first term counts the number of pairs $(i,j)$
with $i,j\in\mathbb Z_m$ and $i+j\neq0$ where $\langle i \rangle=\langle j \rangle=\mathbb Z_m$.
The second term counts the number of pairs $(i,j)$
with $i,j\in\mathbb Z_m$ where either $\langle i\rangle=\mathbb Z_m$ or $\langle j\rangle=\mathbb Z_m$
but not both.
The third term counts the number of pairs $(i,j)$ where 
$\langle i,j \rangle=\mathbb Z_m$ but $\langle i \rangle\subsetneq\mathbb Z_m$ and $\langle j \rangle\subsetneq\mathbb Z_m$.

The number of pairs $(i,j)$
with $i,j\in\mathbb Z_m$ and $i+j\neq0$ where $\langle i \rangle=\langle j \rangle=\mathbb Z_m$
is $\phi(m)(\phi(m)-1)$, since the number of elements $i\in \mathbb Z_m$
of (additive) order $m$ is $\phi(m)$ and for each of them there are
$\phi(m)-1$ elements $j\in \mathbb Z_m$ of order $m$ where $i+j\neq0$.

The number of pairs $(i,j)$
with $i,j\in\mathbb Z_m$ where $\langle i \rangle=\mathbb Z_m$ but $\langle j \rangle\neq\mathbb Z_m$
is $\phi(m)(m-1-\phi(m))$, since there are $\phi(m)$ elements $i$
of order $m$, and for each such $i$ there are $m-1-\phi(m)$ nonzero elements
$j$ not of order $m$. Swapping the roles of $i$ and $j$ doubles the count.

To count the pairs $(i,j)$ where 
$\langle i,j \rangle=\mathbb Z_m$ but $\langle i \rangle\subsetneq\mathbb Z_m$ and $\langle j \rangle\subsetneq\mathbb Z_m$
we use the ring isomorphism
$\mathbb Z_m\cong \oplus_k \mathbb Z_{p_k^{e_k}}$.
Note that elements $i=(i_1,\ldots,i_r), j=(j_1,\ldots,j_r)\in \oplus_k\mathbb Z_{p_k^{e_k}}$ 
together generate the full group under addition if and only if
either $i_k$ or $j_k$ generates $\mathbb Z_{p_k^{e_k}}$ for each $k$.
The counting strategy we will use here is as follows:
for each nonempty proper subset $S\subsetneq\{1,\dots,r\}$,
count the number of elements $i=(i_1,\ldots,i_r)$ such that
$\langle i_k\rangle \subsetneq \mathbb Z_{p_k^{e_k}}$ exactly when $k\in S$,
and for each such $i$, count the number of elements $j=(j_1,\ldots,j_r)$
such that $\langle j_k\rangle\neq \mathbb Z_{p_k^{e_k}}$ 
exactly when $k\not\in S$).

There are $p_k^{e_k}-\phi(p_k^{e_k})=p_k^{e_k-1}$ elements $i_k$ such that
$i_k$ has order less than $p_k^{e_k}$ and $\phi(p_k^{e_k})$ 
elements $i_k$ such that $i_k$ has order $p_k^{e_k}$.
Hence there are $\prod_{k\in S}p_k^{e_k-1}\prod_{k\not\in S}\phi(p_k^{e_k})$
elements $i$, and similarly there are 
$\prod_{k\in S}\phi(p_k^{e_k})\big(\prod_{k\not\in S}p_k^{e_k}-\prod_{k\not\in S}\phi(p_k^{e_k})\big)$ elements $j$, so
the third term is
$$\sum_S\left(\prod_{k\in S}p_k^{e_k-1}\prod_{k\not\in S}\phi(p_k^{e_k})\prod_{k\in S}\phi(p_k^{e_k})\big(\prod_{k\not\in S}p_k^{e_k}-\prod_{k\not\in S}\phi(p_k^{e_k})\big)\right)$$
or
$$\phi(m)\sum_S\left(\big(\prod_{k\in S}p_k^{e_k-1}\big)
\big(\prod_{k\not\in S}p_k^{e_k}-\prod_{k\not\in S}\phi(p_k^{e_k})\big)\right).\eqno\qedhere$$
\end{proof}

For $b$ skew lines $\mathcal L$ with two transversals, 
it's possible for $G_{\mathcal L}$ to be large when $b$ is small (but bigger than 3),
but we can apply the preceding result to see that $G_{\mathcal L}$ cannot be too small
if $b$ is large (as long as $|G_{\mathcal L}|\neq1$). Examples of this behavior were seen in \cite{POLITUS},
which showed (over the complex numbers) that if $b>4$ and $|G_{\mathcal L}|>1$ then $|G_{\mathcal L}|>3$,
and \cite{Fields}, which (over the complex numbers) showed if $b>6$ and $|G_{\mathcal L}|>1$ then $|G_{\mathcal L}|>4$.

\begin{corollary}\label{LineCountCor}
For any characteristic, assume $\mathcal L$ is a set of $b\geq3$ skew lines
in $\PP^3_{\field}$ with two distinct transversals such that $|G_{\mathcal L}|> 1$.
\begin{enumerate}
\item[(a)] If $b>\frac{2r(r-1)(r-2)}{3}+2$, then $|G_{\mathcal L}|>r$. 
\item[(b)] If $b>2(\sum_{m=1}^r n_m)+2$, then $|G_{\mathcal L}|>r$.
\end{enumerate}
\end{corollary}

\begin{proof}
Let  $\mathcal L=\{L_1,\ldots,L_b\}$.
The lines $L_1,L_2,L_3$ define a quadric $Q$. Let $S$ be the subset of the lines $L_i$
not contained in $Q$ and let $s=|S|$. 
Let $S'$ be the set of the remaining $s'=b-s-2$ lines, other than $L_1,L_2$,
thus every line in $S'$ is contained in $Q$.

Recall $G_{\mathcal L}$ is a subgroup of the multiplicative group $\field^*$.
If $|G_{\mathcal L}|\leq r$, then $1<|G_{{\mathcal L}_L}| = m\leq r$ for 
each four line set $\mathcal L_L=\{L_1,L_2,L_3,L\}$, $L\in S$.

(a) Applying Remark \ref{AlgorithmRem}, there are at most $(m-1)(m-2)$ choices for 
$l$ and $t$ defining $L$ such that $G_{{\mathcal L}_L}$ is nontrivial but contained in
a cyclic group of order $m$, so $s\leq\sum_{i=1}^n (m-1)(m-2)=r(r-1)(r-2)/3$. 
Let $S'$ be the set of the remaining $s'=b-s-2$ lines, other than $L_1,L_2$,
thus every line in $S'$ is contained in $Q$.
Pick any line $L\in S$ and let $Q'$ be the quadric defined by $L_1,L_2,L$. Then $Q\cap Q'$ consists
of the two transversals and $L_1\cup L_2$. The same argument as before shows 
$b-s-2=s'\leq r(r-1)(r-2)/3$. Thus $b-2\leq 2r(r-1)(r-2)/3$.

(b) There are $n_m$ lines $L\in S$ such that $|G_{{\mathcal L}_L}| = m$.
Thus there are $\sum_{m=1}^r n_m$ lines $L\in S$ such that $|G_{{\mathcal L}_L}| \leq m$.
Therefore, arguing as before, if $b>2(\sum_{m=1}^r n_m)+2$, then $|G_{\mathcal L}|>r$.
\end{proof}

\begin{remark}
While the bound in Corollary \ref{LineCountCor}(a) is a lot easier to compute than the bound in (b),
the one in (b) is often smaller, but still not sharp (note that the bound in (b) depends on
the characteristic). Whereas over the complex numbers $b$ is at most 4 when $|G_{\mathcal L}|=r=3$ \cite{POLITUS},
and $b$ is at most 6 when $|G_{\mathcal L}|=r=4$ \cite{Fields}, both bounds in Corollary \ref{LineCountCor} are 6 and 18 respectively.
However, in characteristic $p>0$ the Hopf fibration gives examples with $|G_{\mathcal L}|=r=q+1$ and $b=q^2+1$
over a field of order $q$. In characteristic $p=2$ the bound (b) for $r=3$ is 6 compared to the Hopf fibration example of 5 for $q=2$, 
and in characteristic $p=3$ the bound (b) for $r=4$ is 14 compared to the Hopf fibration example of 10 for $q=3$.
\qed\end{remark}

The bounds in Corollary \ref{LineCountCor} come from noting that
the lines $L_4,\ldots, L_b$ not on the quadric defined by $L_1,L_2,L_3$
must all come from a set $S$ (determined by Remark \ref{AlgorithmRem})
of known cardinality. But it ignores the possibility that some lines in $S$ 
might not ever occur together when the group has given order.
For example, the lines in $S$ might not be pairwise skew. 
If they are not necessarily pairwise skew, then choosing some of the lines excludes choosing some of the others, 
so taking that into account might give better bounds.
However, the next result shows the lines one can choose for $L_4$ to get a group of prime order are skew.

\begin{proposition}\label{The t and l values are different}
Let ${\rm char}(\overline{\field})=p$. 
Let $m>2$ be a prime number; if $p>0$, assume $m\neq p$.
Consider three distinct skew lines $L_1,L_2,L_3\subset\PP^3_{\overline{\field}}$ 
and two distinct transversals $T_1,T_2$.  
Let $H_1,\ldots,H_q$ be the $n_m=(m-1)(m-2)$ lines which can be used for $L_4$
in order that $G_1$ has order $m$.
Let $(l_i,t_i)\in{\overline{\field}}^2$ be the values of $l$ and $t$ corresponding to $H_i$.
Then there are $n_m$ distinct values of $l_i$ and $n_m$ distinct values of $t_i$.
\end{proposition}

\begin{proof} We use $\alpha$ and $\gamma$ as defined in Remark \ref{AlgorithmRem}. 
We consider the case of the $l$'s; the proof for the $t$'s works the same way.
Suppose two of the $l$'s are equal, say $l$ and $l'$.
As in the proof of Proposition \ref{LineCount}, we have
$l=\frac{\gamma\alpha-1}{\alpha\gamma-\alpha}$, with
$\gamma$ and $\alpha$ being $m$th primitive roots of 1 with $\alpha\gamma\neq1$, and
$l'=\frac{\gamma'\alpha'-1}{\alpha'\gamma'-\alpha'}$, with
$\gamma'$ and $\alpha'$ being $m$th primitive roots of 1
with $\alpha'\gamma'\neq1$.

So we can write $l=(u^{i+1}-1)/(u^i(u-1))$ and $l'=(v^{j+1}-1)/(v^j(v-1))$,
where $u=\gamma$, $\alpha=u^i$, $v=\gamma'$ and $\alpha'=v^j$, 
for some $i,j$ with $0<i< m-1$, $0<j< m-1$.
(We have $i< m-1$ since $u^{i+1}=\alpha\gamma\neq 1$,
and likewise for $j$.)
Assume $l=l'$; we need to show that $u=v$ and $i=j$.

Suppose first that $j>i$.
By the division algorithm we have 
$l=1+u^{-1}+\cdots+u^{-i}$ and $l'=1+v^{-1}+\cdots+v^{-j}$.
There is no harm in replacing $u$ by $1/u$ and $v$ by $1/v$, so that we have
$l=1+u+\cdots+u^i$ and $l'=1+v+\cdots+v^j$.
But $v=u^k$ for some $0<k< m-1$, so
$1+u+\cdots+u^i=1+v+\cdots+v^j$ is equivalent to
$1+u+\cdots+u^i=1+u^{[k]}+\cdots+u^{[kj]}$, where 
$[kj]$ means the exponent has been reduced modulo $m$.
Since $j>i$, after cancelling common terms
in $1+x^{[k]}+\cdots+x^{[kj]}-(1+x+\cdots+x^i)$
we will have terms left over, giving a polynomial in $x$
with no constant term and, after dividing out by $x$, degree at most $m-2$.
But this polynomial vanishes for $x=u$.
Since the minimal polynomial for $u$ is
$1+x+\cdots+x^{m-1}$, this is impossible.
Thus $j\le i$, and, by symmetry, $j\ge i$, so $i=j$.

Now we want to show $k=1$. Since $i=j$, we have $1+u^{[k]}+\cdots+u^{[ki]}-(1+u+\cdots+u^i)=0$.
Reasoning as before, this implies $\{u^{[k]},\ldots,u^{[ki]}\}=\{u,\ldots,u^i\}$. 
Assume $k\neq1$; in order for $\{u^{[k]},\ldots,u^{[ki]}\}\subseteq\{u,\ldots,u^i\}$,
we must have $k\geq m-i$, since $k<m-i$ means the gap from $i$ to $m$ will be
big enough that the least multiple $ks$ of $k$ with $ks>i$ 
has $1\leq s\leq i$ and $ks<m$ and hence $u^{[ks]}=u^{ks}$ is in $\{u^{[k]},\ldots,u^{[ki]}\}$ but could not be in $\{u,\ldots,u^i\}$.

On the other hand, $\{u^{[k]},\ldots,u^{[ki]}\}=\{u,\ldots,u^i\}$ also implies that
$\{u^{[k(i+1)]},\ldots,u^{[k(m-1)]}\}=\{u^{i+1},\ldots,u^{m-1}\}$
and thus that $\{u^{[-k(m-i-1)]},\ldots,u^{[-k]}\}=\{u^{-(m-i-1)},\ldots,u^{-1}\}$, so
similarly we must have $k\geq m-(m-i-1)$, or $k\geq i+1$.
Thus $2k\geq m+1$.

But $u^k\in \{u,\ldots,u^i\}$ implies $k\leq i$, while
$u^{-k}\in\{u^{-(m-i-1)},\ldots,u^{-1}\}$ implies $k\leq m-i-1$,
so $2k\leq m-1$. I.e., $m-1\geq 2k\geq m+1$; this contradiction implies $k=1$ so $v=u$.
\end{proof}

The following lemma is elementary but useful. It is Lemma 2.4 in 
\cite{POLITUS2};
we include the proof for the convenience of the reader. 

\begin{lemma}\label{LiftingLemma}
For any characteristic, given skew lines $T_1$ and $T_2$ and skew lines $T'_1, T'_2$ 
in $\PP^3$ and automorphisms $f_i\colon T_i\to T'_i$,
there is an automorphism  $f\colon \PP^3\to \PP^3$ restricting on $T_i$ to $f_i$.
\end{lemma}

\begin{proof}
Choose coordinates such that $T_1$ is $x=y=0$ and $T_2$ is $z=w=0$.
The first two columns of the matrix for $f$ are determined by $f_1$ and the second two columns
of the matrix for $f$ are determined by $f_2$. The matrix we get is invertible since
any two points on $T_1$ and any two points on $T_2$ give four points which are not contained in any
plane, and their images under $f$ are not contained in any plane.
\end{proof}

\begin{remark}\label{PosDimModuli}
Formerly it was unknown whether nongrid $(a,b)$-geproci 
sets can have positive dimensional moduli. Here we show that they can. 
For any 4 skew lines having an $[m,4]$-half grid with $m>2$, recall the group $G_1$ is finite, and so
over the complex numbers there are two transversals.
For any point $p$ on one of the skew lines but not on one of the two transversals,
the orbit $[p]$ is an $[r,4]$-half grid with $r$ being the order of $G_1$. Unions of $t$ different such orbits
are $[tr,4]$-half grids and the relative positions of the $t$ orbits give positive dimensional moduli.
This also shows that half grids $Z$ with the same matroid need not be projectively equivalent.
(The matroid here is the usual one given by the linear spans of the subsets of $Z$.)
To be specific, take $r\ge 2$ distinct orbits, each of which is a $[3,4]$-half grid on 4 lines
as in Example \ref{Example:D4}, to get uncountably many nonprojectively equivalent
$[3r, 4]$-half grids. 
Something similar occurs for grids.
An $(a,b)$-grid is a union of $a$ orbits,
where each orbit is a $(1,b)$-grid on the $b$ disjoint grid lines (we could also take $b$ orbits where each 
orbit is a $(1,a)$-grid on the $a$ other grid lines). The relative positions of these orbits gives rise to positive dimensional moduli for grids.
\qed\end{remark} 

When $m<b$, an $[m,b]$-half grid $Z\subset\PP^3_{\overline{\field}}$ 
can be a half grid with respect to more than one set of $b$ skew lines
(see the end of the proof of Theorem \ref{curve theorem}).
The next lemma addresses the case $m\geq b$.
It holds in every characteristic.

\begin{lemma}\label{projEqLemma}
Let $Z, Z'\subset\PP^3_{\overline{\field}}$ be $[m,b]$-half grids
with respect to $b$ skew lines $\mathcal L$ and $\mathcal L'$, respectively. 
\begin{enumerate}
\item[(a)] If $m\geq b$, then $\mathcal L$ is the unique set of $b$ skew lines
with respect to which $Z$ is an $[m,b]$-half grid.
\item[(b)] If $m\geq b$ and $Z$ and $Z'$ are projectively equivalent,
then so are $\mathcal L$ and $\mathcal L'$.
\item[(c)] If $m\geq b$, $G_{\mathcal L}$ and $G_{\mathcal L'}$ are abelian,
$Z$ and $Z'$ are single orbits and $\mathcal L$ and $\mathcal L'$ are 
projectively equivalent, then so are $Z$ and $Z'$.
\end{enumerate}
\end{lemma}

\begin{proof}
(a) Let $\mathcal L''$ be $b$ skew lines with respect to which $Z$ is an $[m,b]$-half grid.
Let $L\in\mathcal L''$. First assume $m>b$.
Then $L$ contains $m>b$ points of $Z$. Each point is on one of the
$b$ lines in $\mathcal L$, and hence at least 2 of the points are on the same line of
$\mathcal L$, so $L\in \mathcal L$, hence $\mathcal L=\mathcal L''$.
Now assume $m=b$. If $L$ is not in
$\mathcal L$, then each of the $m=b$ points of $L\cap Z$ is on a different line in 
$\mathcal L$, so $L$ is a transversal for $\mathcal L$. Since the lines in $\mathcal L''$
are skew, none of the lines in $\mathcal L''$ are in $\mathcal L$. Thus $Z$ is
an $(m,m)$-grid whose grid lines are $\mathcal L$ 
and $\mathcal L''$, but we assumed $Z$ was a half-grid and thus not a grid.

(c) If $Z$ and $Z'$ are projectively equivalent, then there is a linear automorphism
$\psi \colon \PP^3_{\overline{\field}}\to \PP^3_{\overline{\field}}$ with $\psi(Z)=Z'$.
Due to uniqueness of the lines they are half grids with respect to, we also have
$\psi(\mathcal L)=\mathcal L'$, so $\mathcal L$ and $\mathcal L'$ are projectively equivalent.

(d) If  $\mathcal L$ and $\mathcal L'$ are projectively equivalent, there is a 
linear automorphism $\psi \colon \PP^3_{\overline{\field}}\to \PP^3_{\overline{\field}}$ with 
$\psi(\mathcal L)=\mathcal L'$. But $Z'=[q']$ for some $q'\in Z'$ and $Z=[q]$ for some $q\in Z$,
so $\psi(Z)=[\psi(q)]$. By Theorem \ref{TransThm} and either
Proposition \ref{twoTrans}(a) (third bullet point) or Corollary \ref{oneTransCor}
(depending on whether $\mathcal L$ has two transversals or one of multiplicity 2),
$[\psi(q)]$ and $[q']$ are projectively equivalent, hence so are $Z$ and $Z'$.
\end{proof}

\begin{remark}\label{FinalRem}
Here we display a table of the number of equivalence classes
when $4\leq m\leq 20$ of single orbit complex $[m,4]$-half grids.
The cases of $m=3, 4$ were done in \cite{POLITUS, Fields};
the results of the table are purely computational. 
(The computations can be carried out for any characteristic
but so far we have done them only in characteristic 0, which allows us to assume that the four
half grid lines have two distinct transversals.) 
The computation works as follows.
If $Z$ is a single orbit $[m,4]$-half grid with $m\geq 4$, 
then it is so with respect to a unique set of 4 lines $\mathcal L$ by 
Lemma \ref{projEqLemma}(a). Four skew lines always have at least one transversal
(Remark \ref{r.2 transversal-a}),
but in characteristic 0, they support a half grid only when they have exactly
two transversals (Proposition \ref{oneTransProp}(c)), hence 
$G_{\mathcal L}$ is abelian (Theorem \ref{TransThm}).
Thus by Lemma \ref{projEqLemma}(b,c), counting projective equivalence classes of single orbit $[m,4]$-half grids 
is the same as counting projective equivalence classes of sets 
$\mathcal L$ of 4 skew lines with $|G_{\mathcal L}|=m$. 
Given four lines $\mathcal L=\{L_1,L_2,L_3,L_4\}$ with two transversals $T_1,T_2$, we 
can assume (up to projective equivalence) that the lines $L_1,L_2,L_3$ are in standard position
and the transversals $T_i$ are the standard ones. Using Remark \ref{AlgorithmRem}
we can list all lines $L_4$ such that $\mathcal L=\{L_1,L_2,L_3,L_4\}$ has 
$|G_{\mathcal L}|=m$. Then using the cross ratios of the points of intersection
of the $L_i$ with the $T_j$ we can determine the projective equivalence 
classes in the set $\Lambda$ of these finitely many sets $\mathcal L$ of 4 lines.
What we find computationally is that the projective equivalence classes in $\Lambda$  
partition $\Lambda$ into subsets each of cardinality either 6 or 12.
The notation we use to indicate the partition is
$6^a 12^b$ to indicate there are
$a$ sets of projective equivalence classes of size 6 and
$b$ of size 12, and hence $|\Lambda|=6a+12b$
(i.e., there are $6a+12b$ choices for $L_4$), giving $a+b$ 
projective equivalence classes.
(The projective equivalence classes of size 6 come from cases where the 4 lines
have an involutory projective automorphism which transposes the two transversals).

\begin{center}
\begin{tabular}{|r|r|r|r|}
\hline
$m$ & \# classes & $|\Lambda|$ & partition of $\Lambda$\\
\hline
4 & 1 & 6 & $6^1 12^0$\\
\hline
5 & 2 & 12 & $6^2 12^0$\\
\hline
6 & 2 & 18 & $6^1 12^1$\\
\hline
7 & 4 & 30 & $6^3 12^1$\\
\hline
8 & 4 & 36 & $6^2 12^2$\\
\hline
9 & 6 & 54 & $6^3 12^3$\\
\hline
10 & 6 & 60 & $6^2 12^4$\\
\hline
11 & 10 & 90 & $6^5 12^5$\\
\hline
12 & 8 & 84 & $6^2 12^6$\\
\hline
13 & 14 & 132 & $6^6 12^8$\\
\hline
14 & 12 & 126 & $6^3 12^9$\\
\hline
15 & 16 & 168 & $6^4 12^{12}$\\
\hline
16 & 16 & 168 & $6^4 12^{12}$\\
\hline
17 & 24 & 240 & $6^8 12^{16}$\\
\hline
18 & 18 & 198 & $6^3 12^{15}$\\
\hline
19 & 30 & 306 & $6^9 12^{21}$\\
\hline
20 & 24 & 264 & $6^4 12^{20}$\\
\hline
\end{tabular}
\end{center}
\qed\end{remark}

The proof of the next result follows the procedure described in
Remark \ref{FinalRem}. We assume the complex numbers for two reasons. 
One is that in characteristic 0, $[m,4]$-half grids always have two
distinct transversals (see Proposition \ref{oneTransProp}(c)). 
Without assuming characteristic 0 we would need to make
having two transversals a hypothesis, since
in positive characteristics there could be a single transversal of multiplicity 2,
a situation requiring a somewhat different analysis which we have not
yet carried out. In fact, most of the proof
of the theorem goes through in any characteristic for half grids on lines with two transversals.
But one step of the proof uses the geometry of the complex numbers specifically.
Replacing this step by a characteristic free argument looks messy and will require
additional study.

\begin{theorem}\label{PrimeThm}
If $m\geq 5$ is prime, then there are exactly $(m^2-1)/12$ projective equivalence classes of 
single orbit complex $[m,4]$-half grids.
\end{theorem}

\begin{proof} 
As explained in Remark \ref{FinalRem}, it is enough to count projective equivalence classes
of sets $\mathcal L=\{L_1,L_2,L_3,L_4\}$ of 4 lines with two transversals $T_1,T_2$
where $L_1,L_2,L_3$ are in standard position, the $T_i$ are standard and $|G_{\mathcal L}|=m$.
By Proposition \ref{LineCount}, there are $n_m$ choices for $L_4$, and since $m$ is prime we have
$n_m=(m-1)(m-2)$. The choices for $L_4$ are given in Remark \ref{AlgorithmRem}.
They are the lines $L_{t,l}$ through the points $(0:0:t:1)\in T_1$ and $(l:1:0:0)\in T_2$,
where $(t,l)=\left(\frac{\gamma-1}{\alpha\gamma-1},\frac{\alpha\gamma-1}{\alpha\gamma-\alpha}\right)$
come from all ordered pairs $(\alpha,\gamma)$ where $\alpha,\gamma$ are $m$th roots of unity which together generate 
the multiplicative subgroup of $\overline{\field}^*$ of order $m$ and which satisfy
$\alpha,\gamma,\alpha\gamma\neq1$.

Since $m$ is prime this means $\alpha$ and $\gamma$
are primitive and that there are $(m-1)(m-2)$ choices of the ordered pair $(\alpha,\gamma)$.
Let $\mathcal R$ be the set of these $(m-1)(m-2)$ ordered pairs. 
Let $\mathcal P$ be the corresponding set of all ordered pairs $(t,l)$.
Note that we can recover $\alpha$ and $\gamma$ from $t$ and $l$ (see Remark \ref{AlgorithmRem}),
so $|\mathcal R|=|\mathcal P|=(m-1)(m-2)$. Let
$\lambda$ be the set of the lines $L_{t,l}$ (so $|\lambda|=(m-1)(m-2)$ by Proposition \ref{LineCount}), 
and let $\Lambda$ be the set of all
sets $\Lambda_{t,l}=\{L_1,L_2,L_3,L_{t,l}\}$ of the four lines $L_1,L_2,L_3,L_{t,l}$
(so $|\Lambda|=(m-1)(m-2)$). 

Now consider some $\Lambda_{t,l}\in\Lambda$. 
The lines in $\Lambda_{t,l}$ meet $T_1$ at the points
$(0:0:0:1)$, $(0:0:1:1)$, $(0:0:1:0)$ and $(0:0:t:1)$ (which as points of $\PP^1$
we can think of as, respectively, $0,1,\infty,t$) and they meet $T_2$ at
$(1:0:0:0)$, $(1:1:0:0)$, $(0:1:0:0)$ and $(l:1:0:0)$
(i.e., at $0,1,\infty,1/l$). 

If $\Lambda_{t',l'}\in\Lambda$ and if $\Lambda_{t,l}$ and $\Lambda_{t',l'}$
are projectively equivalent, then there is an automorphism $\psi$ of $\PP^3$
(i.e., a projectivity) with $\psi(\Lambda_{t,l})=\Lambda_{t',l'}$.
Since $\psi$ is linear and $T_1,T_2$ are the only transversals for 
$\Lambda_{t,l}$ and $\Lambda_{t',l'}$, we have
$\psi(T_1\cup T_2)=T_1\cup T_2$. Hence $\psi$ either preserves each transversal
or swaps them. 

Suppose $\psi$ preserves each transversal. Then 
$\psi(\{0,1,\infty,t)\}=\{0,1,\infty,t')\}$ and $\psi(\{0,1,\infty,1/l)\}=\{0,1,\infty,1/l')\}$,
so some permutation of $0,1,\infty,t$ has the same cross ratio as does $0,1,\infty,t'$, and the same permutation of
$0,1,\infty,1/l$ has the same cross ratio as does $0,1,\infty,1/l'$. 

Conversely, if some permutation $\sigma$ of 
$0,1,\infty,t$ has the same cross ratio as does $0,1,\infty,t'$, and if the same permutation of
$0,1,\infty,1/l$ has the same cross ratio as does $0,1,\infty,1/l'$, then 
there are automorphisms $f_i\colon T_i\to T_i$ such that
$f_1(\sigma(0))=0$,
$f_1(\sigma(1))=1$,
$f_1(\sigma(\infty))=\infty$ and
$f_1(\sigma(t))=t'$, and
$f_2(\sigma(0))=0$,
$f_2(\sigma(1))=1$,
$f_2(\sigma(\infty))=\infty$ and
$f_2(\sigma(1/l))=1/l'$.
Then by Lemma \ref{LiftingLemma}
there is an automorphism $\Psi$ of $\PP^3$
which restricts to $f_i$ on $T_i$, $i=1,2$, and hence
$\Lambda_{t,l}$ and $\Lambda_{t',l'}$
are projectively equivalent.

Given $(t,l)\in\mathcal P$, one computes that the cross ratios of $0,1,\infty,t$ and of 
$0,1,\infty,1/l$ using the same permutation each time are:
\begin{align}\label{chiDisplay1}
\chi(t,l)=\left\{\left(t,\frac{1}{l}\right),
\left(\frac{1}{t},l\right), 
\left(1-t,\frac{l-1}{l}\right), 
\left(\frac{1}{1-t},\frac{l}{l-1}\right), 
\left(\frac{t-1}{t},1-l\right), 
\left(\frac{t}{t-1},\frac{1}{1-l}\right)\right\}.
\end{align}
Thus $\Lambda_{t,l}$ and $\Lambda_{t',l'}$ are projectively equivalent by a projectivity
preserving each $T_i$ if and only if $(t',1/l')\in \chi(t,l)$.

Now we consider the case of projectivities that swap the transversals.
For this we want to use the fact that if $(t,l)\in\mathcal P$, then $\chi(t,l)\subset \mathcal P$.

Regarding $t=\frac{\gamma-1}{\alpha\gamma-1}$ as defining a function
$f(\alpha,\gamma)=\frac{\gamma-1}{\alpha\gamma-1}$
and $1/l=\frac{\alpha\gamma-\alpha}{\alpha\gamma-1}$ as defining a function
$g(\alpha,\gamma)=\frac{\alpha\gamma-\alpha}{\alpha\gamma-1}$,
we can rewrite the cross ratios for $0,1,\infty,t$ and for $0,1,\infty,1/l$ as
\begin{align}\label{chiDisplay2}
t=f(\alpha,\gamma), \frac{1}{t}=f\left(\frac{1}{\alpha},\alpha\gamma\right), 1-t=f\left(\frac{1}{\gamma},\frac{1}{\alpha}\right), 
\frac{1}{1-t}=f\left(\gamma,\frac{1}{\alpha\gamma}\right),
\end{align}
\begin{align*}
\frac{t-1}{t}=f\left(\frac{1}{\alpha\gamma},\alpha\right), 
\frac{t}{t-1}=f\left(\alpha\gamma,\frac{1}{\gamma}\right)
\end{align*}
\begin{align*}
1/l=g(\alpha,\gamma), l=g\left(\frac{1}{\alpha},\alpha\gamma\right), \frac{l-1}{l}=g\left(\frac{1}{\gamma},\frac{1}{\alpha}\right), 
\frac{l}{l-1}=g\left(\gamma,\frac{1}{\alpha\gamma}\right),
\end{align*}
\begin{align*}
1-l=g\left(\frac{1}{\alpha\gamma},\alpha\right), 
\frac{1}{1-l}=g\left(\alpha\gamma,\frac{1}{\gamma}\right).
\end{align*}
Notice that each cross ratio is obtained as the value of $f$ (resp. $g$)
at an ordered pair $(\delta,\epsilon)$ where $\delta$ and $\epsilon$
are primitive $m$th roots of unity, neither equal to 1 and such that
$\delta\epsilon\neq 1$; i.e., each of these 6 cross ratios is itself one of the
allowed values for $t$ (and likewise for $l$) and thus $\chi(l,t)\subset \mathcal P$.
Moreover, the ordered pairs at which we evaluate $f$ to get the 6 cross ratios
for $t$ above (and likewise for $l$) are distinct.
(For example, if $(\alpha,\gamma)=\left(\frac{1}{\alpha},\alpha\gamma\right)$,
then $\alpha=\frac{1}{\alpha}$, hence $\alpha^2=1$, contrary to the assumption that
$\alpha$ is a primitive $m$th root of unity for $m\geq5$. Or, if you prefer,
$\gamma=\alpha\gamma$ hence $\alpha=1$.)
Now, by Proposition \ref{The t and l values are different},
the $t$'s (and the $l$'s) coming from the $(m-1)(m-2)$ allowed choices
of the ordered pairs $(\alpha,\gamma)$ are distinct.
Thus $|\chi(t,l)|=6$.

Now we claim $(t,l)\in\mathcal P$ if and only if $(1/l,1/t)\in\mathcal P$ and that $\Lambda_{t,l}$ is
projectively equivalent to $\Lambda_{1/l,1/t}$ by a projectivity that swaps the transversals.
(This is because the projective equivalence given by
$(a:b:c:d)\mapsto (d:c:b:a)$ preserves $L_1,L_2,L_3$, but swaps $T_1$ and $T_2$ 
and swaps $L_{t,l}$ with $L_{1/l,1/t}$.) Thus $\Lambda_{t,l}$ and $\Lambda_{t',l'}$ are projectively equivalent 
by a projectivity that swaps the transversals if and only if $(t',1/l')\in \chi(1/l,1/t)$.

Putting the two cases together we get that $\Lambda_{t,l}$ and $\Lambda_{t',l'}$ are projectively equivalent
if and only if $(t',1/l')\in \chi(t,l)\cup \chi(1/l,1/t)$. We will show that the subsets 
$\chi(t,l)$ partition $\mathcal P$, which implies 
there are two cases: either $\chi(t,l) = \chi(1/l,1/t)$ 
(in which case $| \chi(t,l)\cup \chi(1/l,1/t)|=6$)
or $\chi(t,l)\cap \chi(1/l,1/t)=\varnothing$
(in which case $| \chi(t,l)\cup \chi(1/l,1/t)|=12$).
We now want to count how often each case happens. What we will see is that
$\chi(t,l) = \chi(1/l,1/t)$ for $3(m-1)$ choices of $(t,l)$
and $|\chi(t,l)\cap \chi(1/l,1/t)|=0$ in the rest of the cases (so for $(m-1)(m-2)-3(m-1)=(m-1)(m-5)$
choices of $(t,l)$). Given this, we get that there are
$(m-1)(m-5)/12$ projective equivalence classes of cardinality 12 and $3(m-1)/6=(m-1)/2$
projective equivalence classes of cardinality 6, for a total of $(m-1)(m-5)/12+(m-1)/2=(m^2-1)/12$
equivalence classes.

All that's left is to confirm that the subsets $\chi(t,l)$ partition $\mathcal P$
and then to show $|\chi(t,l)\cup \chi(1/l,1/t)|=6$ in $3(m-1)$ cases and
$|\chi(t,l)\cup \chi(1/l,1/t)|=12$ in $(m-1)(m-5)$ cases.

If $\chi(t,l)\cap\chi(t',l')$ is not empty, then some ordered pair
$(t^*,l^*)$ is in both, hence $\chi(t^*,1/l^*)$ is in both, so $\chi(t,l)=\chi(t^*,1/l^*)=\chi(t',l')$.
Thus the sets $\chi(t,l)$ partition $\mathcal P$ into subsets of cardinality 6.

Now say $|\chi(t,l)\cup \chi(1/l,1/t)|=6$; i.e., $\chi(t,l) = \chi(1/l,1/t)$. 
This means $1/l=\frac{\alpha\gamma-\alpha}{\alpha\gamma-1}$
is either $t,1/t,1-t,1/(1-t),(t-1)/t$ or $t/(1-t)$. 

If $1/l=t$, then from Display \eqref{chiDisplay2}
we have $\frac{\alpha\gamma-\alpha}{\alpha\gamma-1}=\frac{\gamma-1}{\alpha\gamma-1}$
so $\alpha=1$ which is excluded by hypothesis.

If $1/l=1/t$, then 
we have $\frac{\alpha\gamma-\alpha}{\alpha\gamma-1}=\frac{\alpha\gamma-1}{\gamma-1}$
so $\alpha\gamma^2=1$. This happens for $(m-1)$ choices of $(\alpha,\gamma)$.

If $1/l=1-t$, then 
we have $\frac{\alpha\gamma-\alpha}{\alpha\gamma-1}=\frac{\alpha\gamma-\gamma}{\alpha\gamma-1}$
so $\alpha=\gamma$. This happens for $(m-1)$ choices of $(\alpha,\gamma)$.

If $1/l=1/(1-t)$, then 
we have $\frac{\alpha\gamma-\alpha}{\alpha\gamma-1}=\frac{\alpha\gamma-1}{\alpha\gamma-\gamma}$
which gives $\alpha+\gamma=3-\frac{1}{\alpha\gamma}$. Here's where we use the complex numbers:
$\alpha$ and $\gamma$ have norm 1. The set of complex numbers of the form
$\alpha+\gamma$ where $\alpha$ and $\gamma$ have norm 1 is a disk of radius 2 centered at 0, 
while the set of points of the form $3-\frac{1}{\alpha\gamma}$
is a circle centered at 3. These two sets have only the number 2 in common, which means
$3-\frac{1}{\alpha\gamma}=2$ hence $\alpha\gamma=1$ which is excluded by hypothesis.

If $1/l=(t-1)/t$, then 
we have $\frac{\alpha\gamma-\alpha}{\alpha\gamma-1}=\frac{\gamma-\alpha\gamma}{\gamma-1}$
which gives $2=\alpha\gamma+\frac{1}{\alpha}$. Arguing as in the previous case,
this happens only when $\alpha=\gamma=1$, which is excluded.

If $1/l=t/(t-1)$, then 
we have $\frac{\alpha\gamma-\alpha}{\alpha\gamma-1}=\frac{\gamma-1}{\gamma-\alpha\gamma}$
so $\alpha^2\gamma=1$. This happens for $(m-1)$ choices of $(\alpha,\gamma)$.

Thus the $3(m-1)$ ordered pairs $(\alpha,\gamma)$ where either $\alpha=\gamma$,
$\alpha^2\gamma=1$ or $\alpha\gamma^2=1$ are the ones with $|\chi(t,l)\cup \chi(1/l,1/t)|=6$,
the rest then have $|\chi(t,l)\cup \chi(1/l,1/t)|=12$.
\end{proof}

\paragraph*{\bf Acknowledgement.}
Our work was partly carried out during the Workshop on 
Lefschetz Properties in Algebra, Geometry, Topology and Combinatorics 
held at the Fields Institute in Toronto in the period May 15-19, 2023. 
We are grateful to the Institute for hosting our group and providing stimulating working conditions.

This work was completed while all authors participated in the 
Research Group ``Intersections in Projective Spaces'' hosted by the BIRS 
program in Kelowna in summer 2023. We thank BIRS for generous support 
and in particular Chad Davis for making our stay in British Columbia so pleasant.

{Chiantini and Favacchio are members of GNSAGA-INDAM.}
{Farnik was partially supported by National Science Centre, Poland, grant 2018/28/C/ST1/00339.
}
{Favacchio was partially supported by Fondo di
	Finanziamento
	per la Ricerca di
	Ateneo, Università degli studi di Palermo and by GNSAGA-INDAM.}
{Harbourne was partially supported by Simons Foundation grant \#524858.}
{Migliore was partially supported by Simons Foundation grant \#839618.
}
{Szemberg and Szpond were partially supported by National Science Centre, Poland, grant 2019/35/B/ST1/00723.}

We also thank Susan Hermiller and Carl Wang-Erickson for helpful consultations.

\section*{ Statements and Declarations.}

\noindent{\bf Competing interests:} The authors have no  potential conflicts 
of interest (financial or non-financial) to declare that are relevant to the content of this article.


\begin{thebibliography}{10}

\bibitem{Alperin}
R.~C. Alperin.
\newblock An elementary account of Selberg's lemma.
\newblock {\em Enseign. Math. (2)}, 33:269--273, 1987.

\bibitem{BGM}
A. Bigatti, A.~V. Geramita and J. Migliore.
\newblock Geometric consequences of extremal behavior in a theorem of Macaulay.
\newblock Trans. Amer. Math. Soc.346 (1994), 203--235.

\bibitem{BB}
R.~Bruck and R.~Bose. 
\newblock The construction of translation planes from projective spaces. 
\newblock Journal of Algebra 1 (1964), 85--102.

\bibitem{POLITUS}
L.~Chiantini, {\L }.~Farnik, G.~Favacchio, B.~Harbourne, J.~Migliore,
  T.~Szemberg, and J.~Szpond.
\newblock Configurations of points in projective space and their projections,
  arXiv:2209.04820.

\bibitem{POLITUS2}
L.~Chiantini, {\L }.~Farnik, G.~Favacchio, B.~Harbourne, J.~Migliore,
  T.~Szemberg, and J.~Szpond.
\newblock On the classification of certain geproci sets,  Lefschetz Properties. Current and New Directions, Springer, 2024, 81--96.

\bibitem{Fields}
L. Chiantini, \L. Farnik, P. De Poi, G. Favacchio, B. Harbourne, G. Ilardi, J. Migliore, T. Szemberg, J. Szpond,
\newblock Geproci sets on skew lines in ${\bf P}^3$ with two transversals,
\newblock {\em Journal of Pure and Applied Algebra}, 229 (2025), 107809, 
pp. 14.

\bibitem{CM}
L.~Chiantini and J.~Migliore.
\newblock Sets of points which project to complete intersections, and
  unexpected cones.
\newblock {\em Trans. Amer. Math. Soc.}, 374(4):2581--2607, 2021.
\newblock With an appendix by A. Bernardi, L.~Chiantini, G. Denham, G.
  Favacchio, B. Harbourne, J. Migliore, T. Szemberg and J. Szpond.

\bibitem{CDFPGR}
C.~Ciliberto, T.~Dedieu, F.~Flamini, R.~Pardini, C.~Galati, and S.~Rollenske.
\newblock Open problems.
\newblock {\em Boll. Unione Mat. Ital.}, 11(1):5--11, 2018.

\bibitem{Davis}
E.~D. Davis. 
\newblock Complete Intersections of Codimension 2 in ${\mathbb P}^r$: The Bezout-Jacobi-Segre Theorem Revisited.
\newblock Rend. Sem. Mat. Univers. Politecn. Torino, 43, 4 (1985), 333–353.

\bibitem{DGO}
E.~D. Davis, A.~V. Geramita and F. Orecchia. 
\newblock Gorenstein algebras and the Cayley-Bacharach theorem.
\newblock Proc. Amer. Math. Soc. 93 (1985), 593--597.

\bibitem{F}
X.~Faber.
\newblock Finite $p$-Irregular Subgroups of $PGL(2,k)$.
\newblock La Matematica 2, 479--522 (2023). 

\bibitem{Ganger}
A.~Ganger.
\newblock Spreads and transversals and their connection to geproci sets.
\newblock PhD thesis, University of Nebraska, 2024.

\bibitem{M2}
D.~R. Grayson and M.~E. Stillman.
\newblock Macaulay2, a software system for research in algebraic geometry.
\newblock Available at \url{http://www.math.uiuc.edu/Macaulay2/}.

\bibitem{hartley1993}
R.~I. Hartley.
\newblock Invariants of lines in space.
\newblock In {\em Proc. DARPA Image Understanding Workshop}, pages 737--744,
  1993.

\bibitem{hartley1994}
R.~I. Hartley.
\newblock Projective reconstruction and invariants from multiple images.
\newblock {\em IEEE Transactions on Pattern Analysis and Machine Intelligence},
  16(10):1036--1041, 1994.

\bibitem{huang1991}
R.~Q. Huang.
\newblock Invariants of sets of lines in projective 3-space.
\newblock {\em Journal of Algebra}, 143(1):208--218, 1991.

\bibitem{Kaplansky}
I.~Kaplansky.
\newblock {\em Fields and rings.}
\newblock Chicago, IL: University of Chicago Press, 2nd ed. edition, 1972.

\bibitem{kettinger2023}
J.~Kettinger.
\newblock The geproci property in positive characteristic, Proc. Amer. Math. Soc. 152 (2024), 3229--3242.

\bibitem{L}
J.~Lubin.
\newblock StackExchange, 2020, \url{https://math.stackexchange.com/users/17760/lubin}

\bibitem{LW}
G.~J. Leuschke and R.~Wiegand. 
\newblock  {\em Cohen-Macaulay representations}. 
Mathematical Surveys and Monographs, vol 181. 
American Mathematical Society, Providence, RI, 2012.

\bibitem{JuanBook}
J.~C. Migliore.
\newblock {\em Introduction to liaison theory and deficiency modules}, volume  165 of {\em Prog. Math.}
\newblock Boston, MA: Birkh{\"a}user, 1998.

\bibitem{PSS}
P.~Pokora, T.~Szemberg, and J.~Szpond.
\newblock Unexpected properties of the {K}lein configuration of $60$ points in
  ${\PP}^3$, Michigan Math. J. 74 (3), 599--615.

\bibitem{Polizzi}
F.~Polizzi and D.~Panov.
\newblock When is a general projection of $d^2$ points in ${\PP}^3$ a complete
  intersection?
\newblock {\em MathOverflow Question 67265}, 2011.

\bibitem{srinivasan2021}
P.~Srinivasan and K.~Wickelgren.
\newblock An arithmetic count of the lines meeting four lines in
  $\mathbb{P}^3$.
\newblock {\em Transactions of the American Mathematical Society},
  374(5):3427--3451, 2021.

\bibitem{WZ}
 P.~Wi\'sniewska and M.~Zi{\c e}ba, {\em Generic projections of the H4 configuration of points}, 
Adv. Applied Mathematics (2023), 142(4):102432.

\end{thebibliography}
\end{document}